\documentclass[review, 10pt]{elsarticle}

\usepackage{amssymb, amsmath, color}
\usepackage{graphicx, amsthm}
\usepackage{subfigure}
\usepackage{hyphenat}
\usepackage{latexsym, bm}
\usepackage{booktabs}
\usepackage{cases}
\usepackage{graphicx}
\usepackage{epstopdf}
\usepackage{amsthm}
\usepackage{fancyhdr}
\usepackage{float}
\usepackage{bm}
\usepackage{amssymb}
\usepackage{caption}
\usepackage{wrapfig}
\usepackage{verbatim}
\usepackage{graphicx}
\usepackage[colorlinks,
            linkcolor=blue,
            anchorcolor=blue,
            citecolor=blue]{hyperref}

\biboptions{numbers, sort&compress}
\allowdisplaybreaks[4]

\journal{}

\begin{document}
 \oddsidemargin = -5pt  \headsep=-10pt \footskip=50pt
\newtheorem{The}{Theorem}[section]
\newtheorem{lem}{Lemma}[section]

\renewcommand{\baselinestretch}{1.25}
\def\a{\alpha}
\def\b{\beta}
\def\g{\gamma}
\def\d{\delta}
\def\D{\Delta}
\def\l{\lambda}
\def\o{\omega}

\def\s{\zeta}
\def\sk{G}
\def\e{\eta}

\def\no{{\nonumber}}
\def\r{arrow}
\def\p{\partial}
\def\f{\frac}
\def\div{\nabla \cdot}
\def\u{{u}}
\def\w{{\omega}}
\def\v{{v}}
\def\i{{\varphi}}
\def\n{\mathcal{N}_1}
\def\bn{{\mathcal{N}}_2}
\def\k{\kappa}
\def\m{\mathcal{S}}
\def\T{\mathcal{T}}
\fancyhf{}

\begin{frontmatter}

\title{A Pressure Correction Projection Finite Element Method for The 2D/3D Time-Dependent Thermomicropolar Fluid Problem}

\author{Yuhang Ren}
\ead{Ryuhangs@163.com}
\author{Demin Liu \corref{mycorrespondingauthor}}
\cortext[mycorrespondingauthor]{Research Fund from the Key Laboratory of Xinjiang Province, Grant Number: 2022D04014; National Natural Science Foundation of China, Grant Number:  12061075.}
\ead{ Corresponding author: followtime@126.com}

\address{College of Mathematics and System Sciences,  Xinjiang University,  Urumqi 830046,  P.R. China}

\begin{abstract}
 In this paper, the pressure correction finite element method is proposed for the 2D/3D time-dependent thermomicropolar fluid equations. The first-order and second-order backward difference formulas (BDF) are adopted to approximate the time derivative term, stability analysis and error estimation of the first-order semi-discrete scheme are proved. Finally, some numerical examples are given to show the effectiveness and reliability of the proposed method, which can be used to simulate the problem with high Rayleigh number.
\end{abstract}

\end{frontmatter}

\section{Introduction}
The thermomicropolar fluid model can be used to describes the interaction behavior of  the polar fluid and the temperature under  the external source term, which can be looked as an extention of the micropolar fluid to include heat conduction and heat dissipation effects. The governing equations of the thermomicropolar fluid model express the conservation laws of mass, linear momentum, angular momentum and energy~\cite{1972Theory}. In the following paragraph, let $\Omega\subset\mathbb{R}^{d}$, $d$=2 or 3, be an open bounded domain with a sufficiently smooth boundary $\partial \Omega$. The time-dependent thermomicropolar equations can be written as in the primitive variables formulation~\cite{2018Rayleigh}
\begin{eqnarray}\label{1}
\left\{\begin{array}{ll}
\u_{t}-(\nu+\nu_r) \Delta \u+(\u \cdot \nabla) \u+\nabla p=2 {\nu_r}  \nabla \times \omega+J \hat{e}T+f_1,  \quad &\text { in }\Omega \times(0,  \T], \\
\div \u=0,  \quad &\text { in } \Omega \times(0,  \T],\\
\omega_{t}-\a\Delta \omega+(\u \cdot \nabla) \omega - \b\nabla\nabla\cdot\omega+4 {\nu_r}  \omega=2 {\nu_r} \nabla \times \u+f_2,  \quad &\text { in } \Omega \times(0,  \T],\\
T_{t}-\k\Delta T+(\u \cdot \nabla) T=D \nabla \times \omega \cdot \nabla T+f_3,  \quad &\text { in } \Omega \times(0,  \T],\\
\end{array}\right.
\end{eqnarray}
where the unknown variables are the velocity $u$, pressure $p$, angular velocity $\omega$ and temperature $T$. $J$ is the unit upward vector $(0, 1)^T$ when $d$=2 or $(0, 0, 1)^T$ when $d$=3. The meaning of the other positive constants can be explained as~\cite{1972Theory, jiang2018analysis}: $\nu$ is the kinematic Newtonian viscosity, $\nu_r$ is the kinematic microrotation viscosity, $\hat{e}$ is the Rayleigh number, the positive constants $\a$ and $\b$ are micropolar material viscosities, $\k$ is the thermal conductivity and $D$ is the micropolar thermal conduction. $f_1$,  $f_2$ and $f_3$ are the given source terms. For the simplicity of notation, we do not distinguish the difference between the results of $\nabla \times$ acting on $d=2$ or $d=3$, when space dimension $d$=3,
\begin{equation}
\begin{array}{ccc}
\u=(u_1,u_2,u_3), \quad \nabla \times \u=(\frac{\partial u_3}{\partial y}-\frac{\partial u_2}{\partial z},\frac{\partial u_1}{\partial z}-\frac{\partial u_3}{\partial x},\frac{\partial u_2}{\partial x}-\frac{\partial u_1}{\partial y}),\\
\omega=(\omega_1, \omega_2, \omega_3), \quad \nabla \times \omega =(\frac{\partial \omega_3}{\partial y}-\frac{\partial \omega_2}{\partial z},\frac{\partial \omega_1}{\partial z}-\frac{\partial \omega_3}{\partial x},\frac{\partial \omega_2}{\partial x}-\frac{\partial \omega_1}{\partial y}),\no\\
\end{array}
\end{equation}
or $d$=2,
\begin{equation}
\begin{array}{ccc}
\u=(u_1,u_2, 0),\quad\omega=(0, 0, \omega_3),\quad
\nabla \times \u=(0, 0, {\frac{\partial \u_2}{\partial x}-\frac{\partial \u_1}{\partial y}}), \quad
\nabla \times \omega =(\frac{\partial \omega_3}{\partial y}, -\frac{\partial \omega_3}{\partial x}, 0).\no
\end{array}
\end{equation}
In this paper, we supplied the following boundary condition and initial condition
\begin{eqnarray}\label{1_1}
\left\{\begin{array}{ll}
\u=0, \quad \omega=0, \quad  T=0,  \quad &\text { on } \partial \Omega \times(0,  \T], \\
\u=u^{0}, \quad  \omega=\omega^{0}, \quad T=T^{0},  &\text { in } \Omega \times\{0\}.
\end{array}\right.
\end{eqnarray}

For the above mentioned coupled nonlinear equations~(\ref{1}), the numerical difficulties can be listed as: (i) the restriction of incompressible conditions; (ii) the coupling of many physical variable; (iii) the non-linear property, etc. Obviously, the Galerkin variational problem of~(\ref{1}) is still a saddle-point problem. In order to avoid solving the saddle-point problem, a common strategy is to use the projection scheme for decoupling calculation~\cite{2017Pressure, 2021ha}.

In this paper, we mainly focus on the projection method based on the idea of pressure correction. The pressure correction method has being widely used to solve the incompressible Navier-Stokes equations~\cite{1992On, Shen1996On, Weinan1995Projection, 2001Accurate, 2020Error}. In the pressure correction  method, a series elliptical equations about the decoupled velocity and pressure only needed to be solved at each time step, so it will very effective for large-scale numerical simulations.
However, the non-physical Neumann boundary condition imposed on the pressure boundary by the standard pressure correction (SPC) method result in a numerical boundary layer, which limits the accuracy of the format~\cite{2006An}. Shen~\cite{2004On, 2011A} proposed a rotational pressure correction (RPC) method, the numerical boundary layer phenomenon can be effectively reduced.

The structure of this paper is as follows. In Section 2, we introduce some basic knowledge of finite element theory. In Section 3, the first-order and second-order semi-discrete pressure correction projection method is proposed, and the unconditional stability and error estimation of the first-order scheme are given. In Section 4, the corresponding fully discrete version is introduced, and the stability of the first-order fully discrete scheme is analyzed. In Section 5, some numerical experiments are carried out. In the last Section, we draw a conclusion.
\section{Preliminaries}
Let $\|\cdot\|_{m}$ denote the standard norm of Sobolev space $H^{m}(\Omega)$, $m \in \mathbb{Z}$. Specially, when $m=0$, $\|\cdot\|_{0}$ reduces to the standard norm of ${L^2}(\Omega)$, i.e.,
$$
\|u\|_{0}^{2}=\int_{\Omega}|u(x)|^{2} \mathrm{~d} x.
$$
In order to consider the variational problem of problem (\ref{1}), the following standard Sobolev spaces need to be introduced.
\begin{eqnarray}
\begin{aligned}
&{\rm X}=H_{0}^{1}(\Omega)^{d}=\{{v} \in H^{1}(\Omega)^{d}, {v}|_{\partial \Omega}=0\},\\
&{\rm Q}= L^{2}_0(\Omega)=\{q \in L^{2}(\Omega): \displaystyle{\int_{\Omega} q \mathrm{~d} x=0\}},\\
&{\rm M}= H_{0}^{1}(\Omega)=\{{G} \in H^{1}(\Omega), {G}|_{\partial \Omega}=0\},\\
&{\rm H}=\{{v} \in L^{2}(\Omega)^{d}, \nabla \cdot {v}=0, {v} \cdot {n}|_{\partial \Omega}=0\}.
\end{aligned}
\end{eqnarray}
Next recalling the following orthogonality identity,
\begin{align}
&\|\nabla v\|_{0}^{2}=\|\nabla \cdot v\|_{0}^{2}+\|\nabla \times v\|_{0}^{2},  \quad \forall v \in {\rm X},\\ \no
&\|v\|_{0} \leq c\|\nabla v\|_{0},  \quad \forall v \in {\rm X} \text{ or }{\rm M},\no
\end{align}
where $c$ represents a general constant, which value can depend on $\Omega, \nu, \nu_r, \a, \b, \k, D, u^0, \omega^0, T^0, f_1$, $f_2$ and $f_3$ and may be different from the previous value.

The weak formulation of (\ref{1}) reads: find $(u, p, \omega, T) \in {\rm X} \times {\rm Q} \times {\rm X} \times M$ for all $t \in (0,\T]$, such
that, for all $(\v, q, \s, \sk) \in {\rm X} \times {\rm Q} \times {\rm X} \times {\rm M}$,
\begin{equation}\label{weak}
\left\{\begin{array}{l}
(\u_{t},  {\v})+B((\u,  p) ;({\v},  q))+\n(\u ; \u,  {\v})=2\nu_rr(\u, {\v})+ \hat{e}(JT, {\v})+(f_1, \v),\\
(\omega_{t}, \s)+\a a(\omega, \s)+\n(\u;T, \s)+4\nu_r(\omega, \s)+\b(\nabla\cdot\omega,\nabla\cdot \s)=2\nu_r r(\u, \s)+(f_2, \s),\\
(T_{t}, \sk)+\k\bar{a}(T, \sk)+\bn(u;T, \sk)=D\m(u ; T, \sk)+(f_3, \sk),
\end{array}\right.
\end{equation}
where the bilinear and trilinear forms are defined as
\begin{equation}
\begin{aligned}\label{3line}
&a(\u,  {\v})=(\nabla \u,  \nabla {\v}), \quad \bar{a}(T,  \sk)=(\nabla T,  \nabla \sk), \quad d({v},  p)=(p,  \nabla \cdot \v), \\
&B((\u,  p) ;({\v},  q))=(\nu+\nu_r)a(\u,  {\v})-d({\v},  p)+d(\u,  q),\\
&\n(\u ; {\v},  {w})=((\u \cdot \nabla) {\v},  {w}),  \quad \bn(\u ; T,  \sk)=((u \cdot \nabla) T,  \sk), \\
&r(\u, \v)=(\nabla \times\u, \v),  \quad\bar{r}(u, v)=(\nabla \times\u, v),  \quad \m(\omega ; T, \sk)=(\nabla \times \omega \cdot \nabla T,  \sk).
\end{aligned}
\end{equation}
It is easy to show that $\mathcal{N}_1(\cdot ;\cdot,  \cdot)$ and $\mathcal{N}_2(\cdot ;\cdot,  \cdot)$ possess following properties~\cite{1984Navier}.
\begin{equation}\label{3line}
\begin{aligned}
\begin{array}{l}
\n(\u ; {\v},  {w})=-\n(\u ; {w},  {\v}),\quad \n(\u ; {\v},  {\v})=0,\quad \forall \u \in {\rm H} ,\quad \forall {\v}, {w} \in {\rm X},\\
\bn(\u ; T, \sk)=-\bn(\u ; \sk,  T, ),\quad \bn(\u ; \sk,  \sk)=0,\quad \forall \u \in{\rm H},\quad \forall T,  \sk \in{\rm  M}.
\end{array}
\end{aligned}
\end{equation}
\begin{equation}
\begin{aligned}\label{3line2}
| \n(\u ; w, \v)| \leq\left\{\begin{array}{l}
\|\u\|_{1}\|\v\|_{1}\|w\|_{1}, \quad \forall \u, w, \v \in {\rm X}, \\
\|\u\|_{2}\|\v\|_{0}\|w\|_{1}, \quad \forall \u \in {\rm H}^{2}(\Omega) \cap {\rm X}, \quad\forall w, \v \in {\rm X},\\
\|\u\|_{2}\|\v\|_{1}\|w\|_{0}, \quad \forall \u \in {\rm H}^{2}(\Omega) \cap {\rm X}, \quad\forall w, \v \in {\rm X},\\
\|\u\|_{1}\|\v\|_{2}\|w\|_{0}, \quad \forall \v \in {\rm H}^{2}(\Omega) \cap {\rm X}, \quad\forall \u, w \in {\rm X},
\end{array}\right.
\end{aligned}
\end{equation}
further, we use the skew-symmetric convection term \cite{2009Convergence}
$$
\n(\u ; {\v},  {w})=\n(\u ; {\v},  {w})+\frac{1}{2}((\nabla\cdot\u) \v, w)=\frac{1}{2} \n(\u ; {\v},  {w})-\frac{1}{2}\n(\u ; {w},  {\v})\quad \forall {\u}, {\v}, {w} \in {\rm X}.
$$

Similarly, we can prove the following Lemma.

\textbf{Lemma 2.1} $\m$ is a trilinear continuous form on $ {\rm X} \times {\rm M} \times {\rm M}$.
\begin{align}
&\m({\omega}; T, \sk)\leq \|\nabla {\omega}\|_{0}\|\nabla T\|_{1}\|\sk\|_{1},\label{m norm}\\
&\m(\omega; T, T)=0,\label{3line31}\\
&\m(\omega; T, \sk)=-\m(\omega; \sk, T),\label{3line3}
\end{align}
\textbf{Proof.} By using H\"older inequality and Sobolev imbedding theorems, we can be deduced
\begin{equation}
\begin{aligned}
\m(\omega; T, \sk)&=\sum_{i=1}^{n} \int_{\Omega} \left(\nabla \times\omega\right)_{i}\left(\partial_{i} T \right) \sk \mathrm{~d} x\leq\|\nabla \times {\omega}\|_{L^2}\|\nabla T\|_{L^4}\|\sk\|_{L^4}\\
&\leq\|\nabla {\omega}\|_{L^2}\|\nabla T\|_{L^4}\|\sk\|_{L^4}\leq c\|\nabla {\omega}\|_{0}\|\nabla T\|_{1}\|\sk\|_{1}. \nonumber \\
\end{aligned}
\end{equation}
The above inequality is due to $H^{1}(\Omega)^{d} \hookrightarrow L^{4}(\Omega)^{d}$. By using $\operatorname{div} (\nabla \times a)=0$, $a \in \mathbb{R}^3$, property (\ref{3line3}) is a consequence of (\ref{3line31}) if the variable $T$ is replaced by $T+\sk$, next by using the integration by parts and Gauss' law, then there exist for any $(\omega, T, \sk)$ $\in$ $({\rm X} \times {\rm M} \times {\rm M})$
\begin{equation}
\begin{aligned}
\m(\omega; T, T)&=\sum_{i=1}^{n} \int_{\Omega} \left(\nabla \times\omega\right)_{i}\left(\partial_{i} T\right) T \mathrm{~d} x=\sum_{i=1}^{n} \int_{\Omega} \left(\nabla \times\omega\right)_{i} \partial_{i} \frac{T^{2}}{2} \mathrm{~d}x\\
&=-\frac{1}{2} \sum_{i=1}^{n} \int_{\Omega} \partial_{i} \left(\nabla \times\omega\right)_{i}T^{2} \mathrm{~d} x =-\frac{1}{2}\int_{\Omega} \operatorname{div} (\nabla \times\omega)T^{2} \mathrm{~d} x=0. \nonumber \\
\end{aligned}
\end{equation}

Next the following Gronwall's Lemma is very important in the numerical analysis~\cite{2019Stability}.

\textbf{Lemma 2.2} (Discrete Gronwall's Lemma) Let $y^{n}, h^{n}, g^{n}$, and $f^{n}$ be nonnegative sequence such that:
$$
y^{m}+\tau \sum_{n=0}^{m} h^{n} \leqslant B+\tau \sum_{n=0}^{m}\left(g^{n} y^{n}+f^{n}\right), \quad \tau \sum_{n=0}^{N} g^{n} \leqslant K, \quad 0 \leqslant m \leqslant N=\left[\frac{\T}{\tau}\right],
$$
where $B$, $K$, is a given constant. In addition, assume that $g^{n} \tau<1$ for every $n$. Define $\rho=\max \limits_{0 \leqslant n \leqslant N}\left(1-g^{n} \tau\right)^{-1}$. Then:
$$
y^{m}+\tau \sum_{n=0}^{m} h^{n} \leqslant e^{\rho K}\left(B+\tau \sum_{n=0}^{m} f^{n}\right), \quad 0 \leqslant m \leqslant N.
$$

Next, we introduce the Stokes operator \cite{1992On}.

\textbf{Assumption 1} Given $\u \in {\rm H}$, $v = A^{-1} \u$ is the solution of the following Stokes equations:
$$
-\Delta v+\nabla q=\u, \quad \nabla \cdot v=0 \text { in } \Omega, \quad v|_{\partial \Omega}=0,
$$
we assume
\begin{equation}
\exists c_{1},  c_{2}>0, \text { such that } \forall \u \in {\rm H}: \left\{\begin{array}{l}
\|A^{-1} \u\|_{s} \leq c_{1}\|\u\|_{s-2}, \quad \text { for } s=1, 2, \\
c_{2}\|\u\|_{-1}^{2} \leq(A^{-1} \u,  \u) \leq c_{1}^{2}\|\u\|_{-1}^{2}.
\end{array}\right.
\end{equation}
Hence, we can use $(A^{-1}\u,\u)^{\frac{1}{2}}$ as an equivalent norm of ${\rm H}^{-1}(\Omega)^d$ for $\u \in {\rm H}$.

To simplify our presentation, we will assume that the initial conditions $u^0, \omega^0, T^0$ and right hand side are sufficiently smooth. More precisely, we assume~\cite{2021Global}

\textbf{Assumption 2}
\begin{equation}
\begin{array}{l}
\u^{0} \in {\rm V} \cap(H^{2}(\Omega))^{d},\quad \omega^{0} \in {\rm V} \cap H^{2}(\Omega)^{d},\quad  T^{0} \in {\rm V} \cap H^{2}(\Omega),\\
f_1, f_2 \in L^{\infty}(0,  \T ;(L^{2}(\Omega)^{d}) \cap L^{2}(0,  \T ; H^{1}(\Omega)^{d}),\quad f_3 \in L^{\infty}(0,  \T ;(L^{2}(\Omega)) \cap L^{2}(0,  \T ; H^{1}(\Omega)),
\end{array}
\end{equation}
and
\begin{equation}
\left\{\begin{array}{l}
\mathop{\sup}\limits_{t \in[0,  \T]}\{\|\u(t)\|_{2}+\|\omega(t)\|_{2}+\|T(t)\|_{2}+\|\u_{t}(t)\|_{0}+\|\omega_{t}(t)\|_{0}+\|T_{t}(t)\|_{0}+\|\nabla p(t)\|_{0}\} \leq C, \\[10pt]
\displaystyle\int_{0}^{\T}\|\nabla \u_{t}(t)\|_{0}^{2} d t+\int_{0}^{\T}\|\nabla \omega_{t}(t)\|_{0}^{2} d t+\int_{0}^{\T}\|\nabla T_{t}(t)\|_{0}^{2} d t \leq C, \\[10pt]
\displaystyle\int_{0}^{\T}\|\nabla \u_{t t}(t)\|_{-1}^{2} d t+\int_{0}^{\T}\|\nabla \omega_{t t}(t)\|_{-1}^{2} d t +\int_{0}^{\T}\|\nabla T_{t t}(t)\|_{-1}^{2} d t \leq C.
\end{array}\right.
\end{equation}
Among them, $C$ is a general normal number, and its value can be determined by $\Omega$, $\nu$, $\nu_r$, $\a$, $\b$, $\k$, $D$, $u$, $p$, $\omega$, $T$, $f_1$, $f_2 $, $f_3$, and the value is different from the last time it appeared.
\section{Pressure correction schemes: temporal semi-discrete schemes}
In this section, we consider the time discretization of the problem (\ref{1}) using the BDF1 and BDF2 to approximate the time derivative term, respectively. Let $\tau>0$ be the time step, $t_n=n\tau>0$, for all $0\leq n \leq N$, the time discrete approximation of $(u(t_n), p(t_n), \omega(t_n), T(t_n))$ will be expressed as $(u^n, p^n, \omega^n, T^n )$.

Pressure correction schemes are time-marching techniques composed of two substeps at each time step~\cite{2004On, 1992On}, the pressure is treated explicitly in the first substep and is corrected in the second one by projecting the velocity onto the space ${\rm H}$. We use the first-order SPC and RPC scheme and the second-order RPC scheme to solve the thermal micropolar fluid equations.
\subsection{First-order pressure correction schemes}
Adopting semi-implicit scheme for the nonlinear term and BDF1 to march in time, without causing confusion $n$ is the outward unit normal vector, the first-order pressure correction schemes are as follows:\\
\textbf{Algorithm 3.1(First-order SPC scheme)}\\
Assume that at each time step, $(\u^n, p^n, \omega^n, T^n)$ are given and one seeks $(\u^{n+1}, p^{n+1}, \omega^{n+1}, T^{n+1})$\\
Step1: Find $\tilde{\u}^{n+1}$ as the solution of
\begin{eqnarray}\label{3.1.1}
\left\{\begin{array}{lll}
\frac{1}{\tau}({\tilde{\u}^{n+1}-{\u}^{n}})-(\nu+{\nu_r}) \Delta \tilde{\u}^{n+1}+(\u^n \cdot \nabla)\tilde{\u}^{n+1}+\nabla p^n
=2\nu_r\nabla \times\omega^n+J \hat{e}T^n+f_1^{n+1},\\
\tilde {\u}^{n+1}|_{\partial\Omega}=0.
\end{array}\right.
\end{eqnarray}
Step2: Compute ${\u}^{n+1}$, $p^{n+1}$ as the solution of
\begin{eqnarray}\label{3.1.2}
\left\{\begin{array}{lll}
\frac{1}{\tau}({\u}^{n+1}-\tilde{\u}^{n+1})+\nabla(p^{n+1}-p^n)=0, \\
\nabla \cdot \u^{n+1}=0,\\
\u^{n+1} \cdot n|_{\partial \Omega}=0.
\end{array}\right.
\end{eqnarray}
Step3: Compute $\omega^{n+1}$ as the solution of
\begin{eqnarray}\label{3.1.3}
\left\{\begin{array}{lll}
\frac{1}{\tau}({\omega}^{n+1}-{\omega}^{n}) -\a\Delta \omega^{n+1} +(\u^{n+1} \cdot \nabla)\omega^{n+1} +4\nu_r\omega^{n+1}-\b\nabla\nabla\cdot\omega^{n+1}=2\nu_r\nabla \times\u^{n+1}+f_2^{n+1},\\
\omega^{n+1}|_{\partial\Omega}=0.
\end{array}\right.
\end{eqnarray}
Step4: Compute $T^{n+1}$ as the solution of
\begin{eqnarray}\label{3.1.4}
\left\{\begin{array}{l}
\frac{1}{\tau}(T^{n+1}-T^{n})-\k\Delta T^{n+1}+(u^{n+1} \cdot \nabla) T^{n+1}-D\nabla \times\omega^{n+1} \cdot \nabla T^{n+1}=f_3^{n+1}, \\
T^{n+1}|_{\partial\Omega}=0.
\end{array}\right.
\end{eqnarray}
\textbf{Remark 3.1.} In the actual calculation, we apply $\nabla \cdot$ to both sides of (\ref{3.1.2}), and then set $$\phi^{n+1}=p^{n+1}-p^n,$$ we can get\\
Step2.1: Find $\phi^{n+1}$ as the solution of
\begin{eqnarray}\label{3.1.2.1}
\left\{\begin{array}{l}
\Delta \phi^{n+1}=\frac{1}{\tau} \nabla \cdot \tilde {\u}^{n+1}, \\
\frac{\partial \phi}{\partial n}|_{\partial \Omega}=0.
\end{array}\right.
\end{eqnarray}
Step2.2: Update $\u^{n+1}, p^{n+1}$
\begin{eqnarray}\label{3.1.2.2}
\left\{\begin{array}{l}
\u^{n+1}=\tilde{\u}^{n+1}-\tau \nabla \phi^{n+1}, \\
p^{n+1}=\phi^{n+1}+p^{n}.
\end{array}\right.
\end{eqnarray}
\textbf{Algorithm 3.2(First-order RPC scheme)}
Just replace the Step2 in Algorithm 3.1 with\\
Step2(RPC): Compute ${\u}^{n+1}$, $p^{n+1}$ as the solution of
\begin{eqnarray}\label{Rotational 3.2}
\left\{\begin{array}{lll}
\frac{1}{\tau}({\u}^{n+1}-\tilde{\u}^{n+1})+\nabla(p^{n+1}-p^{n}+(\nu+\nu_r)\nabla \cdot \tilde{\u}^{n+1})=0, \\
\nabla \cdot \u^{n+1}=0,\\
\u^{n+1} \cdot n|_{\partial \Omega}=0.
\end{array}\right.
\end{eqnarray}
\textbf{Remark 3.2.} Actually, from (\ref{3.1.2}) we observe that
\begin{align}\label{remark3.3.1}
\left.\nabla\left(p^{k+1}-p^{k}\right) \cdot n\right|_{\partial \Omega}=0,
\end{align}
which implies that
\begin{align}\label{remark3.3.2}
\left.\nabla p^{k+1} \cdot n\right|_{\partial \Omega}=\left.\nabla p^{k} \cdot n\right|_{\partial \Omega}=\left.\cdots \nabla p^{0} \cdot n\right|_{\partial \Omega}.
\end{align}
As mentioned before, the SPC method imposes an additional Neumann boundary condition on the pressure, and this non-physical boundary condition results in an artificial numerical boundary layer phenomenon.\\
For RPC method, sum of (\ref{3.1.1}) and (\ref{Rotational 3.2}), using $\nabla \times \nabla \times \tilde{\u} ^{ k+1}=\nabla \times \nabla \times \u^{k+1}$, we have
\begin{equation}\label{3.2.2.2}
\left\{\begin{array}{l}
\frac{1}{\tau}(\u^{n+1}-\u^{n})+(\nu+\nu_r)\nabla \times \nabla \times \u^{n+1}+\nabla p^{n+1}=2\nu_r\nabla \times\omega^n+J \hat{e}T^n+f_1^{n+1}, \\
\nabla \cdot \u^{n+1}=0,\\
\u^{n+1} \cdot n|_{\partial \Omega}=0.
\end{array}\right.
\end{equation}
We observe from (\ref{3.2.2.2}) that
\begin{align}
\left.\partial_{n} p^{k+1}\right|_{\partial \Omega}=\left.\left(2\nu_r\nabla \times\omega^n+J \hat{e}T^n+f\left(t^{k+1}\right)-(\nu+\nu_r) \nabla \times \nabla \times u^{k+1}\right) \cdot n\right|_{\partial \Omega},
\end{align}
which, unlike (\ref{remark3.3.2}), is a consistent pressure boundary condition. In view of (\ref{3.2.2.2}), where the operator $\nabla\times\nabla\times$ plays a key role, the Algorithm 3.2 is called the rotational pressure correction scheme~\cite{2006An}.\\
\textbf{Remark 3.3.} Similarly, we apply $\nabla \cdot $ to both sides of (\ref{Rotational 3.2}), and then set $$\phi^{n+1}=p^{n+1}-p^n+(\nu+{\nu_r}) \nabla \cdot \tilde{\u}^{n+1},$$ we can get\\
Step2.1: Find $\phi^{n+1}$ as the solution of
\begin{eqnarray}\label{3.1.2.1}
\left\{\begin{array}{l}
\Delta \phi^{n+1}=\frac{1}{\tau} \nabla \cdot \tilde {\u}^{n+1}, \\
\frac{\partial \phi}{\partial n}|_{\partial \Omega}=0.
\end{array}\right.
\end{eqnarray}
Step2.2: Update $\u^{n+1}, p^{n+1}$
\begin{eqnarray}\label{3.1.2.2}
\left\{\begin{array}{l}
\u^{n+1}=\tilde{\u}^{n+1}-\tau \nabla \phi^{n+1}, \\
p^{n+1}=\phi^{n+1}+p^{n}-(\nu+{\nu_r})\nabla \cdot \tilde{\u}^{n+1}.
\end{array}\right.
\end{eqnarray}
\textbf{Theorem3.1(Stability)} Algorithm 3.1 is unconditionally stable under the Assumptions 1-3. If ($\u^{n+1}$, $p^{n+1}$, $\omega^{n+1}$, $T^{n+1}$) is the solution of Algorithm 3.1, then for all $0\leq n\leq N-1$ has a constant $c$ such that
\begin{align*}\label{3.7}
&\|\u^{N}\|_{0}^{2}+\tau^2\|\nabla p^{N}\|_{0}^{2}+(1+4\nu_r\tau)\|\omega^{N}\|_{0}^{2}+\|T^{N}\|_{0}^{2}
+\nu \tau \sum_{n=0}^{N-1}\|\nabla \tilde{\u}^{n+1}\|_{0}^{2}+ \a\tau \sum_{n=0}^{N-1}\|\nabla T^{n+1}\|_{0}^{2}\\
+&\k\tau\sum_{n=0}^{N-1}\|\nabla \omega^{n+1}\|_{0}^{2}
+2\b\tau\sum_{n=0}^{N-1}\|\nabla\cdot\omega^{n+1}\|_{0}^{2}
+\sum_{n=0}^{N-1}(\|\tilde{\u}^{n+1}-u^{n}\|_{0}^{2}+\|\omega^{n+1}-\omega^{n}\|_{0}^{2}+\|T^{n+1}-T^{n}\|_{0}^{2})\\
\leq& \|\u^{0}\|_{0}^{2}+\tau^2\|\nabla p^{0}\|_{0}^{2}+(1+4\nu_r\tau)\|\omega^{0}\|_{0}^{2}+\|T^{0}\|_{0}^{2}+c\tau\sum_{n=0}^{N-1}(\|f_1(t_{n+1})\|_{0}^{2}+\|f_2(t_{n+1})\|_{0}^{2}+\|f_3(t_{n+1})\|_{0}^{2}).
\end{align*}
\textbf{Proof.}
For any sequence $\{A^{n}\}_{n=0}^{N}$, then we have
\begin{equation}\label{Lemma 2N}
2(A^{n+1}-A^{n},  A^{n+1})=\|A^{n+1}\|_{0}^{2}-\|A^{n}\|_{0}^{2}+\|A^{n+1}-A^{n}\|_{0}^{2}.
\end{equation}
Taking the inner product of (\ref{3.1.1}), (\ref{3.1.3}), (\ref{3.1.4}) with $(2\tau\tilde{\u}^{n+1}, 2\tau \omega^{n+1}, 2\tau T^{n+1})$, using (\ref{Lemma 2N}), we have
\begin{align}
&\|\tilde{\u}^{n+1}\|_{0}^{2}-\|\u^{n}\|_{0}^{2}+\|\tilde{u}^{n+1}-u^{n}\|_{0}^{2}+2 (\nu+{\nu_r})\tau \|\nabla \tilde{\u}^{n+1}\|_{0}^{2}+2 \tau\left(\nabla p^{n},  \tilde{\u}^{n+1}\right) \no\\
&=4{\nu_r}\tau\left(\nabla \times \omega^{n},  \tilde{\u}^{n+1}\right)+2\hat{e}\tau \left(J T^{n},  \tilde{\u}^{n+1}\right)+2 \tau\left(f_1\left(t_{n+1}\right),  \tilde{\u}^{n+1}\right), \label{20}\\
&\|\omega^{n+1}\|_{0}^{2}-\|\omega^{n}\|_{0}^{2}+\|\omega^{n+1}-\omega^{n}\|_{0}^{2}+2 \a\tau\|\nabla \omega^{n+1}\|_{0}^{2}+8{\nu_r}\tau\|\omega^{n+1}\|_{0}^{2}+2\b\tau\|\nabla\cdot\omega^{n+1}\|_{0}^{2}\no \\
&=4{\nu_r}  \tau\left(\nabla \times \u^{n+1},  \omega^{n+1}\right)+2 \tau\left(f_2\left(t_{n+1}\right),  \omega^{n+1}\right), \label{21}\\
&\|T^{n+1}\|_{0}^{2}-\|T^{n}\|_{0}^{2}+\|T^{n+1}-T^{n}\|_{0}^{2}+2 \k\tau\|\nabla T^{n+1}\|_{0}^{2}=2 \tau\left(f_3\left(t_{n+1}\right),  T^{n+1}\right).\label{22}
\end{align}
According to the Cauchy-Schwarz inequality and Young's inequality, the right side of (\ref{20})-(\ref{22}) can be estimated
\begin{align}
4{\nu_r}\tau(\nabla \times \w^{n},  \tilde{\u}^{n+1})
&=4{\nu_r}\tau(\w^{n},  \nabla \times \tilde{\u}^{n+1}) \no \\
& \leq 4{\nu_r}\tau\|\w^{n}\|_{0}\|\nabla \times\tilde{\u}^{n+1}\|_{0} \no \\
& \leq {\nu_r}\tau\|\nabla \tilde{\u}^{n+1}\|_{0}^{2}+4{\nu_r}\tau\|\w^{n}\|_{0}^{2}, \label{23}\\
2\tau  \hat{e}(J T^{n},  \tilde{\u}^{n+1})
& \leq 2\tau\hat{e}\|T^{n}\|_{0}\|\nabla \tilde{\u}^{n+1}\|_{0} \no \\
& \leq \frac{\tau} {2} \nu \|\nabla \tilde{\u}^{n+1}\|_{0}^{2}+c \tau\| T^{n}\|_{0}^{2},\label{24}\\
2 \tau (f_1(t_{n+1}),  \tilde{\u}^{n+1})
& \leq 2\tau \|f_1(t_{n+1})\|_{0}\|\nabla\tilde{\u}^{n+1}\|_{0} \no \\
& \leq \frac{\tau} {2} \nu \|\nabla\tilde{\u}^{n+1}\|_{0}^{2}+c \tau \|f_1(t_{n+1})\|_{0}^{2}, \label{25}\\
4{\nu_r}\tau (\nabla \times{\u}^{n+1},  \w^{n+1})
& \leq 4{\nu_r}\tau\|\nabla \times {\u}^{n+1}\|_{0}\|\w^{n+1}\|_{0} \no \\
& \leq {\nu_r}\tau \|\nabla {\u}^{n+1}\|_{0}^{2}+4{\nu_r}\tau \|\w^{n+1}\|_{0}^{2}, \label{26}\\
2 \tau (f_2(t_{n+1}),  {\omega}^{n+1})
& \leq 2\tau \|f_2(t_{n+1})\|_{0}\|\nabla{\omega}^{n+1}\|_{0} \no \\
& \leq \a \|\nabla{\omega}^{n+1}\|_{0}^{2}+c \tau \|f_2(t_{n+1})\|_{0}^{2}, \label{27} \\
2 \tau (f_3(t_{n+1}),  T^{n+1})
& \leq 2\tau \|f_3(t_{n+1})\|_{0}\|\nabla T^{n+1}\|_{0} \no \\
& \leq \k \|\nabla T^{n+1}\|_{0}^{2}+c \tau \|f_3(t_{n+1})\|_{0}^{2}.\label{28}
\end{align}
Rewrite (\ref{3.1.2}) as $\u^{n+1}+\tau \nabla p^{n+1}=\tilde{\u}^{n+ 1}+\tau \nabla p^{n}$,  take the inner product of (\ref{3.1.2}) both sides according to $(\nabla p^{n+1},  \u^{n+1})=-(p^{n+1},  \nabla \cdot \u^{n+ 1})=0$. We get
\begin{eqnarray}\label{3.5}
\begin{aligned}
\|\u^{n+1}\|_{0}^{2}-\|\tilde{\u}^{n+1}\|_{0}^{2}+\tau^2(\|\nabla p^{n+1}\|_{0}^{2}-\|\nabla p^{n}\|_{0}^{2})=2\tau(\nabla p^{n},  \tilde{\u}^{n+1}).
\end{aligned}
\end{eqnarray}
Using (\ref{23})-(\ref{28}) to estimate the right side of (\ref{20})-(\ref{22}),  and then add it to the formula (\ref{3.5}),  there is
\begin{eqnarray}\label{3.6}
\begin{aligned}
&\|\tilde{\u}^{n+1}\|_{0}^{2}-\|\u^{n}\|_{0}^{2}+\|\tilde{\u}^{n+1}-u^{n}\|_{0}^{2}+\|\omega^{n+1}\|_{0}^{2}-\|\omega^{n}\|_{0}^{2}+\|\omega^{n+1}-\omega^{n}\|_{0}^{2}\\
&+\|T^{n+1}\|_{0}^{2}-\|T^{n}\|_{0}^{2}+\|T^{n+1}-T^{n}\|_{0}^{2}+\nu\tau \|\nabla \tilde{\u}^{n+1}\|_{0}^{2}+\a\tau\|\nabla \omega^{n+1}\|_{0}^{2}+\k\tau \|\nabla T^{n+1}\|_{0}^{2}\\
&+4{\nu_r} \tau\|\omega^{n+1}\|_{0}^{2}+\tau^2(\|\nabla p^{n+1}\|_{0}^{2}-\|\nabla p^{n}\|_{0}^{2})+2 \b\tau\|\nabla\cdot\omega^{n+1}\|_{0}^{2}\\
&\leq  c\tau \|f_1(t_{n+1})\|_{0}^{2}+c \tau \|f_2(t_{n+1})\|_{0}^{2}+c \tau \|f_3(t_{n+1})\|_{0}^{2}+c \tau\| T^{n}\|_{0}^{2}+4{\nu_r} \tau\|\omega^{n}\|_{0}^{2}.
\end{aligned}
\end{eqnarray}
Using (\ref{28}), summing the above formula from 0 to $N-1$,  we can get the Theorem3.1. $\hfill\Box$

Now we will get present the optimal error estimates of velocity, pressure and temperature for Algorithm 3.1. In order to simplify the descriptions, we denote
\begin{eqnarray}
\begin{array}{ll}
\tilde{e}_{\mathrm{\u}}^{n}=\u(t_{n})-\tilde{\u}^{n},  &\quad e_{\u}^{n}=\u(t_{n})-\u^{n}, \quad \quad e_{p}^{n}=p(t_{n})-p^{n}, \\
e_{\omega}^{n}=\omega(t_{n})-\omega^{n},  &\quad e_{T}^{n}=T(t_{n})-T^{n}, \\
\xi_{p}^{n}=p(t_{n+1})-p^{n},  &\quad\xi_{\omega}^{n}=\omega(t_{n+1})-\omega^{n},  \quad \xi_{T}^{n}=T(t_{n+1})-T^{n}.
\end{array}
\end{eqnarray}
\textbf{Theorem3.2(Convergence analysis)} Under the Assumptions 1-2, Algorithm 3.1 have
\begin{equation}
\begin{aligned}
&\sum_{n=0}^{N-1}(\|e_{u}^{N}\|_{0}^{2}+\|e_{\omega}^{N}\|_{0}^{2}+\|e_{T}^{N}\|_{0}^{2})
+\sum_{n=0}^{N-1}(\frac{1}{2}\|\tilde{e}_{u}^{n+1}-e_{u}^{n}\|_{0}^{2}+\|e_{\omega}^{n+1}-e_{\omega}^{n}\|_{0}^{2}+\|e_{T}^{n+1}-e_{T}^{n}\|_{0}^{2})\\
+&\sum_{n=0}^{N-1}(\nu \tau\|\nabla \tilde{e}_{u}^{n+1}\|_{0}^{2}+\tau\|\nabla e_{\omega}^{n+1}\|_{0}^{2}+\tau\|\nabla e_{T}^{n+1}\|_{0}^{2}+\tau^{2}\|\nabla e_{p}^{n+1}\|_{0}^{2}+2\tau\|\nabla\cdot e_{\omega}^{n+1}\|_{0}^{2})
\leq c \tau^2,\\
&\|\nabla \tilde{e}_{u}^{N}\|_{0}^{2}+\tau \sum_{n=0}^{N-1}(\|d_{t} e_{u}^{n+1}\|_{0}^{2}+\|e_{p}^{n+1}\|_{0}^{2}) \leq c \tau^{2}
, \quad d_{t} e_{u}^{n+1}=\frac{e_{u}^{n+1}-e_{u}^{n}}{\tau} \text {, for } 0 \leq n \leq N-1 .
\end{aligned}
\end{equation}
\textbf{Proof.} The truncation error $R_u^n, R_\omega^n, R_T^n$ is defined as
\begin{align}
&\frac{1}{\tau} ({\u}(t_{n+1})-{\u}(t_{n}))-(\nu+\nu_r)\Delta{\u}(t_{n+1})+(\u(t_{n+1})\cdot \nabla){\u}(t_{n+1})+\nabla p(t_{n+1})\no\\
&=2\nu_r\nabla \times\omega(t_{n+1})+J\hat{e}T(t_{n+1})+f_1(t_{n+1})+R_u^n,\label{34}\\
&\frac{1}{\tau} ({\omega}(t_{n+1})-{\omega}(t_{n}))- \a\Delta \omega(t_{n+1})+(\u(t_{n+1}) \cdot \nabla) \omega(t_{n+1})+4 {\nu_r}  w(t_{n+1})+\b\nabla\nabla\cdot\omega(t_{n+1})\no\\
&=2{\nu_r}  \nabla \times\u(t_{n+1})+f_2(t_{n+1})+R_\omega^n,\label{35}\\
&\frac{1}{\tau}(T(t_{n+1})-T(t_{n}))-\b\Delta T(t_{n+1})+(\u(t_{n+1}) \cdot \nabla) T(t_{n+1})\no\\
&=D\nabla \times \omega(t_{n+1}) \cdot \nabla T(t_{n+1})+f_3(t_{n+1})+R_T^n,\label{36}\\
&\nabla \cdot \u(t_{n+1})=0.\label{37}
\end{align}
where
\begin{equation}
R_{u}^{n}=\frac{1}{\tau} \int_{t_{n}}^{t_{n+1}}(t-t_{n}) \u_{t t}(t) \mathrm{~d}t,
\quad R_{\omega}^{n}=\frac{1}{\tau} \int_{t_{n}}^{t_{n+1}}(t-t_{n}) \omega_{t t}(t) \mathrm{~d}t,
\quad R_{T}^{n}=\frac{1}{\tau} \int_{t_{n}}^{t_{n+1}}(t-t_{n}) T_{t t}(t) \mathrm{~d}t.
\end{equation}
By subtracting (\ref{3.1.1})-(\ref{3.1.4}), from (\ref{34})-(\ref{37}),  we obtain
\begin{align}
&\frac{1}{\tau}(\tilde{e}_{u}^{n+1}-e_{u}^{n})-({\nu}+{\nu_r})\Delta{\tilde{e}_{u}^{n+1}}(t_{n+1})\no \\
&=(\u^{n}\cdot \nabla)\tilde{\u}^{n+1}-(\u(t_{n+1})\cdot \nabla){\u}(t_{n+1})-\nabla
\xi_{p}^{n}+2\nu_r\nabla \times \xi_{\omega}^{n}+J\hat{e}\xi_{T}^{n}+R_u^n,\label{39}\\
&\frac{1}{\tau}(e_{T}^{n+1}-e_{T}^{n})-\a\Delta e_{\omega}^{n+1}\no\\
&=(\u^{n+1} \cdot \nabla) \omega^{n+1}-(\u(t_{n+1}) \cdot \nabla) \omega(t_{n+1})-4\nu_re_{\omega}^{n+1}+2\nu_r\nabla \times e_{u}^{n+1}-\b\nabla\nabla\cdot e_{\omega}^{n+1}+R_{\omega}^{n},\label{40}\\
&\frac{1}{\tau}(e_{T}^{n+1}-e_{T}^{n})-\k\Delta e_{T}^{n+1}\label{41}\\
&=(\u^{n+1} \cdot \nabla) T^{n+1}-(\u(t_{n+1}) \cdot \nabla) T(t_{n+1})+D\nabla \times\omega(t_{n+1}) \cdot \nabla T(t_{n+1})-D\nabla \times\omega^{n+1} \cdot \nabla T^{n+1}+R_{T}^{n},\no \\
&\nabla \cdot \mathbf{e}_{u}^{n+1}=0.
\end{align}
Taking the inner product of (\ref{39})-(\ref{41}) with $(2\tau \tilde e_{{u}}^{n+1}, 2\tau e_{\omega}^{n+1}, 2\tau e_{T}^{n+1})$, using (\ref{Lemma 2N}), we obtain
\begin{align}
&\|\tilde{e}_{u}^{n+1}\|_{0}^{2}-\|e_{u}^{n}\|_{0}^{2}+\|\tilde{e}_{u}^{n+1}-e_{u}^{n}\|_{0}^{2}+2\tau(\nu+{\nu_r})\|\nabla \tilde{e}_{u}^{n+1}\|_{0}^{2}\no \\
=&2\tau \n(\u(t_{n})-\u(t_{n+1}) ; \tilde{\u}^{n+1},  \tilde{e}_{\u}^{n+1})-2\tau \n(e_{u}^{n} ; \tilde{\u}^{n+1},  \tilde{e}_{u}^{n+1}))\no\\
&-2 \tau(\nabla \xi_{p}^{n}, \tilde{e}_{u}^{n+1})+4\nu_r\tau(\nabla \times\xi_{\omega}^{n}, \tilde{e}_{u}^{n+1})+2 \hat{e}\tau (J\xi_{T}^{n},  \tilde{e}_{u}^{n+1})+2 \tau(R_{u}^{n},  \tilde{e}_{u}^{n+1}) \label{43}\\
\leq& \sum_{i=1}^{6} \i_{i}, \no\\
&\|e_{\omega}^{n+1}\|_{0}^{2}-\|e_{\omega}^{n}\|_{0}^{2}+\|e_{\omega}^{n+1}-e_{\omega}^{n}\|_{0}^{2}+2\a\tau\|\nabla e_{\omega}^{n+1}\|_{0}^{2}+2\b\tau\|\nabla\cdot e_{\omega}^{n+1}\|_{0}^{2} \no\\
=&-2 \tau \bn(e_{u}^{n+1} ; \omega(t_{n+1}),  e_{\omega}^{n+1})-8\nu_r\tau\|e_{\omega}^{n+1}\|_{0}^{2}+4\nu_r\tau(\nabla \times e_{u}^{n+1}, e_{\omega}^{n+1})+2\tau(R_{\omega}^{n},  e_{\omega}^{n+1}) \label{44}\\
\leq& \sum_{i=7}^{10} \i_{i}, \no\\
&\|e_{T}^{n+1}\|_{0}^{2}-\|e_{T}^{n}\|_{0}^{2}+\|e_{T}^{n+1}-e_{T}^{n}\|_{0}^{2}+2\k\tau\|\nabla e_{T}^{n+1}\|_{0}^{2} \no\\
=&-2 \tau \bn(e_{u}^{n+1} ; T(t_{n+1}),  e_{T}^{n+1})+2D\tau \m(e_{\omega}^{n+1};T(t_{n+1}), e_{T}^{n+1})+2\tau(R_{T}^{n},  e_{T}^{n+1})\label{45}\\
\leq& \sum_{i=11}^{13} \i_{i}\no.
\end{align}
Now,  we majorize the right-hand side terms of (\ref{20})-(\ref{22}) as follows: by using (\ref{3line}),  (\ref{3line2}), Lemma 2.2 and Cauchy-Schwarz inequality,  we have
\begin{align*}
\i_{1} &={2\tau}\n(\u(t_{n})-\u(t_{n+1}) ; \tilde{\u}^{n+1},  \tilde{e}_{u}^{n+1})=-{2\tau}\n(\u(t_{n})-\u(t_{n+1}) ; \tilde{e}_{u}^{n+1},  \tilde{\u}^{n+1})\\
&=-{2\tau}\n(\u(t_{n})-\u(t_{n+1}) ; \tilde{e}_{u}^{n+1},  \u(t_{n+1})) \leq c \tau\|\u(t_{n})-\u(t_{n+1})\|_{0}\|\nabla \tilde{e}_{u}^{n+1}\|_{0}\|\u(t_{n+1})\|_{2}\\
&\leq C \tau\|\tilde{e}_{u}^{n+1}\|_{0}\|\int_{t_{n}}^{t_{n+1}} \u_{t} \mathrm{~d}t\|_{0}
\leq \frac{\nu\tau}{4}\|\nabla \tilde{e}_{\u}^{n+1}\|_{0}^{2}+C \tau^{2} \int_{t_{n}}^{t_{n+1}}|\u_{t}|^{2} \mathrm{~d}t,
\\
\i_{2} &\leq2 \tau \n(e_{u}^{n} ; \tilde{e}_{u}^{n+1},  \u(t_{n+1})) \leq c \tau\|e_{u}^{n}\|_{0}\|\nabla \tilde{e}_{u}^{n+1}\|_{0}\|u(t_{n+1})\|_{2}\\
& \leq C \tau\|e_{u}^{n}\|_{0}\|\nabla \tilde{e}_{u}^{n+1}\|_{0}\leq \frac{\nu\tau}{4}\|\nabla \tilde{e}_{u}^{n+1}\|_{0}^{2}+C \tau\|e_{u}^{n}\|_{0}^{2},
\\
\i_{3} &\leq2\tau(\nabla \xi_{p}^{n},  \tilde{e}_{u}^{n+1}-e_{u}^{n})\leq \frac{1}{2}\|\tilde{e}_{u}^{n+1}-e_{u}^{n}\|_{0}^{2}+2 \tau^{2}\|\nabla \xi_{p}^{n}\|_{0}^{2},
\\
\i_{4} &\leq4\nu_r\tau(\xi_{\omega}^{n}, \nabla \times\tilde{e}_{u}^{n+1})\leq 4\nu_r \tau\|\xi_{\omega}^{n}\|_{0}\|\nabla\tilde{e}_{u}^{n+1}\|_{0},
\leq \nu_r\tau\|\nabla\tilde{e}_{u}^{n+1}\|_{0}^{2}+4\nu_r\tau\|\xi_{\omega}^{n}\|_{0}^{2},
\\
\i_{5} &\leq 2 \tau \hat{e}\|\xi_{T}^{n}\|_{0}\|\tilde{e}_{u}^{n+1}\|_{0} \leq c \tau\|\xi_{T}^{n}\|_{0}^{2}+\frac{\nu\tau}{4}\|\nabla \tilde{e}_{u}^{n+1}\|_{0}^{2}\\
& \leq c \tau\|e_{T}^{n+1}\|_{0}^{2}+c \tau\|T^{n+1}-T^{n}\|_{0}^{2}+\frac{\nu\tau}{4}\|\nabla \tilde{e}_{u}^{n+1}\|_{0}^{2},
\\
\i_{6} & \leq \frac{\nu\tau}{4}\|\nabla \tilde{e}_{u}^{n+1}\|_{0}^{2}+c \tau\|R_{u}^{n}\|_{-1}^{2} =\frac{\nu\tau}{4}\|\nabla \tilde{e}_{u}^{n+1}\|_{0}^{2}+c
\tau^{-1}\|\int_{t_{n}}^{t_{n+1}}(t-t_{n}) \u_{tt} \mathrm{~d}t\|_{-1}^{2}\\
& \leq \frac{\nu\tau}{4}\|\nabla \tilde{e}_{u}^{n+1}\|_{0}^{2}+c \tau^{-1} \int_{t_{n}}^{t_{n+1}}\|\u_{tt}\|_{-1}^{2} \mathrm{~d}t \int_{t_{n}}^{t_{n+1}}(t-t_{n})^{2} \mathrm{~d}t\\
& \leq \frac{\nu\tau}{4}\|\nabla \tilde{e}_{u}^{n+1}\|_{0}^{2}+c \tau^{2} \int_{t_{n}}^{t_{n+1}}\|\u_{tt}\|_{-1}^{2} \mathrm{~d}t,
\end{align*}
and
\begin{align*}
\i_{7} &=-2 \tau \bn(e_{u}^{n+1} ; \omega(t_{n+1}),  e_{\omega}^{n+1})\leq C \tau\|e_{u}^{n+1}\|_{0}\|\nabla e_{\omega}^{n+1}\|_{0} \leq \frac{\a\tau}{3}\|\nabla e_{\omega}^{n+1}\|_{0}^{2}+C \tau\|e_{u}^{n+1}\|_{0}^{2},
\\
\i_{8} &=-8\nu_r\tau\|e_{\omega}^{n+1}\|_{0}^{2},\\
\i_{9} &\leq4\nu_r(\nabla \times e_{u}^{n+1}, e_{\omega}^{n+1})\leq \nu_r\tau\|\nabla e_{u}^{n+1}\|_{0}^{2}+4\nu_r\tau\|e_{\omega}^{n+1}\|_{0}^{2},
\\
\i_{10}&\leq \frac{\a\tau}{3}\|\nabla e_{\omega}^{n+1}\|_{0}^{2}+c \tau\|R_{\omega}^{n}\|_{-1}^{2}  =\frac{\a\tau}{3}\|\nabla e_{\omega}^{n+1}\|_{0}^{2}+c \tau^{-1}\|\int_{t_{n}}^{t_{n+1}}(t-t_{n}) \omega_{t t} \mathrm{~d}t\|_{-1}^{2}\\
&\leq \frac{\a\tau}{3}\|\nabla e_{\omega}^{n+1}\|_{0}^{2}+c \tau^{2} \int_{t_{n}}^{t_{n+1}}\|\omega_{t t}\|_{-1}^{2} \mathrm{~d}t,
\end{align*}
and
\begin{align*}
\i_{11}&=2 \tau \bn(e_{u}^{n+1} ; e_{T}^{n+1},  T(t_{n+1})) \leq c \tau\|e_{\mathrm{u}}^{n+1}\|_{0}\|\nabla e_{T}^{n+1}\|_{0}\|T(t_{n+1})\|_{2} \\
& \leq C \tau\|e_{u}^{n+1}\|_{0}\|\nabla e_{T}^{n+1}\|_{0} \leq \frac{\k\tau}{2}\|\nabla e_{T}^{n+1}\|_{0}^{2}+C \tau\|e_{u}^{n+1}\|_{0}^{2},
\\
\i_{12}&=2D\tau \m(e_{\omega}^{n+1};T(t_{n+1}), e_{T}^{n+1})\leq c \tau\|\nabla e_{\omega}^{n+1}\|_{0}\|\nabla T(t_{n+1})\|_{1}\|e_{T}^{n+1}\|_{1} \\
& \leq C \tau\|\nabla e_\omega^{n+1}\|_{0}\| e_{T}^{n+1}\|_{1} \leq \frac{C \tau}{2}\| e_{T}^{n+1}\|_{1}^{2}+ \frac{\a\tau}{3}\|\nabla e_\omega^{n+1}\|_{0}^{2},
\\
\i_{13}&\leq \frac{\k\tau}{2}\|\nabla e_{T}^{n+1}\|_{0}^{2}+c \tau\|R_{T}^{n}\|_{-1}^{2}=\frac{ \k\tau}{2}\|\nabla e_{T}^{n+1}\|_{0}^{2}+c \tau^{-1}\|\int_{t_{n}}^{t_{n+1}}(t-t_{n}) T_{t t} \mathrm{~d}t\|_{-1}^{2} \\
\quad & \leq \frac{ \k\tau}{2}\|\nabla e_{T}^{n+1}\|_{0}^{2}+c \tau^{2} \int_{t_{n}}^{t_{n+1}}\|T_{t t}\|_{-1}^{2} \mathrm{~d}t.
\end{align*}
On the other hand,  from $(\ref{3.1.2.2})$, we derive
\begin{equation}
\begin{aligned}\label{46}
e_{\mathrm{u}}^{n+1}&=\tilde{e}_{\mathrm{u}}^{n+1}+\tau \nabla \phi^{n+1}, \\
2 \tau^{2}\|\nabla \xi_{p}^{n}\|_{0}^{2}&=2 \tau^{2}\|\nabla e_{p}^{n+1}+\nabla \phi^{n+1}\|_{0}^{2}\leq c \tau^{2}\|\nabla e_{p}^{n+1}\|_{0}^{2}+\tau^{2}\|\nabla \phi^{n+1}\|_{0}^{2}.
\end{aligned}
\end{equation}
And taking the inner product of the first equation of $(\ref{46})$ with itself on both sides,  we obtain
\begin{equation}
\begin{aligned}\label{3.26}
\|\tilde{e}_{\mathrm{u}}^{n+1}\|_{0}^{2}&=\|e_{u}^{n+1}\|_{0}^{2}+\tau^{2}\|\nabla \phi^{n+1}\|_{0}^{2}.
\end{aligned}
\end{equation}
Using the above inequality and (\ref{43})-(\ref{3.26}), we get
\begin{align}
&\|e_{u}^{n+1}\|_{0}^{2}-\|e_{u}^{n}\|_{0}^{2}+\frac{1}{2}\|\tilde{e}_{u}^{n+1}-e_{u}^{n}\|_{0}^{2}+({\nu+\nu_r}) \tau\|\nabla \tilde{e}_{u}^{n+1}\|_{0}^{2}+\tau^{2}\|\nabla e_{p}^{n+1}\|_{0}^{2} \no\\
\leq & c \tau\left(\|e_{u}^{n}\|_{0}^{2}+\|e_{T}^{n+1}\|_{0}^{2}\right)+4\nu_r\tau\|\xi_{\omega}^{n}\|_{0}^{2}+c\tau^{2}\int_{t_{n}}^{t_{n+1}}\|\u_{tt}\|_{-1}^{2} \mathrm{~d}t+c\tau^{2}\int_{t_{n}}^{t_{n+1}}\left|\u_{t}\right|^{2} \mathrm{~d}t,\label{48}
\\
&\|e_{\omega}^{n+1}\|_{0}^{2}-\|e_{\omega}^{n}\|_{0}^{2}+\|e_{\omega}^{n+1}-e_{\omega}^{n}\|_{0}^{2}+\frac{4\a}{3}\tau\|\nabla e_{\omega}^{n+1}\|_{0}^{2}+2\b\tau\|\nabla\cdot e_{\omega}^{n+1}\|_{0}^{2} \no\\
\leq & c \tau\|e_{u}^{n+1}\|_{0}^{2}-8\nu_r\tau\|e_{\omega}^{n+1}\|_{0}^{2}+4\nu_r\tau\|e_{\omega}^{n+1}\|_{0}^{2}-\nu_r\tau\|\nabla \tilde{e}_{u}^{n+1}\|_{0}^{2}+c\tau^{2} \int_{t_{n}}^{t_{n+1}}\|\omega_{tt}\|_{-1}^{2}\mathrm{~d}t,\label{49}
\\
&\|e_{T}^{n+1}\|_{0}^{2}-\|e_{T}^{n}\|_{0}^{2}+\|e_{T}^{n+1}-e_{T}^{n}\|_{0}^{2}+\k\tau\|\nabla e_{T}^{n+1}\|_{0}^{2} \no\\
\leq &  c \tau\left(\|e_{u}^{n}\|_{0}^{2}+\|e_{T}^{n+1}\|_{1}^{2}\right)+\frac{1}{3}\tau\|\nabla e_\omega^{n+1}\|_{0}^{2}+c \tau^{2} \int_{t_{n}}^{t_{n+1}}\|T_{tt}\|_{-1}^{2} \mathrm{~d}t.\label{50}
\end{align}
Taking the sum of (\ref{48})-(\ref{50}) for $n$ from $0$ to $N-1$, using the regularity of $u$, $\omega$ and $T$, by applying the discrete Gronwall's lemma, we have
\begin{equation}
\begin{aligned}\label{3.33}
&\sum_{n=0}^{N-1}(\|e_{u}^{N}\|_{0}^{2}+\|e_{\omega}^{N}\|_{0}^{2}+\|e_{T}^{N}\|_{0}^{2})
+\sum_{n=0}^{N-1}(\frac{1}{2}\|\tilde{e}_{u}^{n+1}-e_{u}^{n}\|_{0}^{2}+\|e_{\omega}^{n+1}-e_{\omega}^{n}\|_{0}^{2}+\|e_{T}^{n+1}-e_{T}^{n}\|_{0}^{2})\\
+&\sum_{n=0}^{N-1}(\nu \tau\|\nabla \tilde{e}_{u}^{n+1}\|_{0}^{2}+\a\tau\|\nabla e_{\omega}^{n+1}\|_{0}^{2}+\k\tau\|\nabla e_{T}^{n+1}\|_{0}^{2}+\tau^{2}\|\nabla e_{p}^{n+1}\|_{0}^{2}+2\b\tau\|\nabla\cdot e_{\omega}^{n+1}\|_{0}^{2})
\leq c\tau^2.
\end{aligned}
\end{equation}

Then we will prove the convergence of $p$ in the following. In conjunction with $(\ref{3.1.2.2})$, we write $(\ref{39})$ as
\begin{align}\label{3.27}
&\frac{1}{\tau}({{e_{u}^{n+1}}-{e_{u}^{n}}})-(\nu+\nu_r)\Delta{\tilde{e}_{u}^{n+1}}\no\\
=&(\u^{n}\cdot \nabla)\tilde{\u}^{n+1})-(\u(t_{n+1})\cdot \nabla){\u}(t_{n+1}))
-\nabla e_{p}^{n+1}+2\nu_r\nabla \times\xi_{\omega}^{n}+J \hat{e}\xi_{T}^{n}+R_u^n.
\end{align}
We recall the existence of $\exists \b>0$ such that (inf-sup condition)
\begin{equation}\label{3.28}
\b\|q\|_{0}\leq \sup _{{v} \in H_{0}^{1}(\Omega)} \frac{(q,  \nabla \cdot {v})}{\|\nabla {v}\|_{0}}  ,  \quad \forall q \in L_{0}^{2}(\Omega).
\end{equation}
Now,  we use $\|\nabla \tilde{e}_{\mathrm{u}}^{n+1}\|_{0}$ to estimate $(e_{p}^{n+1}, \nabla \cdot \tilde{e}_{u}^{n+1})$,  multiply both sides of $(\ref{3.27})$ $\tilde{e}_{\mathrm {u}}^{n+1}$,  we can get
\begin{equation}
\begin{aligned}\label{3.29}
&(e_{p}^{n+1},  \nabla \cdot \tilde{e}_{\mathrm{u}}^{n+1}) \\
=&\frac{1}{\tau}(({e_{u}^{n+1}}-{e_{u}^{n}}),  \tilde{e}_{u}^{n+1})+(\nu+{\nu_r})(\nabla \tilde{e}_{u}^{n+1}, \nabla \tilde{e}_{u}^{n+1})+((\u(t_{n+1})-\u(t_{n})) \cdot \nabla \tilde{u}^{n+1},  \tilde{e}_{u}^{n+1})
\\
&+((e_{u}^{n} \cdot \nabla) \tilde{\u}^{n+1},  \tilde{e}_{u}^{n+1})-2\nu_r(\nabla \times\xi_{\omega}^{n}, \tilde{e}_{u}^{n+1})-\hat{e}(J(T(t_{n+1})-T^{n}),  \tilde{e}_{u}^{n+1})-(R_{u}^{n},  \tilde{e}_{u}^{n+1})\\
=&\sum_{i=1}^{7} \theta_{i}.
\end{aligned}
\end{equation}
We now proceed to estimate term $\theta_{1}$ to $\theta_{7}$,  respectively
\begin{equation}
\begin{aligned}\label{3.30}
&\theta_{1} \leq \frac{1}{\tau}\|e_{u}^{n+1}-e_{u}^{n}\|_{0}\|\tilde{e}_{u}^{n+1}\|_{0} \leq \frac{c}{\tau}\| e_{u}^{n+1}-e_{u}^{n}\|_{0}\|\nabla \tilde{e}_{u}^{n+1}\|_{0},\\
&\theta_{2}=(\nu+{\nu_r})\|\nabla \tilde{e}_{u}^{n+1}\|_{0}^{2},\\
&\theta_{3} \leq c\|\u(t_{n})-\u(t_{n+1})\|_{0}\|\nabla \tilde{e}_{u}^{n+1}\|_{0}\|\tilde{\u}^{n+1}\|_{2} \leq C\int_{t_{n}}^{t_{n+1}} \|\u_{t}\|_{0} \mathrm{~d} t \|\nabla \tilde{e}_{u}^{n+1}\|_{0},\\
&\theta_{4} \leq c\|e_{u}^{n}\|_{0}\|\nabla \tilde{e}_{u}^{n+1}\|_{0}\|\tilde{u}^{n+1}\|_{2} \leq C\|e_{u}^{n}\|_{0}\|\nabla \tilde{e}_{u}^{n+1}\|_{0},\\
&\theta_{5} \leq 2\nu_r\|\xi_{\omega}^{n}\|_{0}\|\nabla \times\tilde{e}_{u}^{n+1}\|_{0} \leq 2\nu_r\|\xi_{\omega}^{n}\|_{0}\|\nabla\tilde{e}_{u}^{n+1}\|_{0},\\
&\theta_{6}\leq \hat{e}\|\xi_{T}^{n}\|_{0}\|\nabla \tilde{e}_{u}^{n+1}\|_{0},\\
&\theta_{7} \leq c\|R_{u}^{n}\|_{-1}\|\nabla \tilde{e}_{u}^{n+1}\|_{0} \leq c \tau \int_{t_{n}}^{t_{n+1}}\|u_{t t}\|_{-1} \mathrm{~d}t\|\nabla \tilde{e}_{u}^{n+1}\|_{0}.
\end{aligned}
\end{equation}
From $(\ref{3.27})$ and $(\ref{3.28})$, we have
\begin{equation}
\begin{aligned}\label{3.31}
c\|e_{p}^{n+1}\|_{0} \leq & \frac{1}{\tau}\|e_{u}^{n+1}-e_{u}^{n}\|_{0}+(\nu+{\nu_r})\|\nabla \tilde{e}_{u}^{n+1}\|_{0}+C\int_{t_{n}}^{t_{n+1}} \|\u_{t}\|_{0} \mathrm{~d} t+C\|e_{u}^{n}\|_{0}\\
&+\|\xi_{\omega}^{n}\|_{0}+\|\xi_{T}^{n}\|_{0}+c \tau \int_{t_{n}}^{t_{n+1}}\|\u_{tt}\|_{-1} \mathrm{~d}t.
\end{aligned}
\end{equation}
Multiplying by $\tau$ in (\ref{3.31}) and squaring, using Cauchy-Schwarz inequality, from (\ref{3.33})-(\ref{3.30}) we can get
\begin{equation}
\tau \sum_{i=0}^{N-1}\|e_{p}^{n+1}\|_{0}^{2} \leq c \tau^{2}.
\end{equation}
Then we complete the proof.$\hfill\Box$

Adopting semi-implicit scheme for the nonlinear term and BDF2 to march in time, the second-order pressure correction schemes are as follows:\\
\textbf{Algorithm 3.3 (Second-order RPC schemes)}\\
Assume that at each time step, $(\u^{n-1}, p^{n-1}, \omega^{n-1}, T^{n-1})$ and $(\u^{n}, p^{n}, \omega^{n}, T^{n})$ are given, where ($\u^{1}$, $p^{1}$, $\omega^{1}$, $T^{1}$) is obtained by the first-order RPC schemes, solve for $(\u^{n+1}, p^{n+1}, \omega^{n+1}, T^{n+1})$ from\\
Step1: Find $\tilde{\u}^{n+1}$ as the solution of
\begin{eqnarray}\label{3.2.1}
\left\{\begin{array}{lll}
\frac{1}{2\tau}(3\tilde{\u}^{n+1} -4{\u}^{n}+{\u}^{n-1})-(\nu+{\nu_r}) \Delta \tilde{\u}^{n+1}+(\u^n \cdot \nabla)\tilde{\u}^{n+1}+\nabla p^n
=2\nu_r\nabla \times\omega^n+J \hat{e} T^n+f_1^{n+1},\\
\tilde{\u}^{n+1}|_{\partial\Omega}=0.
\end{array}\right.
\end{eqnarray}
Step2.1: Find $\phi^{n+1}$ as the solution of
\begin{eqnarray}\label{3.2.2.1}
\left\{\begin{array}{l}
\Delta \phi^{n+1}=\frac{3}{2\tau} \nabla \cdot \tilde {\u}^{n+1}, \\
\frac{\partial \phi}{\partial n}|_{\partial \Omega}=0.
\end{array}\right.
\end{eqnarray}
Step2.2: Update $\u^{n+1}, p^{n+1}$
\begin{eqnarray}
\left\{\begin{array}{l}
\u^{n+1}=\tilde{\u}^{n+1}-\frac{2\tau}{3}\nabla \phi^{n+1}, \\
p^{n+1}=\phi^{n+1}+p^{n}-(\nu+{\nu_r})\nabla \cdot \tilde{\u}^{n+1}.
\end{array}\right.
\end{eqnarray}
Step3: Compute $\omega^{n+1}$ as the solution of
\begin{eqnarray}\label{3.2.3}
\left\{\begin{array}{lll}
\frac{1}{2\tau}(3\omega^{n+1}-4\omega^{n}+\omega^{n-1}) -\a\Delta \omega^{n+1} +(\u^{n+1} \cdot \nabla)\omega^{n+1}+4\nu_r\omega^{n+1}\\-\b\nabla\nabla\cdot\omega^{n+1} =2\nu_r\nabla \times\u^{n+1}+f_2^{n+1},\\
\omega^{n+1}|_{\partial\Omega}=0.
\end{array}\right.
\end{eqnarray}
Step4: Compute $T^{n+1}$ as the solution of
\begin{eqnarray}\label{3.2.4}
\left\{\begin{array}{l}
\frac{1}{2\tau}(3T^{n+1}-4T^{n}+T^{n-1})-\k\Delta T^{n+1}+(\u^{n+1} \cdot \nabla) T^{n+1}-D\nabla \times\omega^{n+1} \cdot \nabla T^{n+1}=f_3^{n+1},\\
T^{n+1}|_{\partial\Omega}=0.
\end{array}\right.
\end{eqnarray}

\section{The RPC schemes: full discretization}
In this section,  we will introduce the first-order full-discrete scheme and the second-order full-discrete scheme of the pressure correction algorithm for the thermopolar fluid equations. Let $\mathcal{K}_{h}=\{K\}$ be a uniformly regular family of triangulation of $\Omega$,  and define the mesh size $h=\mathop{\max}\limits_{K \in \mathcal{K}_{h}}\{{\rm diam}(K)\}$. The spatial approximation of fluid velocity,  pressure, angular velocity and temperature are applied by mixed element method with $({\rm X}_h, {\rm Q}_h, {\rm W}_h, {\rm M}_h) \in ({\rm X}, {\rm Q}, {\rm X}, {\rm M})$.  Next,  we define the following discrete subspaces:
$$
\begin{aligned}
&{\rm X}_{h}=\left\{v_{h} \in {\rm X} \cap C^{0}(\bar{\Omega})^{d}:\left.{v}_{h}\right|_{K} \in P_{1}(K)^{d} \oplus {\rm span}\left\{b\right\},  \forall K \in \mathcal{K}_{h}\right\}, \\
&{\rm Q}_{h}=\left\{q_{h} \in {\rm Q} \cap C^{0}(\bar{\Omega}):\left.q_{h}\right|_{K} \in P_{1}(K),  \forall K \in \mathcal{K}_{h}\right\}, \\
&{\rm W}_{h}=\left\{v_{h} \in {\rm X} \cap C^{0}(\bar{\Omega})^{d}:\left.{v}_{h}\right|_{K} \in P_{1}(K)^{d},  \forall K \in \mathcal{K}_{h}\right\}, \\
&{\rm M}_{h}=\left\{v_h \in M \cap C^{0}(\bar{\Omega}):\left.v_{h}\right|_{K} \in P_{1}(K),  \forall K \in \mathcal{K}_{h}\right\}.\\
\end{aligned}
$$

In this section,  we also consider using BDF1 and BDF2 to approximate the time derivative of the time discretization of the problem (\ref{1}). Let $\tau>0$ be the time step,  we have $t_n=n\tau>0$.

\subsection{First-order RPC schemes}
For the full discretization, the finite element method is used to discretize the space. Adopting semi-implicit scheme for the nonlinear term and BDF1 to march in time and letting $\phi^{n+1}_h=(p^{n+1}_h-p^n_h)$, the first-order pressure correction schemes are as follows:\\
\textbf{Algorithm 4.1 First-order full discretization SPC schemes}\\
Assume that at each time step, $(\u^n_h, p^n_h, \omega^n_h, T^n_h)$ are given and one seeks $(\u^{n+1}_h, p^{n+1}_h, \omega^{n+1}_h, T^{n+1}_h)$\\
Step1: Find $\tilde{\u}^{n+1}_h$ as the solution of
\begin{eqnarray}\label{4.1.1}
\left\{\begin{array}{lll}
\frac{1}{\tau}(\tilde{\u}^{n+1}_h-{\u}^{n}_h,{\v}_{h})+(\nu+{\nu_r} ) a(\tilde{\u}^{n+1}_h, {\v}_{h})+\n(\u^n_h;\tilde{\u}^{n+1}_h, {\v}_{h})+(\nabla p^n_h, v_h)\\
=2\nu_r\bar r(\omega^n, {\v}_{h})+\hat{e}(J T^n_h, {v}_h)+(f_1^{n+1}, v_h),\\
\tilde {\u}^{n+1}_h|_{\partial\Omega}=0.
\end{array}\right.
\end{eqnarray}
Step2: Compute $\u_{h}^{n+1}, p_h^{n+1}$ as the solution of
\begin{eqnarray}\label{4.1.2}
\left\{\begin{array}{lll}
\frac{1}{\tau}({\u}^{n+1}_h-\tilde{\u}^{n+1}_h)+\nabla(p^{n+1}_h-p^{n}_h)=0, \\
\nabla \cdot \u^{n+1}_h=0,\\
\u^{n+1}_h \cdot n|_{\partial \Omega}=0.
\end{array}\right.
\end{eqnarray}
Step3: Compute $\omega_h^{n+1}$ as the solution of
\begin{eqnarray}\label{4.1.3}
\left\{\begin{array}{lll}
\frac{1}{\tau}({\omega}^{n+1}_h-{\omega}^{n}_h, \s_h)+\a\bar a(\omega^{n+1}_h, \s_h) +\bn(\u^{n+1}_h;\omega^{n+1}_h, \s_h)+4\nu_r\bar r(\omega^{n+1}_h,  \s_h)\\
+\b(\nabla\cdot\omega^{n+1}_h,\nabla\cdot \s_h)=2\nu_rr(\u^{n+1}_h, \s_h)+(f_2^{n+1}, \s_h),\\
\omega^{n+1}_h|_{\partial\Omega}=0.
\end{array}\right.
\end{eqnarray}
Step4: Compute $T_h^{n+1}$ as the solution of
\begin{eqnarray}\label{4.1.4}
\left\{\begin{array}{l}
\frac{1}{\tau}(T^{n+1}_h-T^{n}_h, \sk_h)+\k\bar a( T^{n+1}_h, \sk_h)+\bn(\u^{n+1}_h ;T^{n+1}_h, \sk_h)\\
=D\m(\omega^{n+1}_h;\nabla T^{n+1}_h, \sk_h)+(f_3^{n+1}, \sk_h),\\
T^{n+1}_h|_{\partial\Omega}=0.
\end{array}\right.
\end{eqnarray}
Similarly, Algorithm 4.1 can also be updated to the following rotation form\\
\textbf{Algorithm 4.2(First-order full discretization RPC schemes)}\\
Step2(RPC): Compute $u^{n+1}, p^{n+1}$ as the solution of
\begin{eqnarray}\label{4.1.RPC}
\left\{\begin{array}{lll}
\frac{1}{\tau}({\u}^{n+1}_h-\tilde{\u}^{n+1}_h)+\nabla(p^{n+1}_h-p^{n}_h+(\nu+\nu_r)\nabla \cdot \tilde{\u}^{n+1})=0, \\
\nabla \cdot \u^{n+1}_h=0,\\
\u^{n+1}_h \cdot n|_{\partial \Omega}=0.
\end{array}\right.
\end{eqnarray}
\textbf{Remark 4.1.} In the process of computing,  (\ref{4.1.RPC}) is equivalent to the following equations:\\
Step2.1: Find $\phi^{n+1}_h$ as the solution of
\begin{eqnarray}\label{4.1.2.1}
\left\{\begin{array}{l}
(\nabla \phi^{n+1}_h, \nabla q^{n+1}_h) =\frac{1}{\tau} (\tilde {\u}^{n+1}_h, \nabla q^{n+1}),  \\
\frac{\partial \phi_h}{\partial n}|_{\partial \Omega}=0.
\end{array}\right.
\end{eqnarray}
Step2.2: Update $\u^{n+1}_h, p^{n+1}_h$
\begin{eqnarray}\label{4.1.2.2}
\left\{\begin{array}{l}
\u^{n+1}_h=\tilde{\u}^{n+1}_h-\tau \nabla \phi^{n+1}_h, \\
p^{n+1}_h=\phi^{n+1}_h+p^{n}_h-(\nu+{\nu_r})\nabla \cdot \tilde{\u}^{n+1}_h.
\end{array}\right.
\end{eqnarray}
\textbf{Theorem 4.1(Stability)} Algorithm 4.1 is unconditionally stable under the Assumptions 1-3, if $(\u^{n+1}$, $p^{n+1}$, $\omega^{n+1}$, $T ^{ n+1})$ is the solution of Algorithm 4.1, then for all $0\leq n\leq N-1$ there is a constant $c$ such that:
\begin{align*}
&\|\u^{N}_h\|_{0}^{2}+\tau^2\|\nabla p^{N}_h\|_{0}^{2}+(1+4\nu_r\tau)\|\omega^{N}_h\|_{0}^{2}+\|T^{N}_h\|_{0}^{2}
+\nu \tau \sum_{n=0}^{N-1}\|\nabla \tilde{\u}^{n+1}_h\|_{0}^{2}+ \a\tau \sum_{n=0}^{N-1}\|\nabla T^{n+1}_h\|_{0}^{2}\\
+&\k\tau\sum_{n=0}^{N-1}\|\nabla \omega^{n+1}_h\|_{0}^{2}
+2\b\tau\sum_{n=0}^{N-1}\|\nabla\cdot\omega^{n+1}_h\|_{0}^{2}+\sum_{n=0}^{N-1}(\|\tilde{\u}^{n+1}_h-u^{n}_h\|_{0}^{2}+\|\omega^{n+1}_h-\omega^{n}_h\|_{0}^{2}+
\|T^{n+1}_h-T^{n}_h\|_{0}^{2})\\
\leq& \|\u^{0}\|_{0}^{2}+\tau^2\|\nabla p^{0}\|_{0}^{2}+(1+4\nu_r\tau)\|\omega^{0}\|_{0}^{2}+\|T^{0}\|_{0}^{2}+c\tau\sum_{n=0}^{N-1}(\|f_1(t_{n+1})\|_{0}^{2}+ \|f_2(t_{n+1})\|_{0}^{2}+\|f_3(t_{n+1})\|_{0}^{2}).
\end{align*}
\textbf{Proof.} Just change the $(2\tau\tilde{\u}^{n+1}, 2\tau \omega^{n+1}, 2\tau T^{n+1})$ in the algorithm 3.1 stability proof changed to $(2\tau \tilde{\u}^{n+1}_h, 2\tau \omega^{n+1}_h, 2\tau T^{n+1}_h)$, and use (\ref{23})-(\ref{28}) to get:
\begin{eqnarray}\label{4.6}
\begin{aligned}
&\|\tilde{\u}^{n+1}_h\|_{0}^{2}-\|\u^{n}_h\|_{0}^{2}+\|\tilde{\u}^{n+1}_h-u^{n}_h\|_{0}^{2}+\|\omega^{n+1}_h\|_{0}^{2}-\|\omega^{n}_h\|_{0}^{2}+\|\omega^{n+1}_h-\omega^{n}_h\|_{0}^{2}\\
+&\|T^{n+1}_h\|_{0}^{2}-\|T^{n}_h\|_{0}^{2}+\|T^{n+1}_h-T^{n}_h\|_{0}^{2}
+\nu\tau \|\nabla \tilde{\u}^{n+1}_h\|_{0}^{2}+\a\tau\|\nabla \omega^{n+1}_h\|_{0}^{2}+\k\tau \|\nabla T^{n+1}_h\|_{0}^{2}\\
+&4{\nu_r} \tau\|\omega^{n+1}_h\|_{0}^{2}
+\tau^2(\|\nabla p^{n+1}_h\|_{0}^{2}-\|\nabla p^{n}_h\|_{0}^{2})+2 \b\tau\|\nabla\cdot\omega^{n+1}_h\|_{0}^{2}\\
\leq&c\tau \|f_1(t_{n+1})\|_{0}^{2}+c\tau \|f_2(t_{n+1})\|_{0}^{2}+c\tau \|f_3(t_{n+1})\|_{0}^{2}+c \tau\| T^{n}_h\|_{0}^{2}+4{\nu_r} \tau\|\omega^{n}_h\|_{0}^{2}.
\end{aligned}
\end{eqnarray}
Using (\ref{4.6}), summing the above formula from 0 to $N-1$,  we have
\begin{align*}
&\|\u^{N}_h\|_{0}^{2}+\tau^2\|\nabla p^{N}_h\|_{0}^{2}+(1+4\nu_r\tau)\|\omega^{N}_h\|_{0}^{2}+\|T^{N}_h\|_{0}^{2}
+\nu \tau \sum_{n=0}^{N-1}\|\nabla \tilde{\u}^{n+1}_h\|_{0}^{2}+ \a\tau \sum_{n=0}^{N-1}\|\nabla T^{n+1}_h\|_{0}^{2}\\
+&\k\tau\sum_{n=0}^{N-1}\|\nabla \omega^{n+1}_h\|_{0}^{2}
+2\b\tau\sum_{n=0}^{N-1}\|\nabla\cdot\omega^{n+1}_h\|_{0}^{2}+\sum_{n=0}^{N-1}(\|\tilde{\u}^{n+1}_h-u^{n}_h\|_{0}^{2}+\|\omega^{n+1}_h-\omega^{n}_h\|_{0}^{2}+
\|T^{n+1}_h-T^{n}_h\|_{0}^{2})\\
\leq& \|\u^{0}\|_{0}^{2}+\tau^2\|\nabla p^{0}\|_{0}^{2}+(1+4\nu_r\tau)\|\omega^{0}\|_{0}^{2}+\|T^{0}\|_{0}^{2}+c\tau\sum_{n=0}^{N-1}(\|f_1(t_{n+1})\|_{0}^{2}+ \|f_2(t_{n+1})\|_{0}^{2}+\|f_3(t_{n+1})\|_{0}^{2}).
\end{align*}
Then we complete the proof.$\hfill\Box$
\subsection{Second-order RPC schemes}
Adopting semi-implicit scheme for the nonlinear terms and BDF2 to march in time, the second-order pressure correction schemes are as follows:\\
\textbf{Algorithm 4.2 Second-order full discretization RPC schemes}\\
Assume that at each time step, $(\u^{n-1}_h, p^{n-1}_h, \omega^{n-1}_h, T^{n-1}_h)$ and $(\u^{n}_h, p^{n}_h, \omega^{n}_h, T^{n}_h)$  are given, where ($\u^{1}_h$, $p^{1}_h$, $\omega^{1}_h$, $T^{1}_h$)  is obtained by the first-order RPC schemes, solve for $(\u^{n+1}_h, p^{n+1}_h, \omega^{n+1}_h, T^{n+1}_h)$ from\\
Step1: Find $\tilde{\u}^{n+1}_h$ as the solution of
\begin{eqnarray}\label{4.2.1}
\left\{\begin{array}{ll}
\frac{1}{2\tau}(3\tilde{\u}^{n+1} -4{\u}^{n}+{\u}^{n-1}, {\v}_{h})+(\nu+{\nu_r} ) a(\tilde{\u}^{n+1}_h, {\v}_{h})+\n(\u^n_h;\tilde{\u}^{n+1}_h, {\v}_{h})+(\nabla p^n_h, \v_h)\\
=2\nu_r\bar r(\omega^n, {\v}_{h})+ \hat{e}(J T^n_h, {\v}_{h})+(f_1^{n+1}, {\v}_{h}),\\
\tilde {\u}^{n+1}_h|_{\partial\Omega}=0.
\end{array}\right.
\end{eqnarray}
Step2.1: Find $\phi^{n+1}_h$ as the solution of
\begin{eqnarray}\label{4.2.2}
\left\{\begin{array}{l}
(\nabla \phi^{n+1}_h, \nabla q_h) =\frac{3}{2\tau} ({\u}^{n+1}_h, \nabla q_h),  \\
\frac{\partial \phi_h}{\partial n}|_{\partial \Omega}=0.
\end{array}\right.
\end{eqnarray}
Step2.2: Update $\u^{n+1}_h, p^{n+1}_h$
\begin{eqnarray}
\left\{\begin{array}{l}
\u^{n+1}_h=\tilde{\u}^{n+1}_h-\frac{2\tau}{3} \nabla \phi^{n+1}_h, \\
p^{n+1}_h=\phi^{n+1}_h+p^{n}_h-(\nu+{\nu_r})\nabla \cdot \tilde{\u}^{n+1}_h.
\end{array}\right.
\end{eqnarray}
Step3: Compute $\omega_h^{n+1}$ as the solution of
\begin{eqnarray}\label{4.2.3}
\left\{\begin{array}{lll}
\frac{1}{2\tau}({3{\omega}^{n+1}_h -4{\omega}^{n}_h+{\omega}^{n-1}_h}, {\s}_{h})+\a\bar a(\omega^{n+1}_h, {\s}_{h}) +\bn(\u^{n+1}_h;\omega^{n+1}_h, {\s}_{h})+4\nu_r\bar r(\omega^{n+1}_h,  {\s}_{h}),\\
+\b(\nabla\cdot\omega^{n+1}_h,\nabla\cdot {\s}_{h})=2\nu_rr(\u^{n+1}_h, {\s}_{h})+(f_2^{n+1}, {\s}_{h}),\\
\omega^{n+1}_h|_{\partial\Omega}=0.
\end{array}\right.
\end{eqnarray}
Step4: Compute $T_h^{n+1}$ as the solution of
\begin{eqnarray}\label{4.2.4}
\left\{\begin{array}{l}
\frac{1}{2\tau}({3{T}^{n+1}_h-4{T}^{n}_h+{T}^{n-1}_h}, \sk_h)+\k\bar a( T^{n+1}_h, \sk_h)+\bn(\u^{n+1}_h ;T^{n+1}_h, \sk_h)\\=D\m(\omega^{n+1}_h;\nabla T^{n+1}_h, \sk_h)+(f_3^{n+1}, \sk_h),\\
T^{n+1}_h|_{\partial\Omega}=0.
\end{array}\right.
\end{eqnarray}
\section{Numerical experiments}
\subsection{Examples with analytical solutions, 2D and 3D}
In this section, the convergence rate of the RPC scheme for the time-dependent thermomicropolar fluid problem with a given true solution is calculated. The calculation formula for the convergence rate relative to the grid size $h$ is $\frac{\log \left(E_{i} / E_{i+1}\right)}{\log \left(h_{i} / h_ {i +1}\right)}$, where $E_i$ and $E_{i+1}$ are the relative errors corresponding to the grid size of $h_i$ and $h_{i+1}$. And choose $\tau=h^2$ and $\tau=h$ as the time step of the first-order format and the second-order format, respectively, which verifies the convergence order of time and space at the same time.
\subsubsection{First-order RPC scheme of 2D}
In the 2D case, the following smooth analytical solution problem is considered, $\nu=\nu_r=\hat{e}=\a=\b=\k=D=1, \T=0.1$, and the initial value in (\ref{1}) is set to the exact solution.
\begin{eqnarray}
\begin{array}{lll}
u_1=  10x^2(x-1)^2y(y-1)(2y-1)\cos(t),\\
u_2=  10y^2(y-1)^2x(1-x)(2x-1)\cos(t),\\
p=  10(2x-1)(2y-1)\cos(t),\\
\omega=  u_1-u_2,\\
T=  u_1+u_2.
\end{array}
\end{eqnarray}

For the first-order format, select $\tau =h^{2}$ as the time step, and use the $P_{1}b-P_{1}-P_{1}-P_{1}$ element to approximate the finite element space. Therefore, the optimal relative error order of $L^2$ in the first-order format $u$, $\omega$ and $T$ is $O(h^{2})$, and $H^1$ relative error optimal order is $O(h^{1})$.
\begin{table}[H]
\centering
\caption{The error of first-order scheme with $\tau =h^{2}$ at $\T=0.1$.}
\begin{tabular*}{16cm}{@{\extracolsep{\fill}}c c c c c c c c c}
\hline
$1/h$ & $\|u-u_h\|_0$  & $\|p-p_h\|_0$  & $\|\omega-\omega_h\|_0$ & $\|T-T_h\|_0$ & $|u-u_h|_1$ & $|\omega-\omega_h|_1$ & $|T-T_h|_1$\\
\hline
10&2.90E-03&8.36E-02&1.86E-03&3.67E-03&7.62E-02&6.09E-02&9.16E-02\\
20&7.31E-04&2.89E-02&4.74E-04&9.51E-04&3.77E-02&3.06E-02&4.66E-02\\
30&3.23E-04&1.51E-02&2.11E-04&4.25E-04&2.50E-02&2.05E-02&3.12E-02\\
40&1.81E-04&9.57E-03&1.19E-04&2.40E-04&1.87E-02&1.53E-02&2.34E-02\\
50&1.15E-04&6.76E-03&7.61E-05&1.54E-04&1.50E-02&1.23E-02&1.87E-02\\
60&8.00E-05&5.10E-03&5.28E-05&1.07E-04&1.25E-02&1.02E-02&1.56E-02\\
70&5.87E-05&4.02E-03&3.88E-05&7.85E-05&1.07E-02&8.77E-03&1.34E-02\\
80&4.49E-05&3.28E-03&2.97E-05&6.01E-05&9.33E-03&7.68E-03&1.17E-02\\
90&3.55E-05&2.74E-03&2.35E-05&4.75E-05&8.29E-03&6.82E-03&1.04E-02\\
100&2.87E-05&2.33E-03&1.90E-05&3.85E-05&7.46E-03&6.14E-03&9.38E-03\\
\hline
\end{tabular*}
\end{table}
\begin{table}[H]
\centering
\caption{The convergence rates of first-order scheme with $\tau =h^{2}$ at $\T=0.1$.}
\begin{tabular*}{16cm}{@{\extracolsep{\fill}}c c c c c c c c c c}
\hline
$1/h$ & $u_{L2}$  & $p_{L2}$  & $\omega_{L2}$ & $T_{L2}$ & $u_{H1}$ & $\omega_{H1}$ & $T_{H1}$ & $RPC(s)$& $FEM(s)$\\
\hline
10&-&-&-&-&-&-&-&0.27&0.27\\
20&1.99 &1.53 &1.97 &1.95 &1.01 &0.99 &0.97 &1.04 &1.34\\
30&2.01 &1.61 &1.99 &1.98 &1.01 &1.00 &0.99 &6.25 &8.92\\
40&2.01 &1.58 &2.00 &1.99 &1.01 &1.00 &1.00 &23.80 &40.66\\
50&2.01 &1.56 &2.00 &2.00 &1.01 &1.00 &1.00 &66.63 &110.92\\
60&2.01 &1.55 &2.00 &2.00 &1.01 &1.00 &1.00 &169.18 &293.29\\
70&2.01 &1.54 &2.00 &2.00 &1.00 &1.00 &1.00 &375.68 &634.49\\
80&2.01 &1.53 &2.00 &2.00 &1.00 &1.00 &1.00 &750.19 &1390.34\\
90&2.01 &1.53 &2.00 &2.00 &1.00 &1.00 &1.00 &1464.60 &2494.04\\
100&2.01 &1.53 &2.00 &2.00 &1.00 &1.00 &1.00 &2661.36 &4454.48\\
\hline
\end{tabular*}
\end{table}
\subsubsection{Second-order RPC scheme of 2D}
For the second-order format, select $\tau =h$ as the time step, and use the $P_{1}b-P_{1}-P_{1}-P_{1}$ element to approximate the finite element space. Therefore, the optimal relative error order of $L^2$ in the first-order format $u$, $\omega$ and $T$ is $O(h^{2})$, and $H^1$ relative error optimal order is $O(h^{1})$.
\begin{table}[H]
\centering
\caption{The error of second-order scheme with $\tau =h$ at $\T=0.1$.}
\begin{tabular*}{16cm}{@{\extracolsep{\fill}}c c c c c c c c c}
\hline
$1/h$ & $\|u-u_h\|_0$  & $\|p-p_h\|_0$  & $\|\omega-\omega_h\|_0$ & $\|T-T_h\|_0$ & $|u-u_h|_1$ & $|\omega-\omega_h|_1$ & $|T-T_h|_1$\\
\hline
10&3.64E-03&9.22E-02&1.86E-03&3.50E-03&7.73E-02&6.09E-02&9.16E-02\\
20&1.01E-03&4.85E-02&4.94E-04&9.42E-04&3.85E-02&3.07E-02&4.66E-02\\
30&4.54E-04&3.29E-02&2.25E-04&4.25E-04&2.56E-02&2.05E-02&3.12E-02\\
40&2.55E-04&2.48E-02&1.28E-04&2.40E-04&1.92E-02&1.53E-02&2.34E-02\\
50&1.63E-04&1.98E-02&8.22E-05&1.54E-04&1.53E-02&1.23E-02&1.87E-02\\
60&1.13E-04&1.64E-02&5.72E-05&1.07E-04&1.28E-02&1.02E-02&1.56E-02\\
70&8.28E-05&1.40E-02&4.21E-05&7.85E-05&1.09E-02&8.77E-03&1.34E-02\\
80&6.33E-05&1.22E-02&3.22E-05&6.01E-05&9.55E-03&7.68E-03&1.17E-02\\
90&5.00E-05&1.08E-02&2.54E-05&4.75E-05&8.49E-03&6.82E-03&1.04E-02\\
100&4.04E-05&9.72E-03&2.06E-05&3.85E-05&7.64E-03&6.14E-03&9.38E-03\\
\hline
\end{tabular*}
\end{table}
\begin{table}[H]
\centering
\caption{The convergence rates of second-order scheme with $\tau =h$ at $\T=0.1$.}
\begin{tabular*}{16cm}{@{\extracolsep{\fill}}c c c c c c c c c}
\hline
$1/h$ & $u_{L2}$  & $p_{L2}$  & $\omega_{L2}$ & $T_{L2}$ & $u_{H1}$ & $\omega_{H1}$ & $T_{H1}$ & $time(s)$\\
\hline
10&-&-&-&-&-&-&-&0.07\\
20&1.86 &0.93 &1.91 &1.90 &1.01 &0.99 &0.97 &0.12\\
30&1.97 &0.96 &1.94 &1.96 &1.01 &1.00 &0.99 &0.38\\
40&2.00 &0.99 &1.96 &1.99 &1.01 &1.00 &1.00 &0.96\\
50&2.01 &1.01 &1.98 &2.00 &1.00 &1.00 &1.00 &1.98\\
60&2.01 &1.02 &1.99 &2.00 &1.00 &1.00 &1.00 &4.02\\
70&2.01 &1.03 &1.99 &2.00 &1.00 &1.00 &1.00 &6.88\\
80&2.01 &1.03 &2.00 &2.00 &1.00 &1.00 &1.00 &11.72\\
90&2.01 &1.03 &2.00 &2.00 &1.00 &1.00 &1.00 &19.21\\
100&2.01 &1.02 &2.01 &2.00 &1.00 &1.00 &1.00 &30.72\\
\hline
\end{tabular*}
\end{table}
\subsubsection{First-order RPC scheme of 3D}
In the case of 3D, the following smooth analytical solution problem is considered, $\nu=\nu_r=\hat{e}=\a=\b=\k=D=1, \T=0.5$, and the initial value in (\ref{1}) is set to the exact solution.
\begin{eqnarray}
\begin{array}{lll}
\ (u_1, u_2, u_3) =  ((y^4+z^2)\cos(t), (z^4+x^2)\cos(t), (x^4+y^2)\cos(t)),\\
\ (\omega_1, \omega_2, \omega_3) =  ((\sin(y)+z)\cos(t), (\sin(z)+x)\cos(t), (\sin(x)+y)\cos(t)),\\
\ p = (2x-1)(2y-1)(2z-1)\cos(t), \quad \ T = u_1+u_2+u_3.\\
\end{array}
\end{eqnarray}

For the first-order format, select $\tau =h^{2}$ as the time step, and use the $P_{1}b-P_{1}-P_{1}-P_{1}$ element to approximate the finite element space. Therefore, the optimal relative error order of $L^2$ in the first-order format $u$, $\omega$ and $T$ is $O(h^{2})$, and $H^1$ relative error optimal order is $O(h^{1})$.
\begin{table}[H]
\centering
\caption{The error of first-order RPC scheme with 3D $\tau =h^{2}$ at $\T=0.5$.}
\begin{tabular*}{16cm}{@{\extracolsep{\fill}}c c c c c c c c c}
\hline
$1/h$ & $\|u-u_h\|_0$  & $\|p-p_h\|_0$  & $\|\omega-\omega_h\|_0$ & $\|T-T_h\|_0$ & $|u-u_h|_1$ & $|\omega-\omega_h|_1$ & $|T-T_h|_1$\\
\hline
4&5.77E-02&5.15E-01&7.94E-03&3.46E-02&5.97E-01&6.35E-02&7.58E-01\\
6&2.59E-02&2.76E-01&3.77E-03&1.55E-02&4.01E-01&4.04E-02&5.08E-01\\
8&1.46E-02&1.76E-01&2.18E-03&8.76E-03&3.01E-01&2.96E-02&3.82E-01\\
10&9.34E-03&1.25E-01&1.41E-03&5.62E-03&2.41E-01&2.34E-02&3.06E-01\\
12&6.49E-03&9.40E-02&9.84E-04&3.90E-03&2.01E-01&1.94E-02&2.55E-01\\
14&4.77E-03&7.42E-02&7.26E-04&2.87E-03&1.72E-01&1.66E-02&2.19E-01\\
16&3.65E-03&6.06E-02&5.57E-04&2.20E-03&1.51E-01&1.45E-02&1.91E-01\\
\hline
\end{tabular*}
\end{table}
\begin{table}[H]
\centering
\caption{The convergence rates of first-order RPC scheme with 3D $\tau =h^{2}$ at $\T=0.5$.}
\begin{tabular*}{16cm}{@{\extracolsep{\fill}}c c c c c c c c c c}
\hline
$1/h$ & $u_{L2}$  & $p_{L2}$  & $\omega_{L2}$ & $T_{L2}$ & $u_{H1}$ & $\omega_{H1}$ & $T_{H1}$ & $RPC(s)$& $FEM(s)$\\
\hline
4&-&-&-&-&-&-&-&0.52&0.88\\
6&1.98 &1.54 &1.83 &1.98 &0.98 &1.12 &0.98 &2.94 &4.15\\
8&1.99 &1.56 &1.91 &1.99 &0.99 &1.08 &0.99 &13.44 &19.09\\
10&2.00 &1.56 &1.95 &1.99 &1.00 &1.05 &1.00 &54.74 &68.90\\
12&2.00 &1.54 &1.97 &2.00 &1.00 &1.04 &1.00 &147.85 &197.84\\
14&2.00 &1.53 &1.98 &2.00 &1.00 &1.03 &1.00 &390.21 &510.42\\
16&2.00 &1.52 &1.98 &2.00 &1.00 &1.02 &1.00 &880.55 &1186.04\\
\hline
\end{tabular*}
\end{table}
\subsubsection{Second-order RPC scheme of 3D}
For the second-order format, select $\tau =h$ as the time step, and use the $P_{1}b-P_{1}-P_{1}-P_{1}$ element to approximate the finite element space. Therefore, the optimal relative error order of $L^2$ in the first-order format $u$, $\omega$ and $T$ is $O(h^{2})$, and $H^1$ relative error optimal order is $O(h^{1})$.
\begin{table}[H]
\centering
\caption{The error of second-order RPC scheme with 3D $\tau =h$ at $\T=0.5$.}
\begin{tabular*}{16cm}{@{\extracolsep{\fill}}l l l l l l l l l}
\hline
$1/h$ & $\|u-u_h\|_0$  & $\|p-p_h\|_0$  & $\|\omega-\omega_h\|_0$ & $\|T-T_h\|_0$ & $|u-u_h|_1$ & $|\omega-\omega_h|_1$ & $|T-T_h|_1$\\
\hline
4&5.77E-02&4.54E-01&7.44E-03&3.46E-02&5.97E-01&6.20E-02&7.58E-01\\
6&2.58E-02&2.55E-01&3.39E-03&1.55E-02&4.01E-01&3.96E-02&5.08E-01\\
8&1.46E-02&1.67E-01&1.92E-03&8.76E-03&3.01E-01&2.93E-02&3.82E-01\\
10&9.31E-03&1.20E-01&1.23E-03&5.62E-03&2.41E-01&2.32E-02&3.06E-01\\
12&6.47E-03&9.19E-02&8.59E-04&3.90E-03&2.01E-01&1.93E-02&2.55E-01\\
14&4.75E-03&7.38E-02&6.31E-04&2.87E-03&1.72E-01&1.65E-02&2.19E-01\\
16&3.63E-03&6.14E-02&4.82E-04&2.20E-03&1.51E-01&1.44E-02&1.91E-01\\
\hline
\end{tabular*}
\end{table}
\begin{table}[H]
\centering
\caption{The convergence rates of second-order RPC scheme with 3D $\tau =h$ at $\T=0.5$.}
\begin{tabular*}{16cm}{@{\extracolsep{\fill}}l l l l l l l l l}
\hline
$1/h$ & $u_{L2}$  & $p_{L2}$  & $\omega_{L2}$ & $T_{L2}$ & $u_{H1}$ & $\omega_{H1}$ & $T_{H1}$ & $time(s)$\\
\hline
4&-&-&-&-&-&-&-&0.27\\
6&1.98 &1.42 &1.94 &1.98 &0.98 &1.10 &0.98 &1.18\\
8&1.99 &1.47 &1.97 &1.99 &0.99 &1.06 &0.99 &3.39\\
10&2.00 &1.48 &1.98 &1.99 &1.00 &1.03 &1.00 &10.16\\
12&2.00 &1.46 &1.99 &2.00 &1.00 &1.02 &1.00 &21.94\\
14&2.00 &1.42 &2.00 &2.00 &1.00 &1.02 &1.00 &44.66\\
16&2.00 &1.38 &2.01 &2.00 &1.00 &1.01 &1.00 &83.51\\
\hline
\end{tabular*}
\end{table}
It can be seen that the first-order and second-order RPC methods proposed in this paper have reached the optimal order of theoretical analysis in the calculation of the above smooth analytical solutions, and for the second-order methods of other articles, the second-order RPC methods proposed in this paper allows to obtain better results with a coarser space step and time step.
\begin{figure}[H]
\centering
\subfigure{\includegraphics[width=4cm]{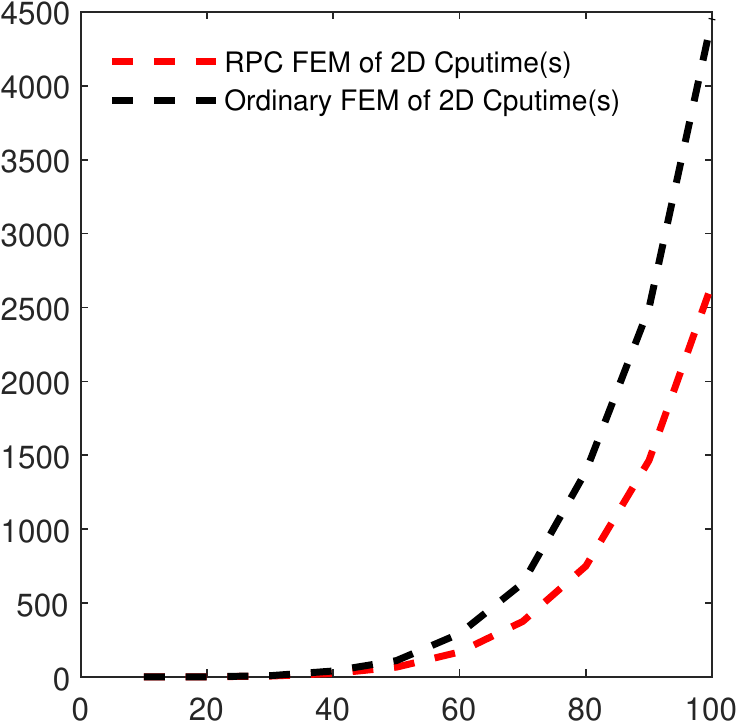}}
\quad
\subfigure{\includegraphics[width=4cm]{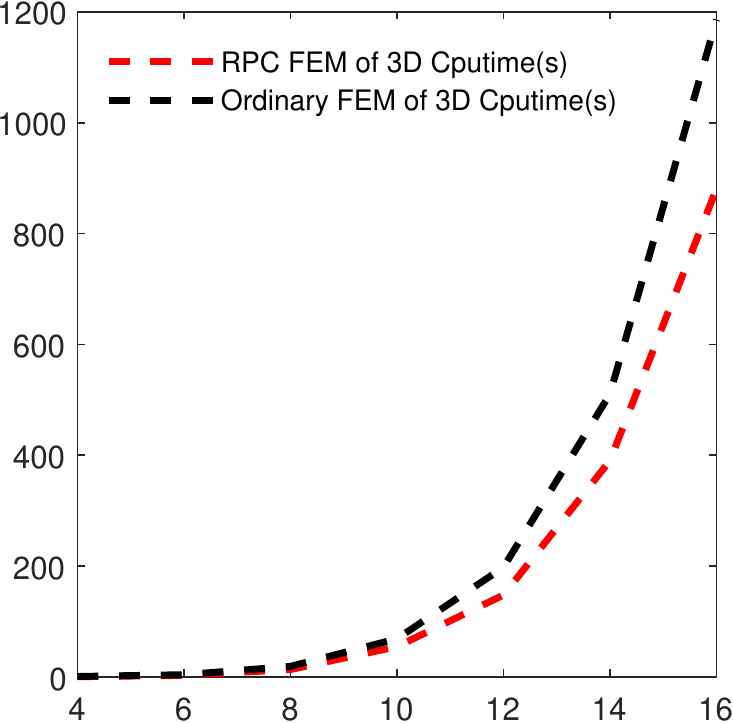}}
\caption{RPC FEM and Ordinary FEM cpu time-consuming comparison}
\label{Figure time}
\end{figure}
We can compare that RPC method is obviously faster than ordinary finite element method, and the advantages of PRC method become more and more obvious with the increase of time layer.
\subsection{Benard convection problem}
In this section, we consider the Benard convection problem, we consider the domain $\Omega=[0, 5]\times[0, 1]$ and given $f_1, f_2, f_3=0$, $\hat{e}=10^4, 10^5, 10^6$, $\nu=\nu_r=\a=\b=\k=D=1$, $h$=1/20, $\tau=h^2$, the boundary conditions are given below, and the initial values are the same as the boundary.
\begin{eqnarray}
\begin{array}{l}
T(x,0,t)=1, \quad T(x,1,t)=0, \quad \frac{\partial T}{\partial n}=0\quad\mathrm{when}\quad x=0, 5,\\
u|_{\partial \Omega}=0, \quad \omega |_{\partial \Omega}=0, \quad f_1=0, \quad f_2=0, \quad f_3=0. \no
\end{array}
\end{eqnarray}
\begin{figure}[H]
\centering
\subfigure[$\T$=1, temperature]{\includegraphics[scale=0.12]{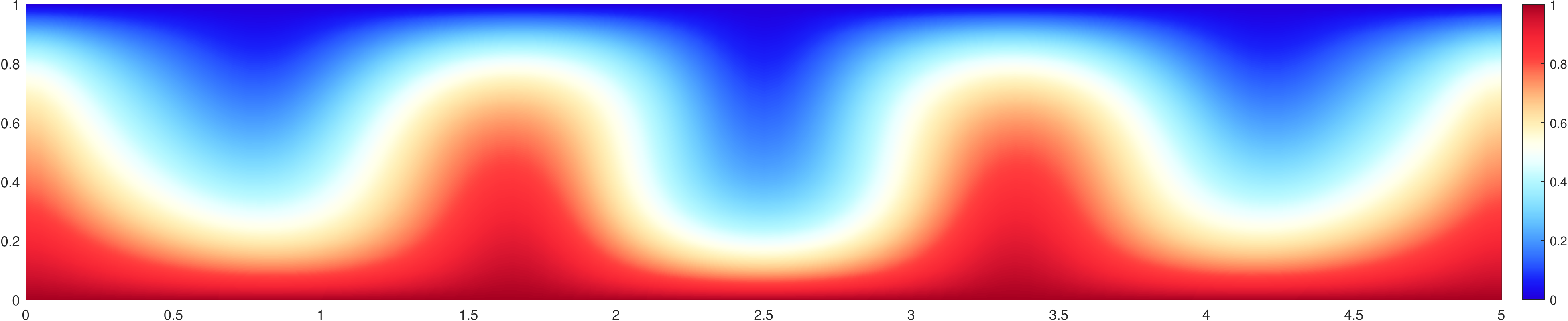}}
\quad
\subfigure[$\T$=1.5, temperature]{\includegraphics[scale=0.12]{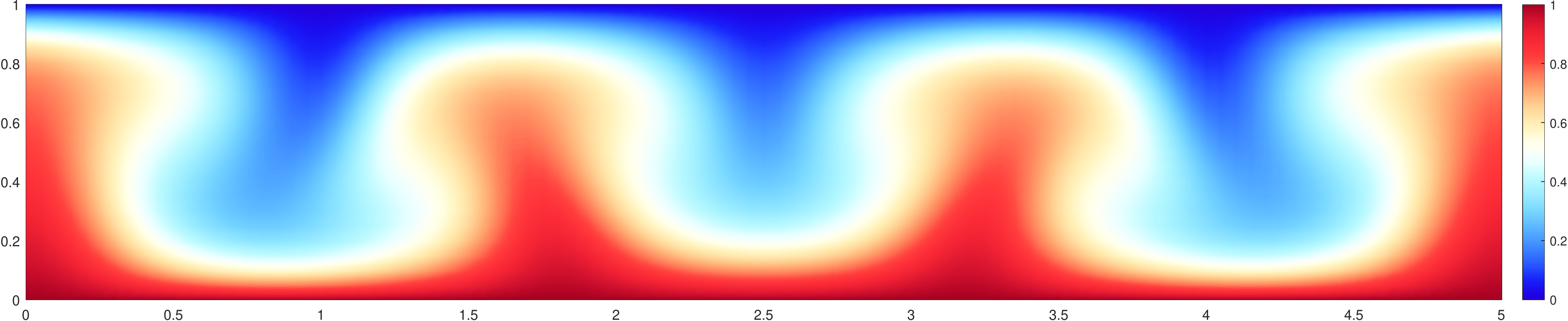}}
\quad
\subfigure[$\T$=2, temperature]{\includegraphics[scale=0.12]{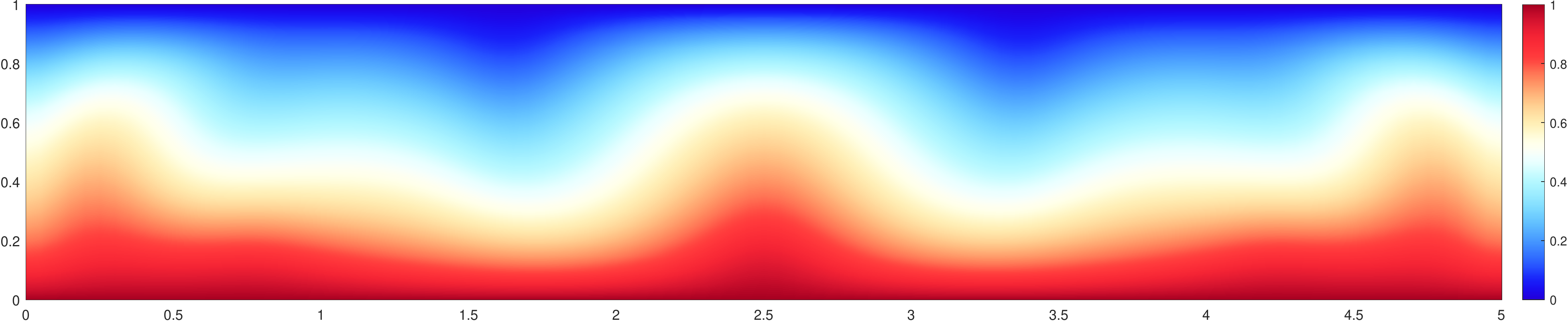}}
\end{figure}
\vspace{-8mm}
\begin{figure}[H]
\centering
\subfigure[$\T$=1, pressure]{\includegraphics[scale=0.12]{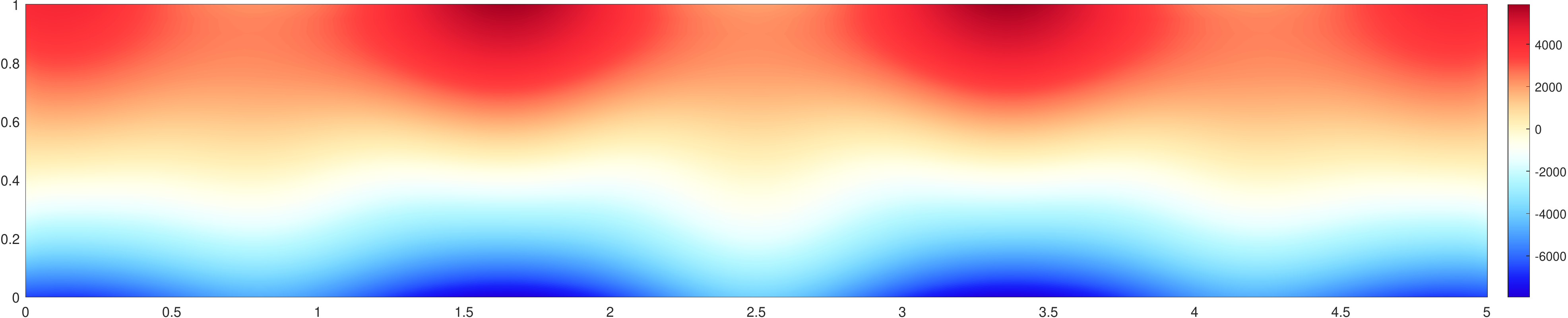}}
\quad
\subfigure[$\T$=1.5, pressure]{\includegraphics[scale=0.12]{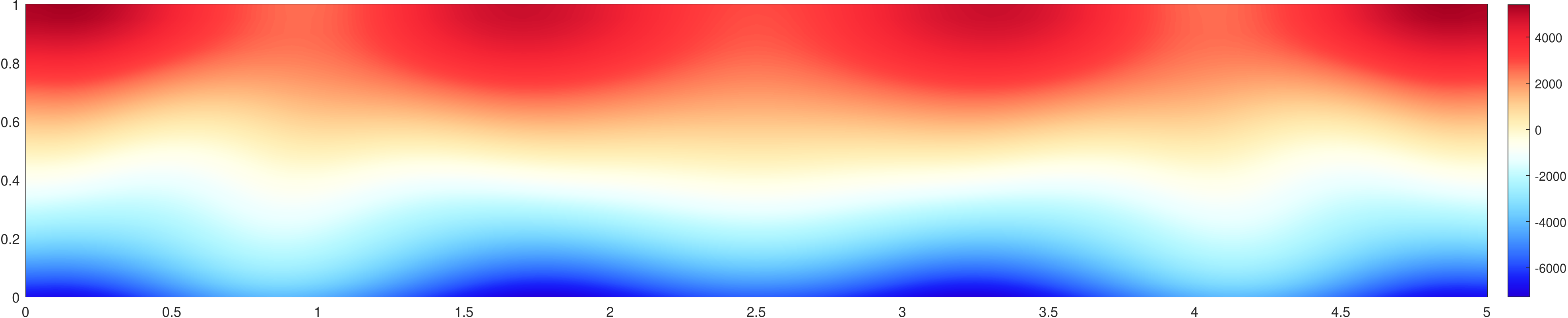}}
\quad
\subfigure[$\T$=2, pressure]{\includegraphics[scale=0.12]{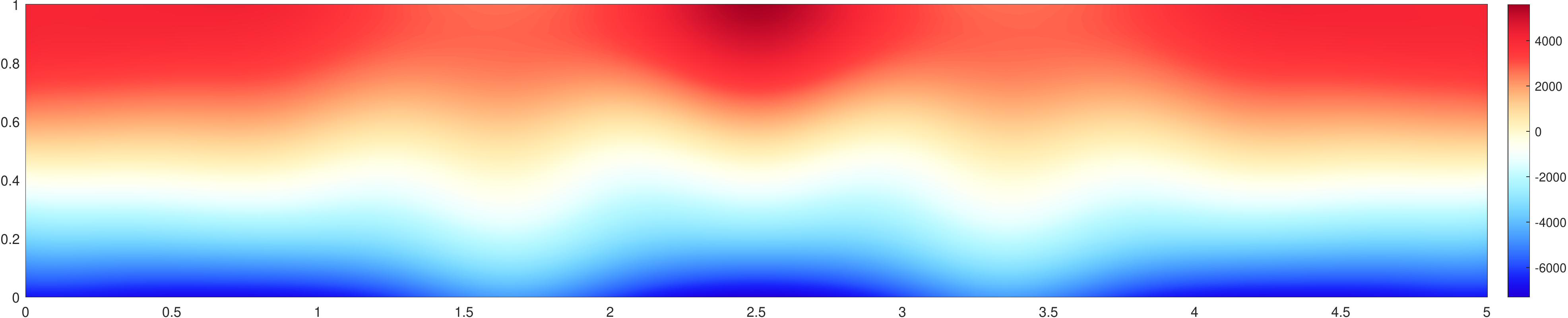}}
\end{figure}
\vspace{-8mm}
\begin{figure}[H]
\centering
\subfigure[$\T$=1, angular velocity]{\includegraphics[scale=0.12]{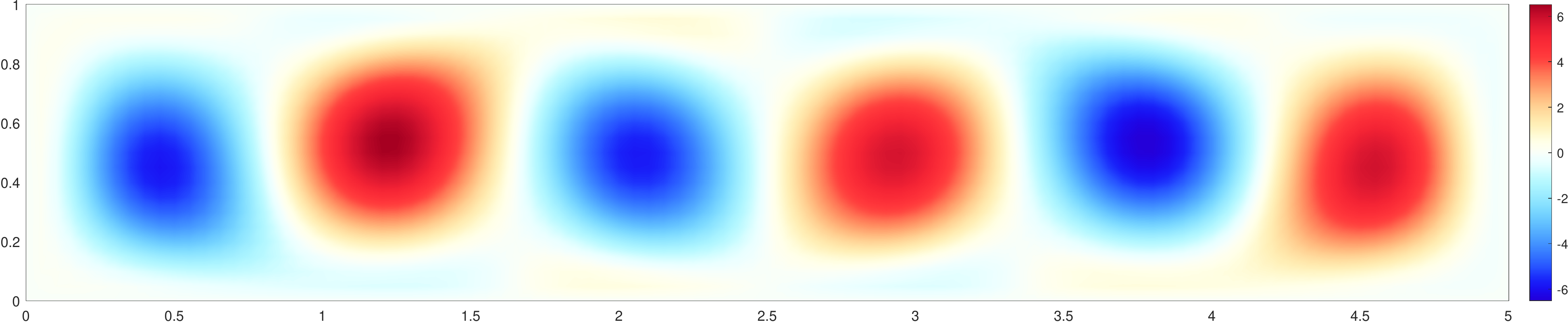}}
\quad
\subfigure[$\T$=1.5, angular velocity]{\includegraphics[scale=0.12]{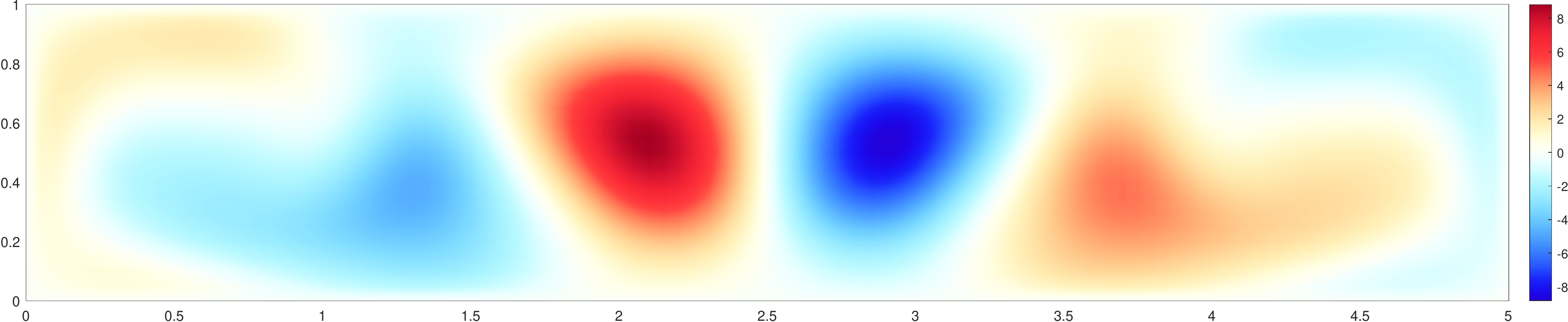}}
\quad
\subfigure[$\T$=2, angular velocity]{\includegraphics[scale=0.12]{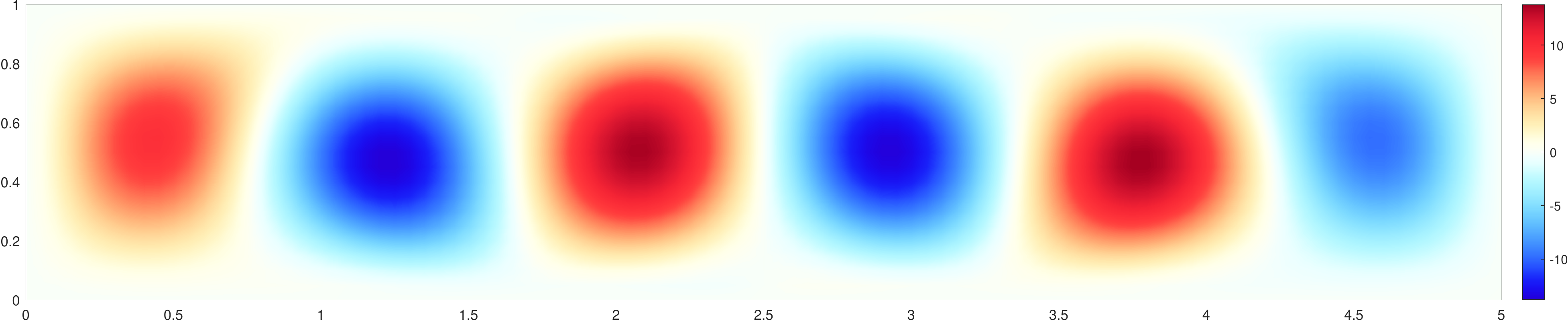}}
\end{figure}
\vspace{-8mm}
\begin{figure}[H]
\centering
\subfigure[$\T$=1, velocity]{\includegraphics[scale=0.12]{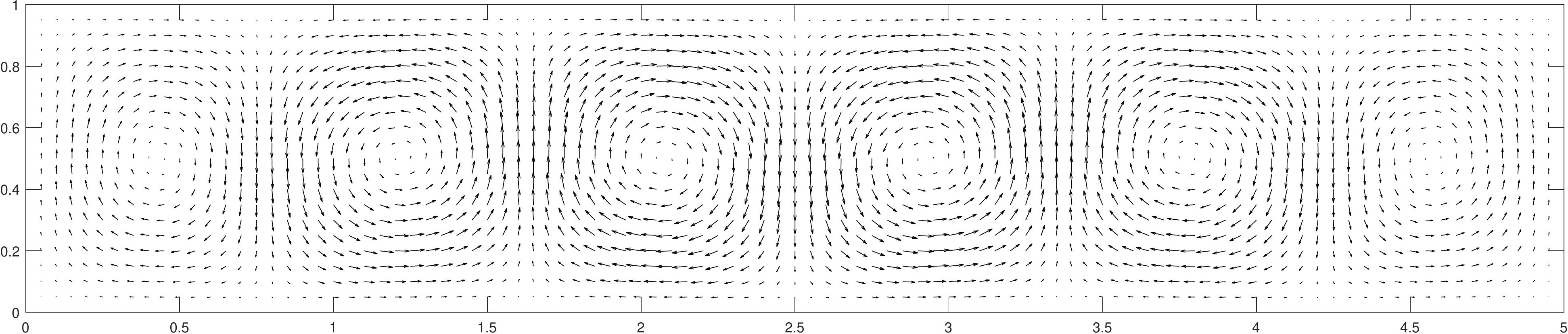}}
\quad
\subfigure[$\T$=1.5, velocity]{\includegraphics[scale=0.12]{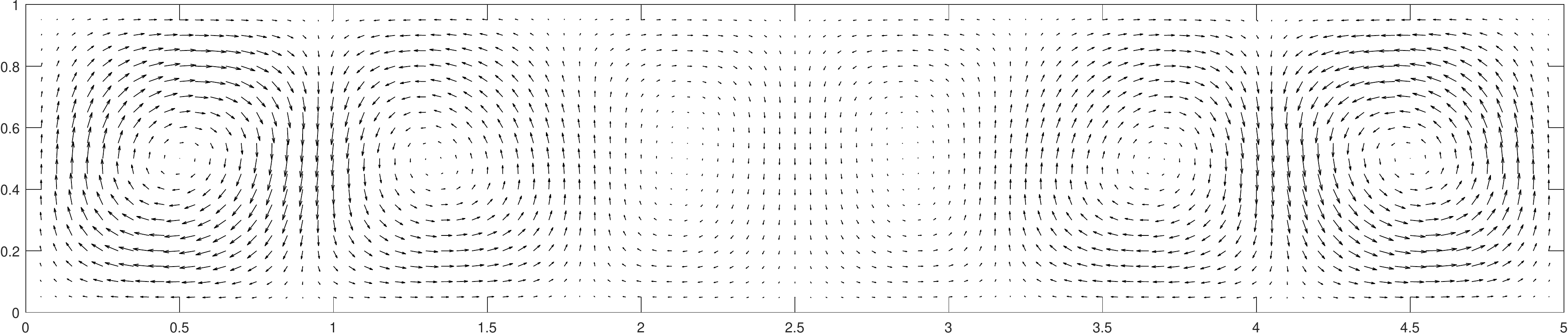}}
\quad
\subfigure[$\T$=2, velocity]{\includegraphics[scale=0.12]{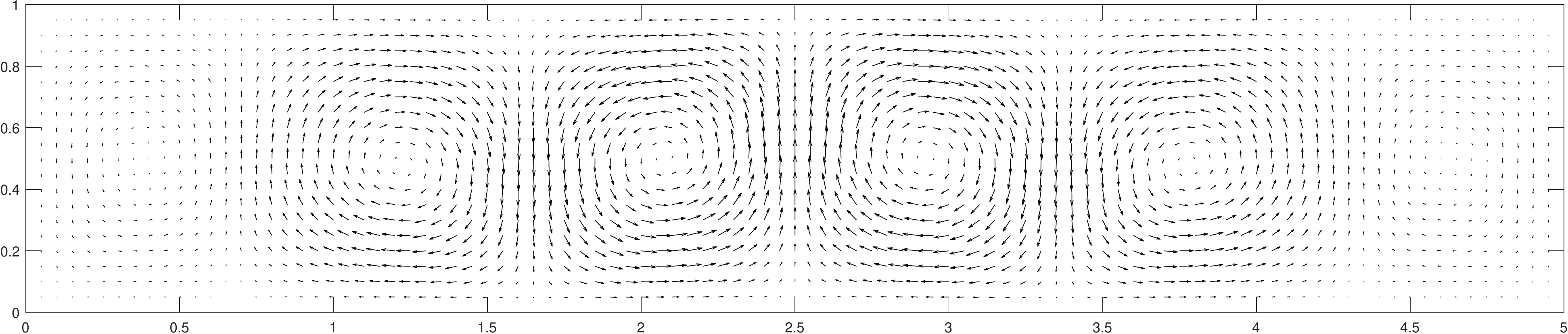}}
\caption{First-order RPC schemes, $P_2-P_1-P_2-P_2$, $\hat{e}=10^4$, $h$=1/20, $\tau=h^2$.}
\label{Figure 2.1}
\end{figure}
\begin{figure}[H]
\centering
\subfigure[$\T$=1, temperature]{\includegraphics[scale=0.12]{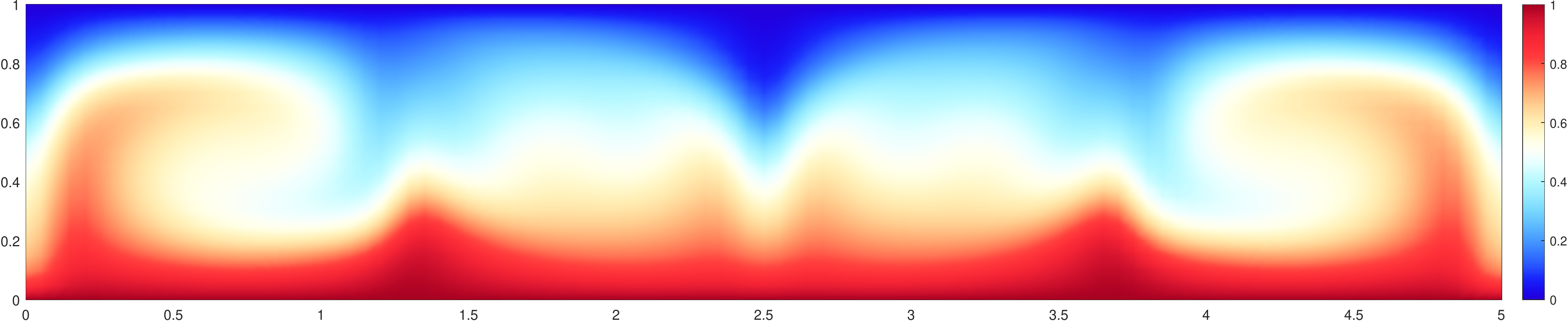}}
\quad
\subfigure[$\T$=1.5, temperature]{\includegraphics[scale=0.12]{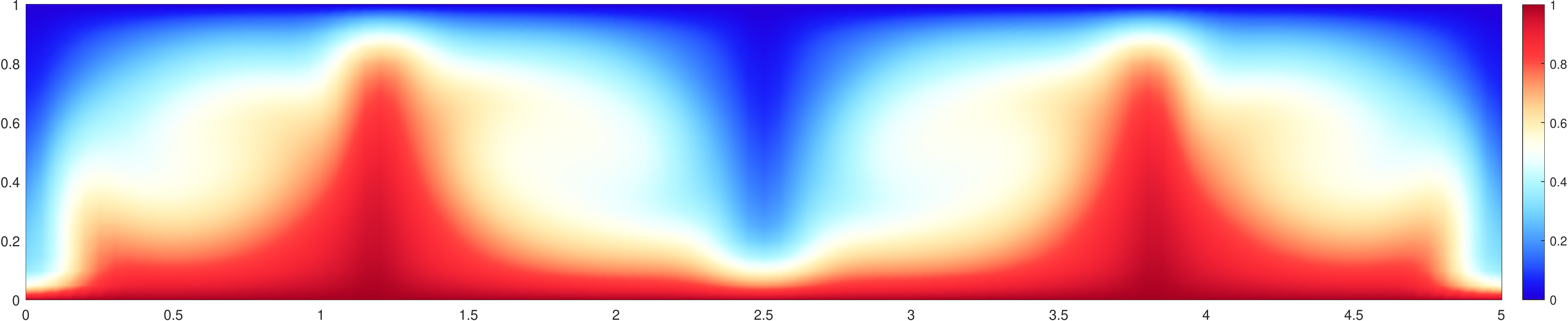}}
\quad
\subfigure[$\T$=2, temperature]{\includegraphics[scale=0.12]{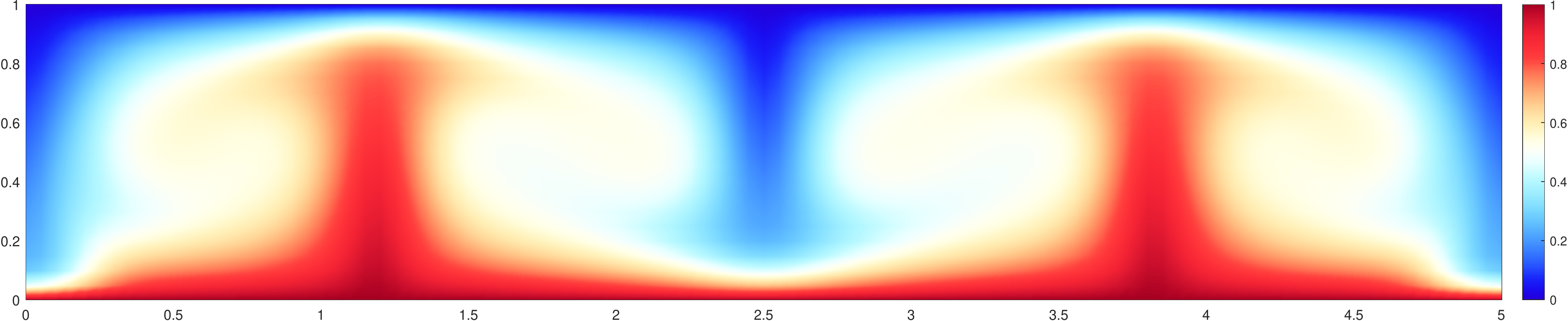}}
\end{figure}
\vspace{-8mm}
\begin{figure}[H]
\centering
\subfigure[$\T$=1, pressure]{\includegraphics[scale=0.12]{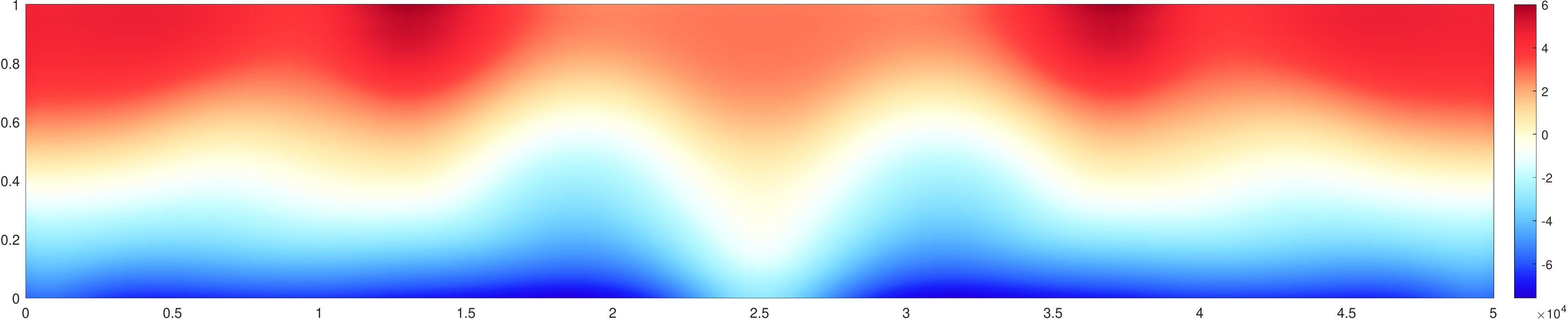}}
\quad
\subfigure[$\T$=1.5, pressure]{\includegraphics[scale=0.12]{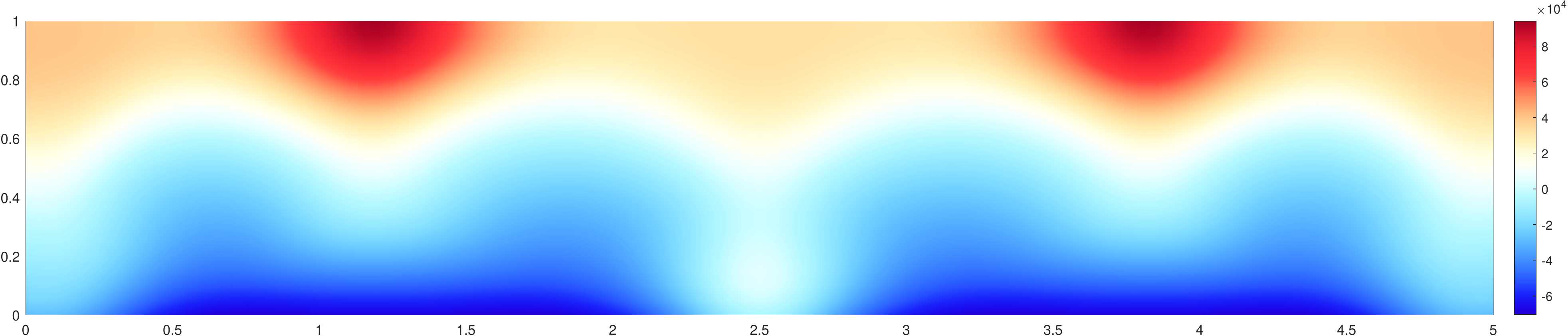}}
\quad
\subfigure[$\T$=2, pressure]{\includegraphics[scale=0.12]{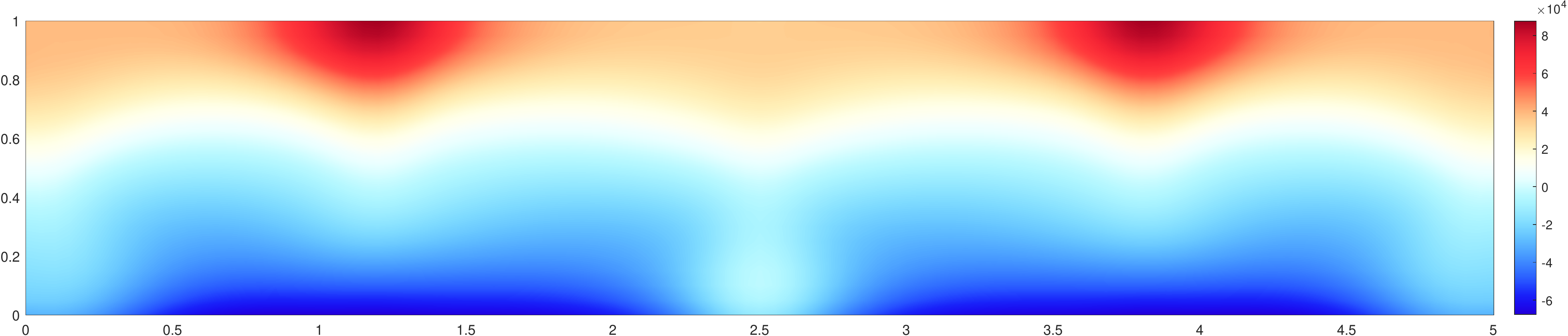}}
\end{figure}
\vspace{-8mm}
\begin{figure}[H]
\centering
\subfigure[$\T$=1, angular velocity]{\includegraphics[scale=0.12]{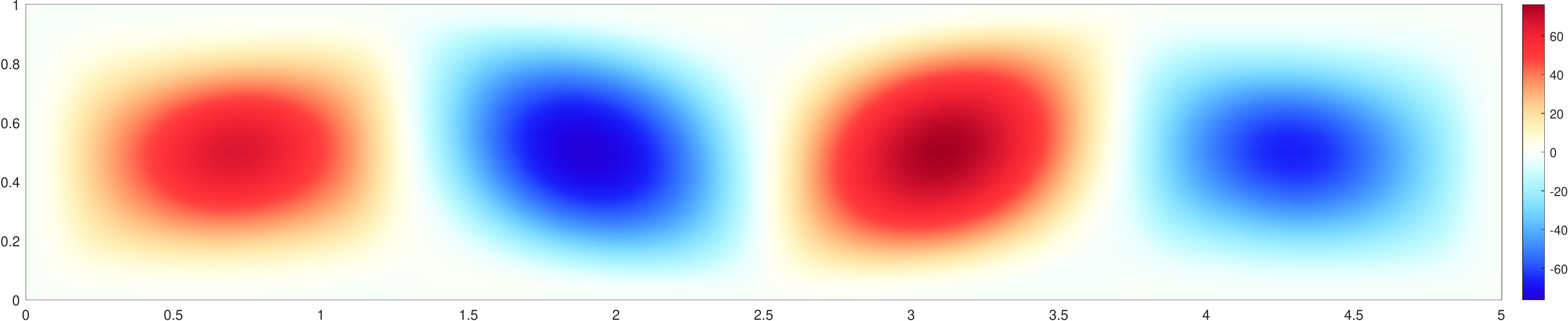}}
\quad
\subfigure[$\T$=1.5, angular velocity]{\includegraphics[scale=0.12]{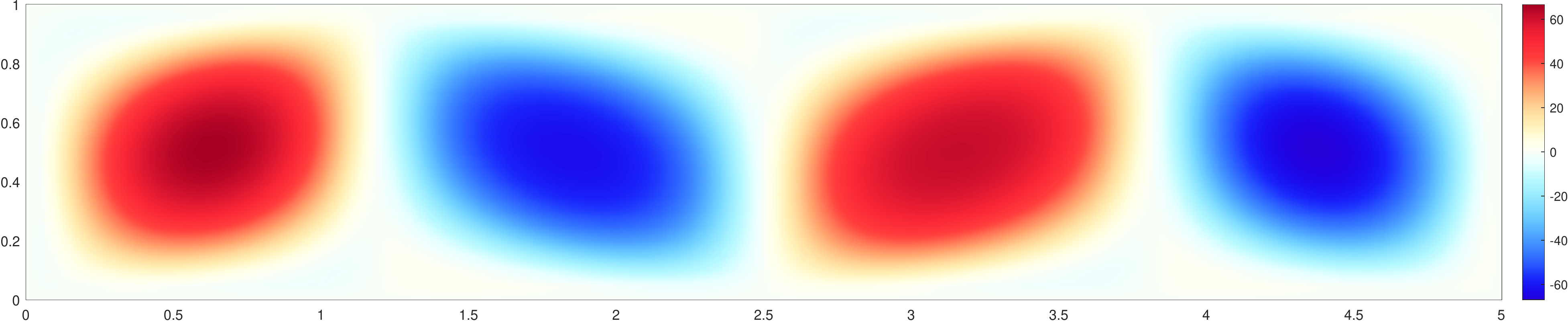}}
\quad
\subfigure[$\T$=2, angular velocity]{\includegraphics[scale=0.12]{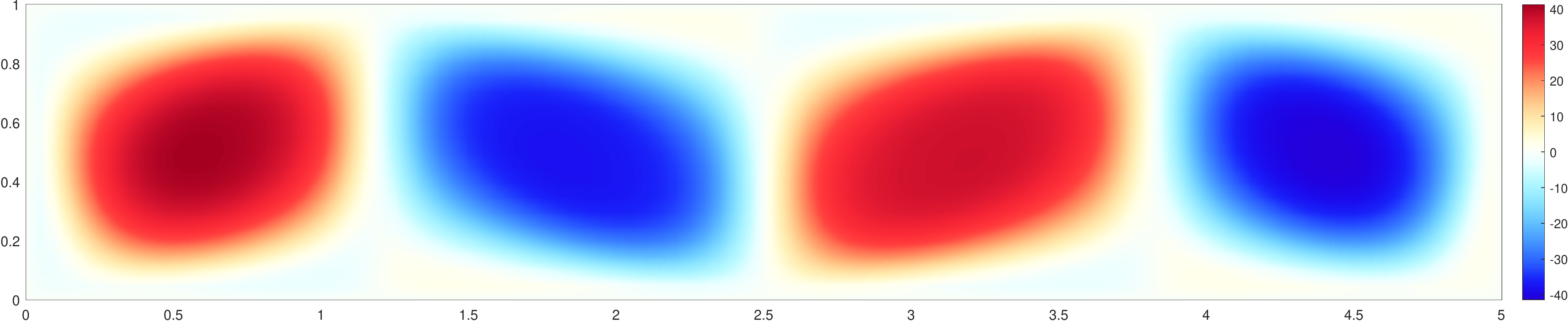}}
\end{figure}
\vspace{-8mm}
\begin{figure}[H]
\centering
\subfigure[$\T$=1, velocity]{\includegraphics[scale=0.12]{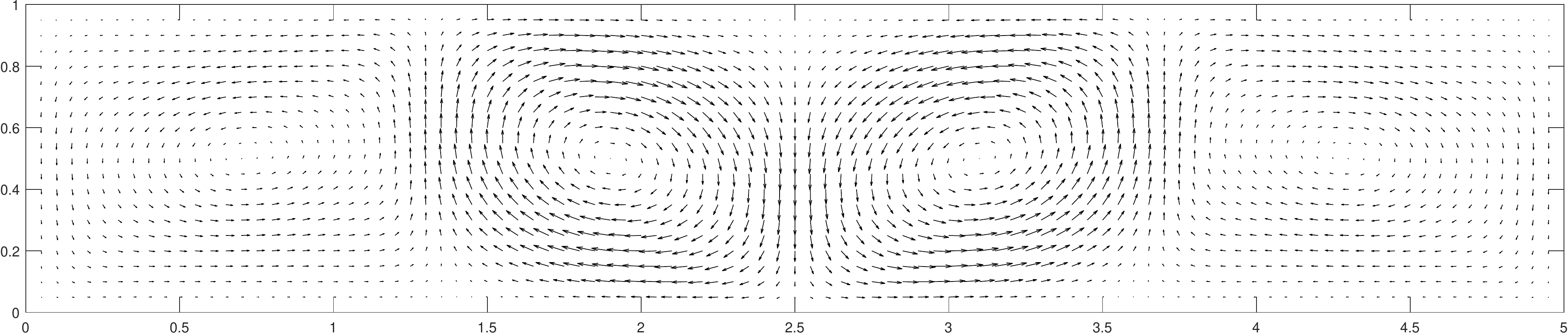}}
\quad
\subfigure[$\T$=1.5, velocity]{\includegraphics[scale=0.12]{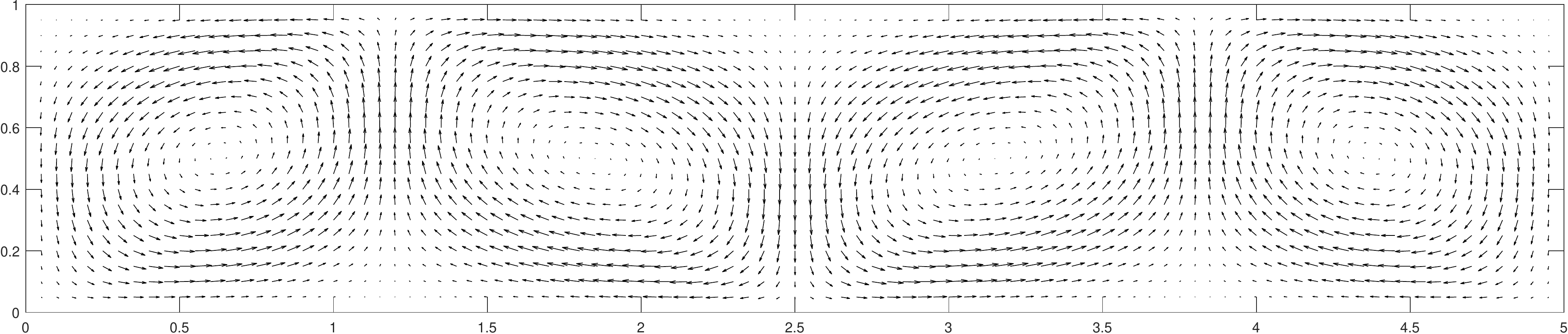}}
\quad
\subfigure[$\T$=2, velocity]{\includegraphics[scale=0.12]{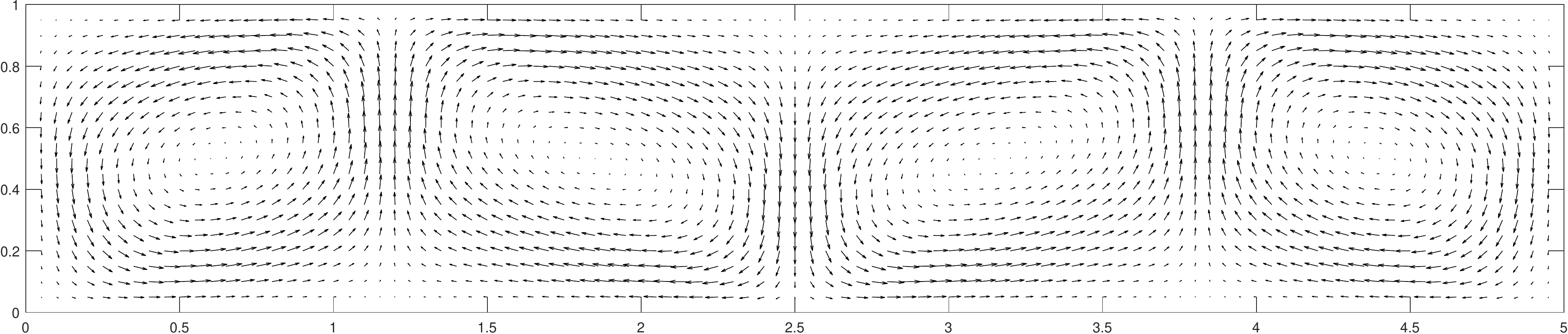}}
\caption{First-order RPC schemes, $P_2-P_1-P_2-P_2$, $\hat{e}=10^5$, $h$=1/20, $\tau=h^2$.}
\label{Figure 2.2}
\end{figure}
\begin{figure}[H]
\centering
\subfigure[$\T$=1, temperature]{\includegraphics[scale=0.12]{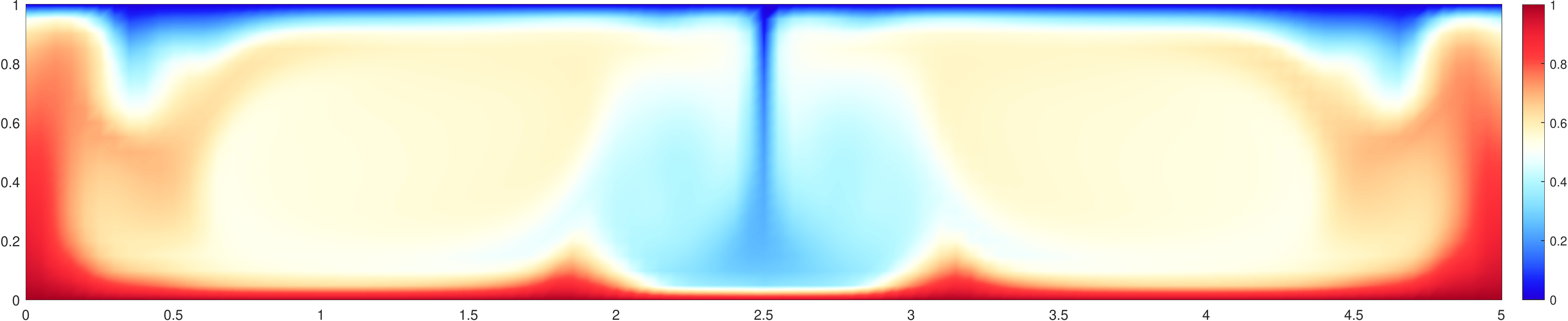}}
\quad
\subfigure[$\T$=1.5, temperature]{\includegraphics[scale=0.12]{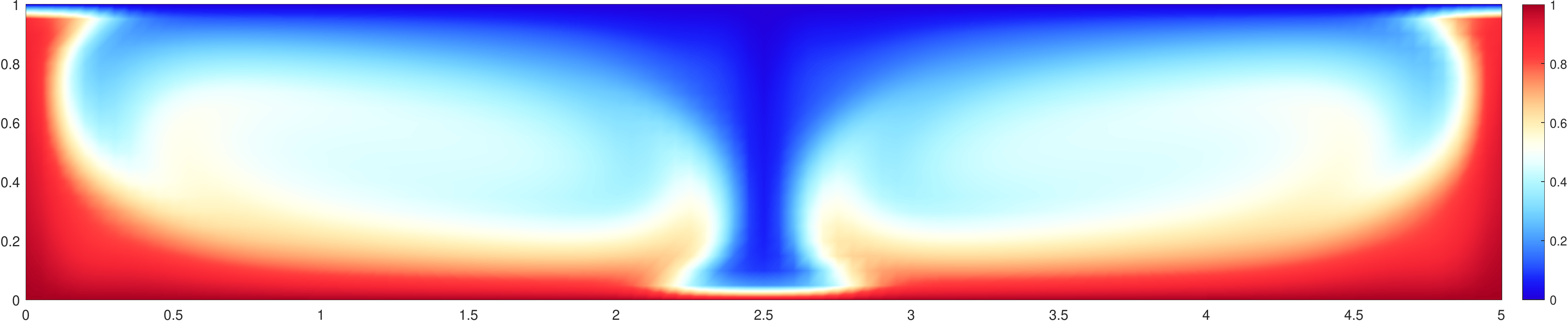}}
\quad
\subfigure[$\T$=2, temperature]{\includegraphics[scale=0.12]{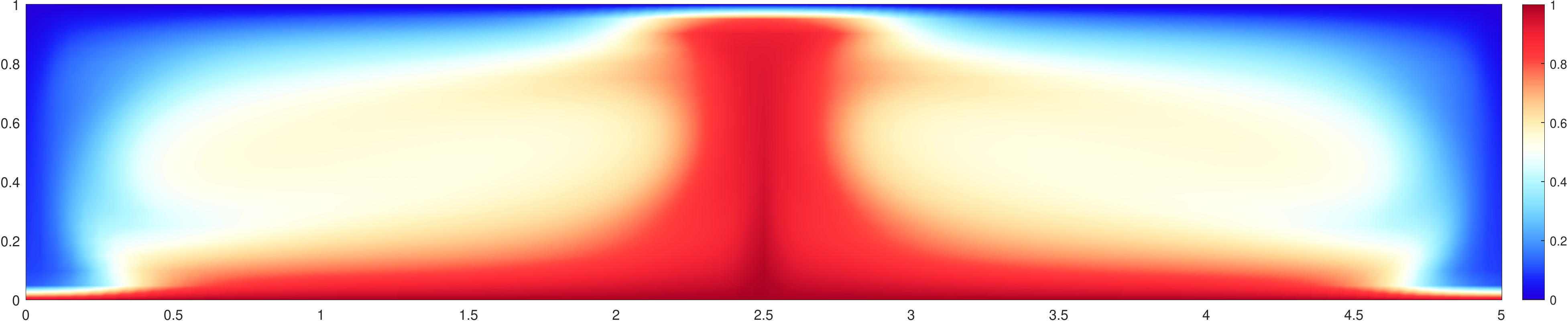}}
\end{figure}
\vspace{-8mm}
\begin{figure}[H]
\centering
\subfigure[$\T$=1, pressure]{\includegraphics[scale=0.12]{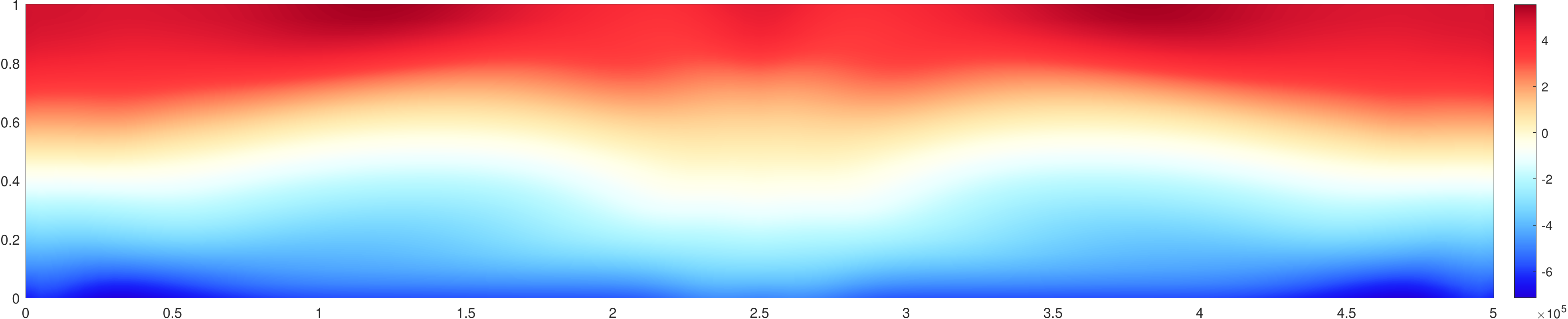}}
\quad
\subfigure[$\T$=1.5, pressure]{\includegraphics[scale=0.12]{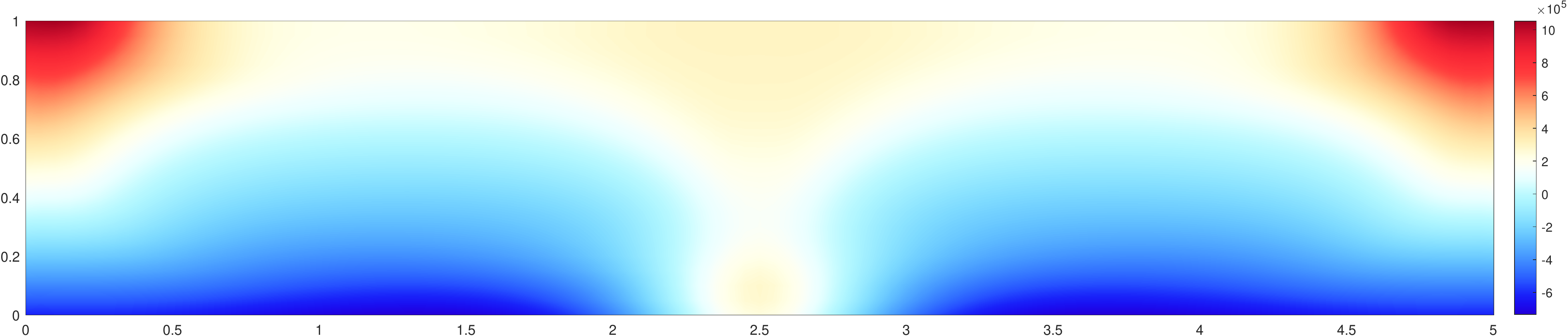}}
\quad
\subfigure[$\T$=2, pressure]{\includegraphics[scale=0.12]{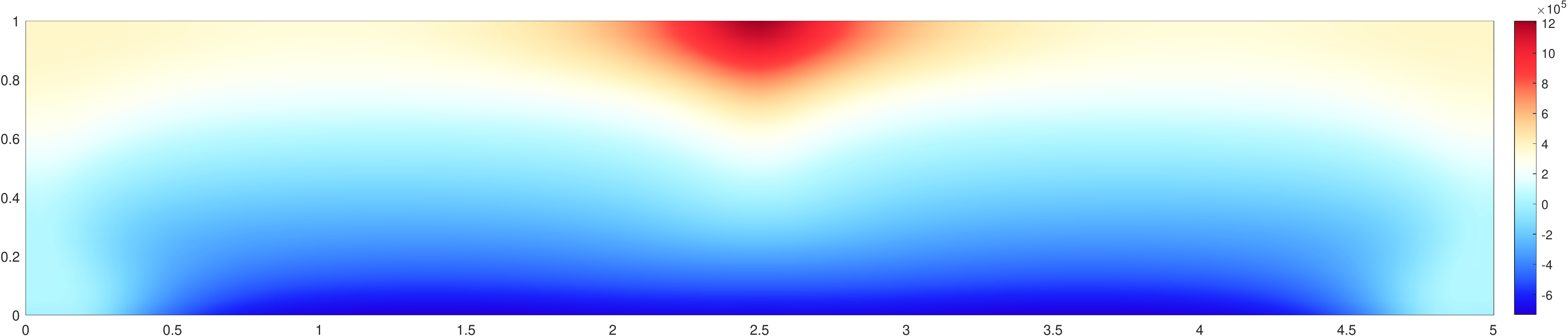}}
\end{figure}
\vspace{-8mm}
\begin{figure}[H]
\centering
\subfigure[$\T$=1, angular velocity]{\includegraphics[scale=0.12]{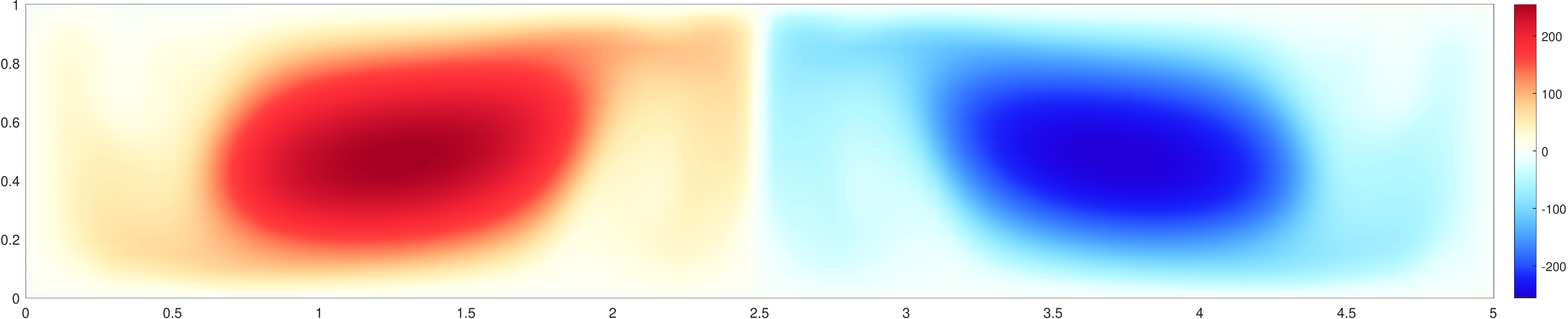}}
\quad
\subfigure[$\T$=1.5, angular velocity]{\includegraphics[scale=0.12]{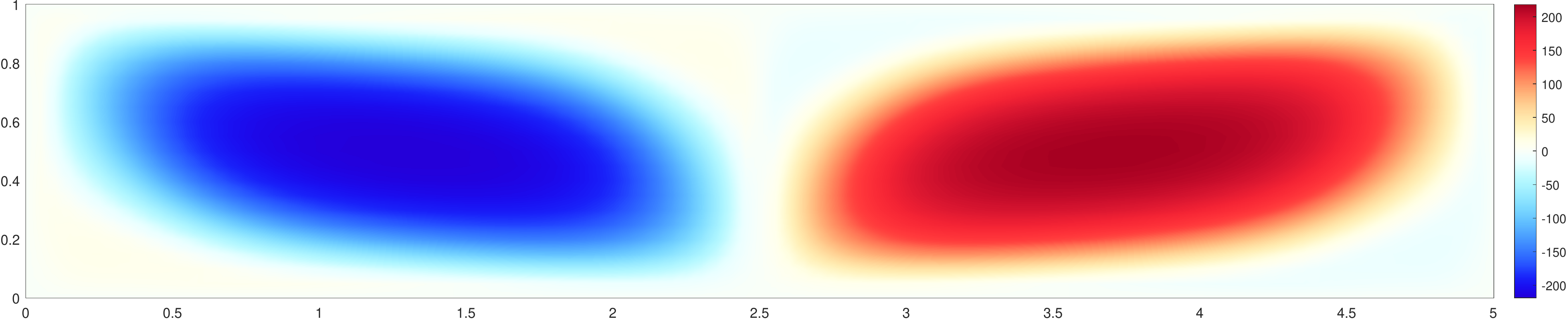}}
\quad
\subfigure[$\T$=2, angular velocity]{\includegraphics[scale=0.12]{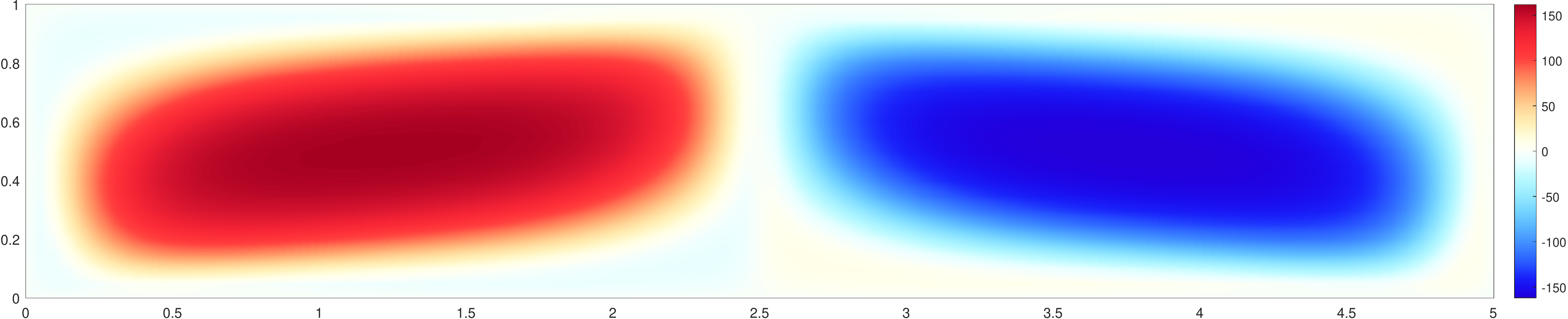}}
\end{figure}
\vspace{-8mm}
\begin{figure}[H]
\centering
\subfigure[$\T$=1, velocity]{\includegraphics[scale=0.12]{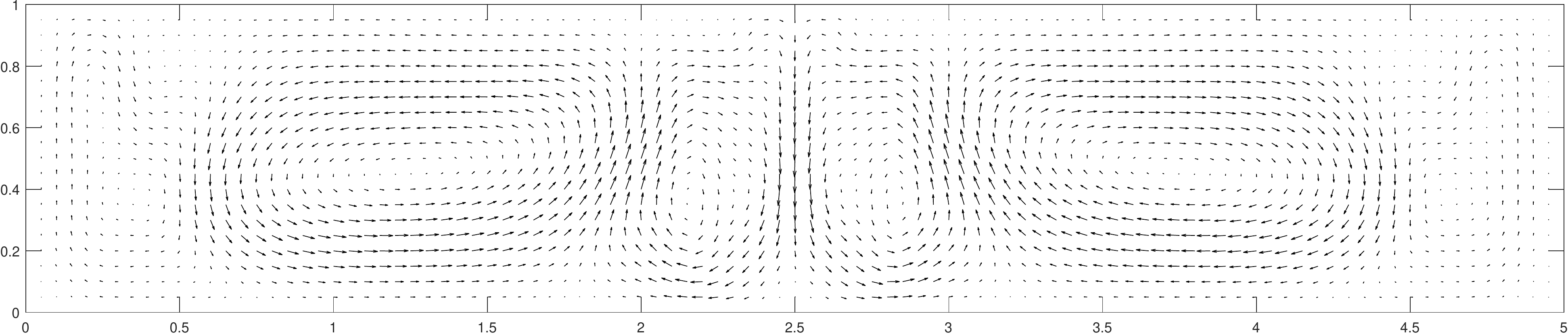}}
\quad
\subfigure[$\T$=1.5, velocity]{\includegraphics[scale=0.12]{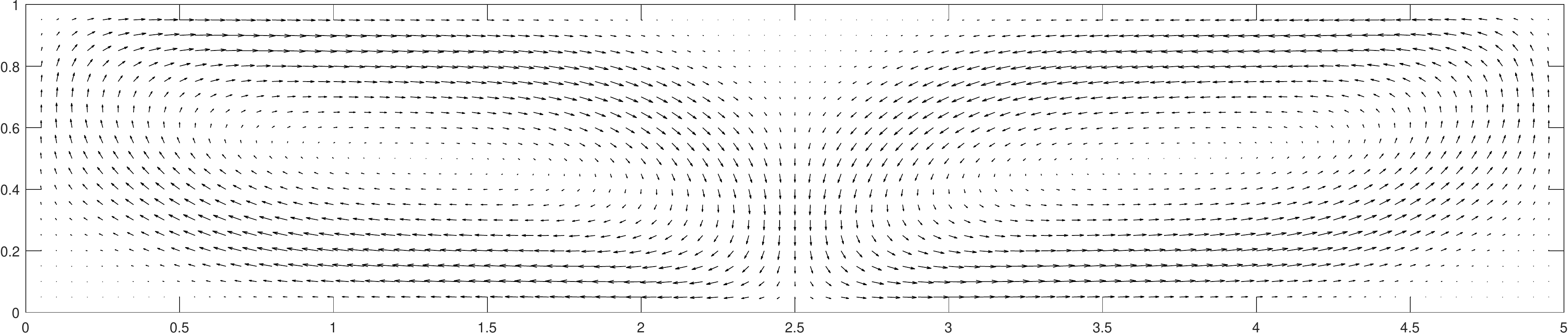}}
\quad
\subfigure[$\T$=2, velocity]{\includegraphics[scale=0.12]{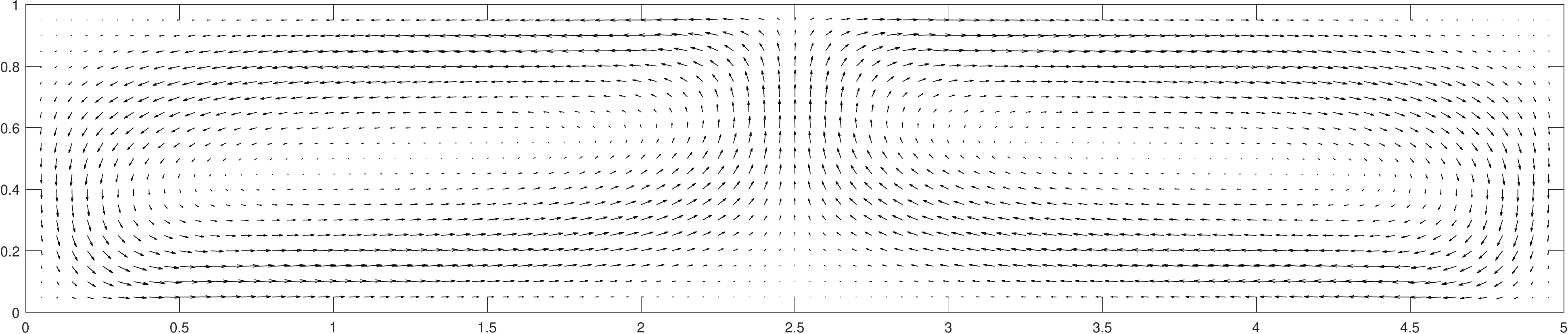}}
\caption{First-order RPC schemes, $P_2-P_1-P_2-P_2$, $\hat{e}=10^6$, $h$=1/20, $\tau=h^2$.}
\label{Figure 2.3}
\end{figure}
It can be seen that even for high Rayleigh numbers, the RPC method can still simulate the temperature, pressure, angular velocity and velocity well, see Figs.~\ref{Figure 2.1}-\ref{Figure 2.3}. Compared with the results in other papers, the results in this paper have better symmetry, which may be better in some ways~\cite{2014Modified, 2017Pressure}.
\subsection{Thermal driven cavity flow problem of 2D}
We consider the domain $\Omega=[0, 1]\times[0, 1]$ and given $f_1, f_2, f_3=0$, $\hat{e}=10^4, 10^5, 10^6$, $\nu=\nu_r=\a=\b=\k=D=1$, $h$=1/30, $\T=1$, $\tau=h^2$. The boundary conditions are given below, and the initial values are the same as the boundary.
\begin{eqnarray}
\begin{array}{l}
T(0,y,t)=1, \quad T(1,y,t)=0, \quad \frac{\partial T}{\partial n}=0\quad when\quad y=0, 1,\\
u|_{\partial \Omega}=0, \quad \omega |_{\partial \Omega}=0, \quad f_1=0, \quad f_2=0, \quad f_3=0. \no
\end{array}
\end{eqnarray}
\begin{figure}[H]
\centering
\subfigure[temperature, $\bar{e}=10^4$]{\includegraphics[scale=0.5]{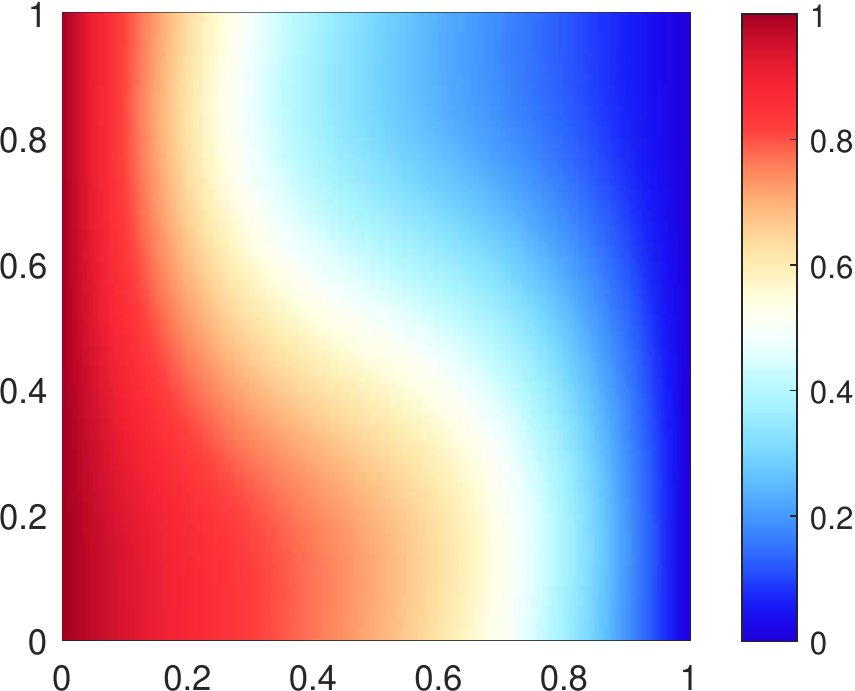}}
\quad
\subfigure[temperature, $\bar{e}=10^5$]{\includegraphics[scale=0.5]{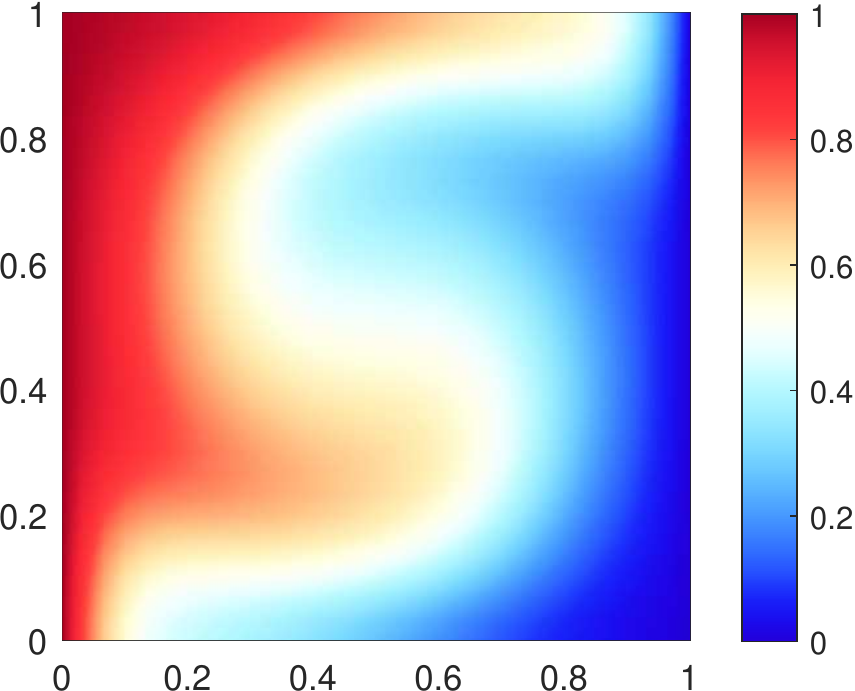}}
\quad
\subfigure[temperature, $\bar{e}=10^6$]{\includegraphics[scale=0.5]{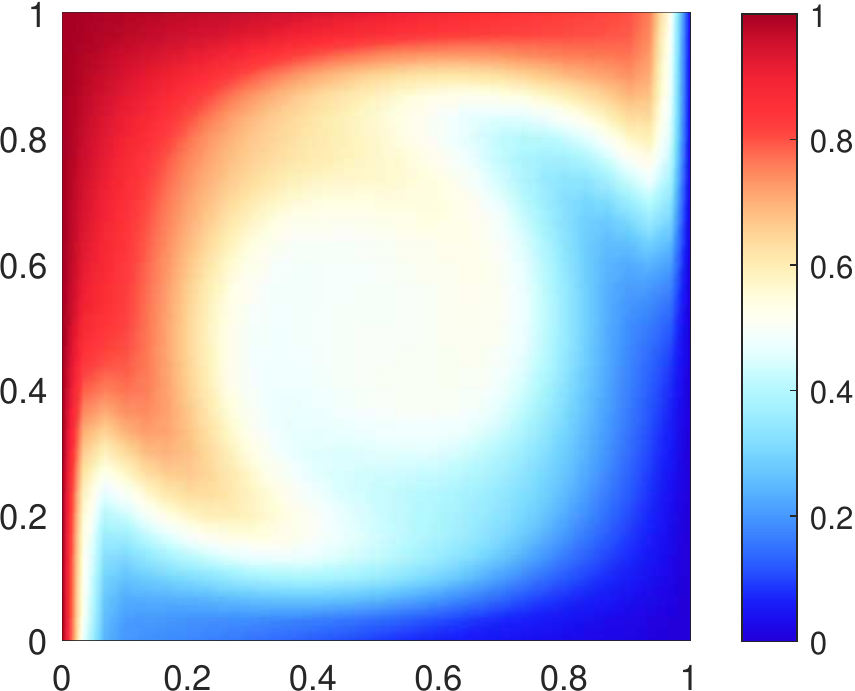}}
\caption{First-order RPC schemes, $P_2-P_1-P_2-P_2$, $\bar{e}=10^4, 10^5, 10^6$, $h$=1/30, $\tau=h^2$, $\T$=1.}
\label{Figure 3.1}
\end{figure}
\vspace{-8mm}
\begin{figure}[H]
\centering
\subfigure[pressure, $\bar{e}=10^4$]{\includegraphics[scale=0.5]{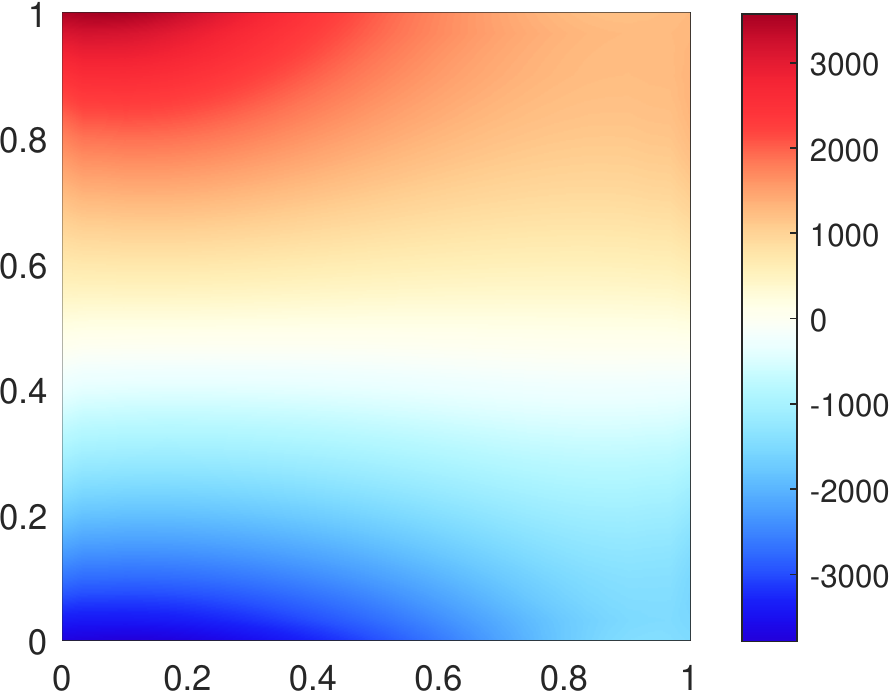}}
\quad
\subfigure[pressure, $\bar{e}=10^5$]{\includegraphics[scale=0.5]{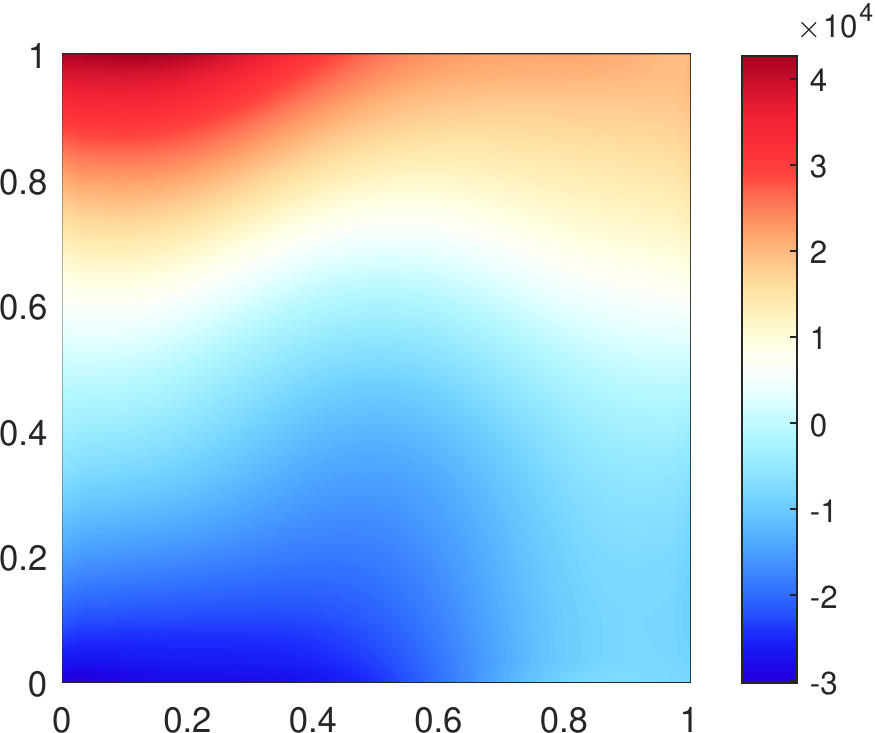}}
\quad
\subfigure[pressure, $\bar{e}=10^6$]{\includegraphics[scale=0.5]{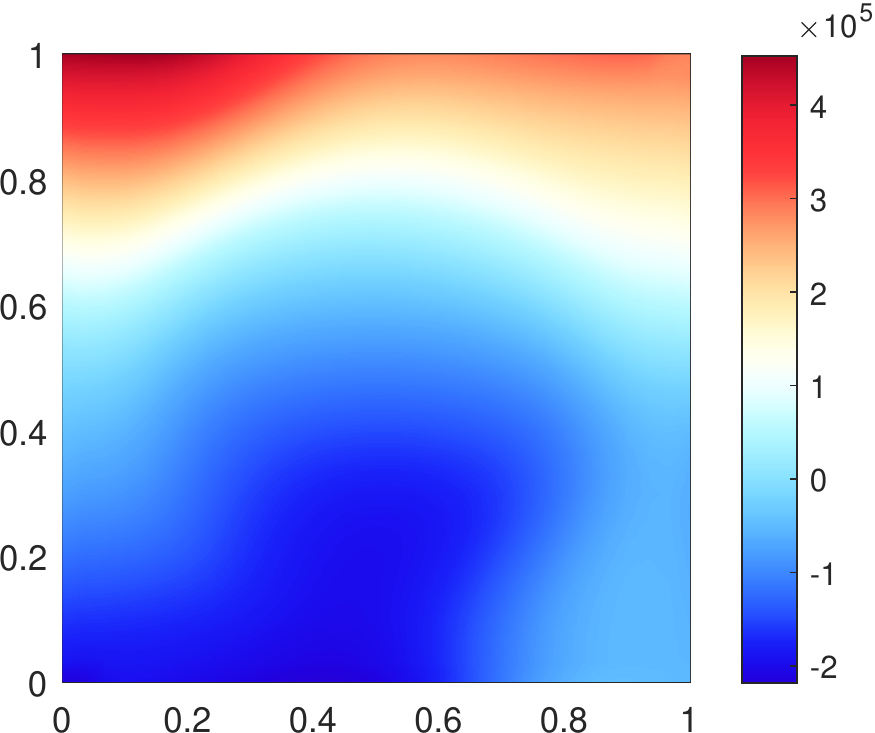}}
\caption{First-order RPC schemes, $P_2-P_1-P_2-P_2$, $\bar{e}=10^4, 10^5, 10^6$, $h$=1/30, $\tau=h^2$, $\T$=1.}
\label{Figure 3.2}
\end{figure}
\vspace{-8mm}
\begin{figure}[H]
\centering
\subfigure[angular velocity, $\bar{e}=10^4$]{\includegraphics[scale=0.5]{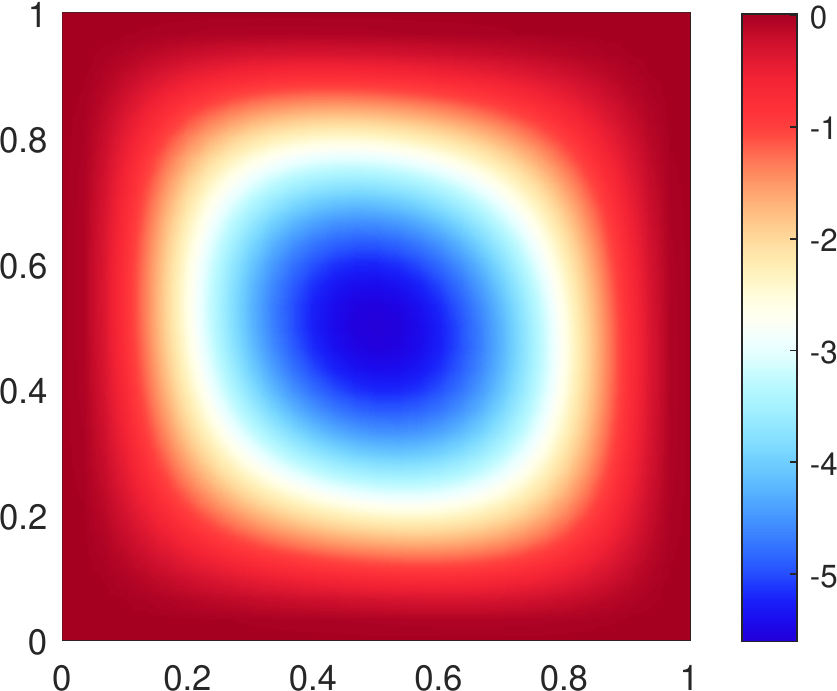}}
\quad
\subfigure[angular velocity, $\bar{e}=10^4$]{\includegraphics[scale=0.5]{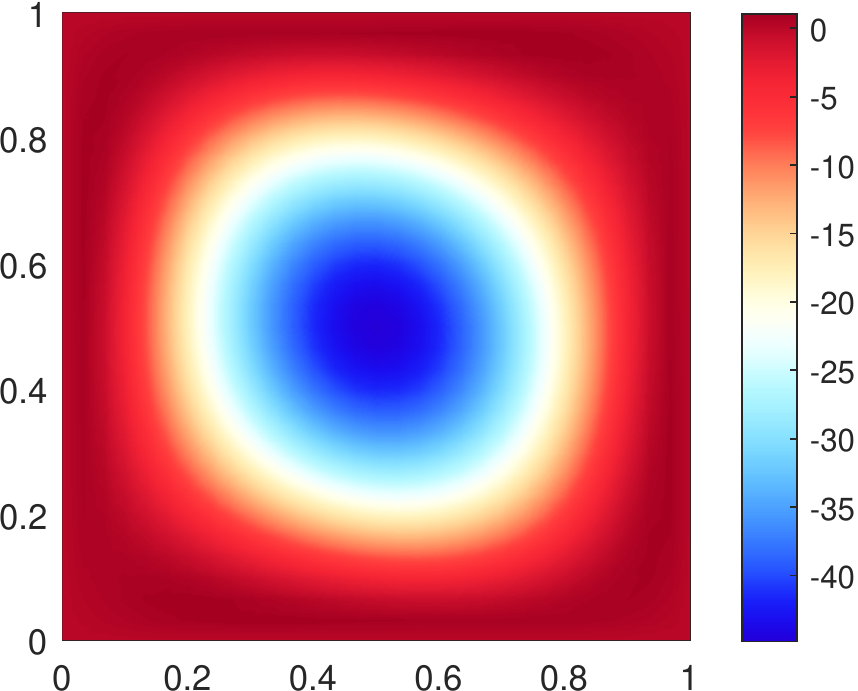}}
\quad
\subfigure[angular velocity, $\bar{e}=10^4$]{\includegraphics[scale=0.5]{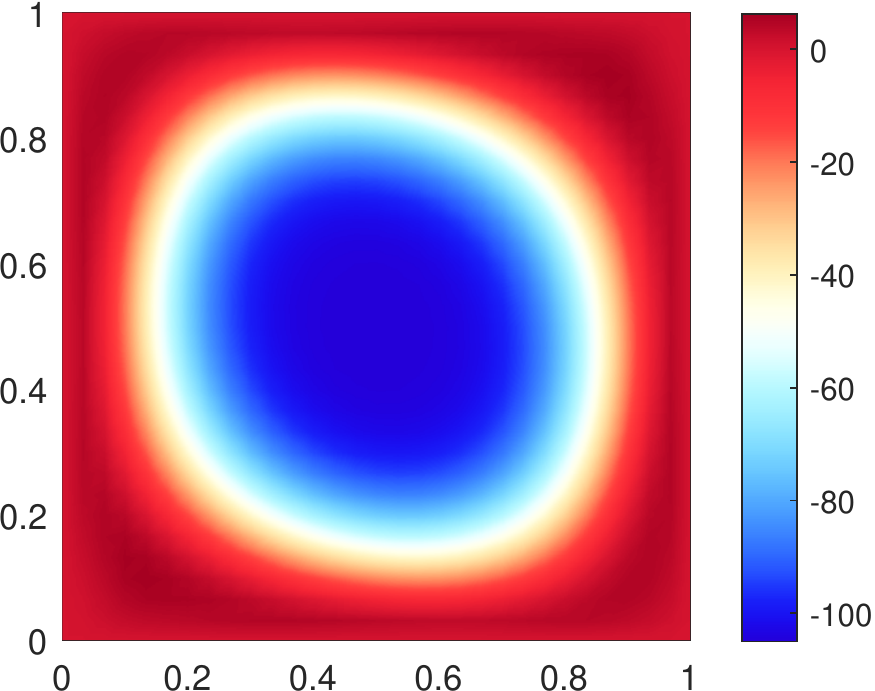}}
\caption{First-order RPC schemes, $P_2-P_1-P_2-P_2$, $\bar{e}=10^4, 10^5, 10^6$, $h$=1/30, $\tau=h^2$, $\T$=1.}
\label{Figure 3.3}
\end{figure}
\vspace{-8mm}
\begin{figure}[H]
\centering
\subfigure[velocity, $\bar{e}=10^4$]{\includegraphics[scale=0.5]{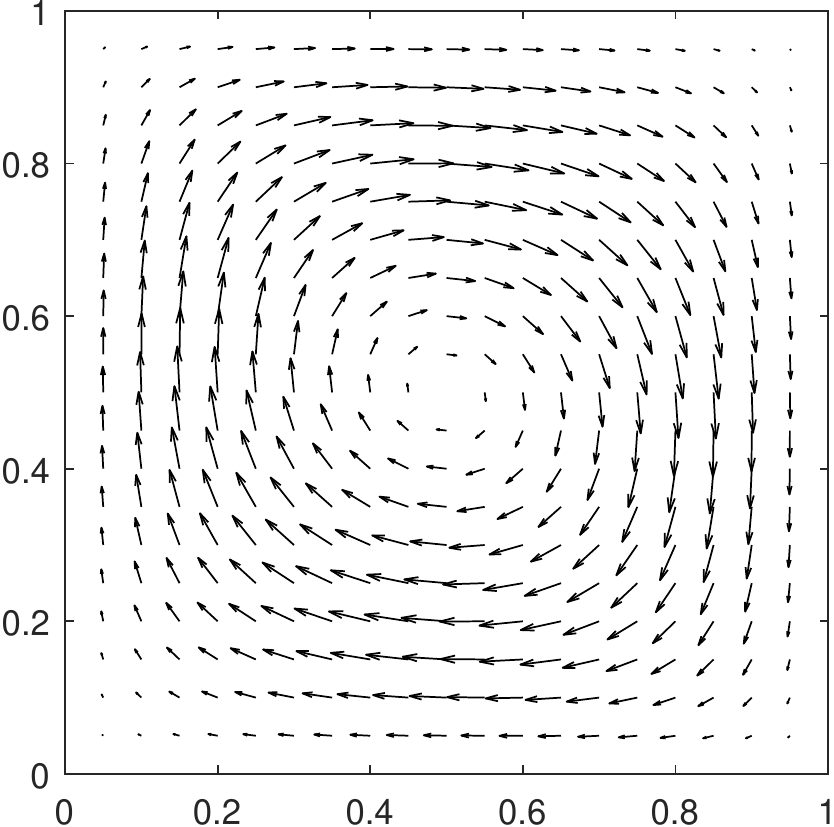}}
\quad
\subfigure[velocity, $\bar{e}=10^5$]{\includegraphics[scale=0.5]{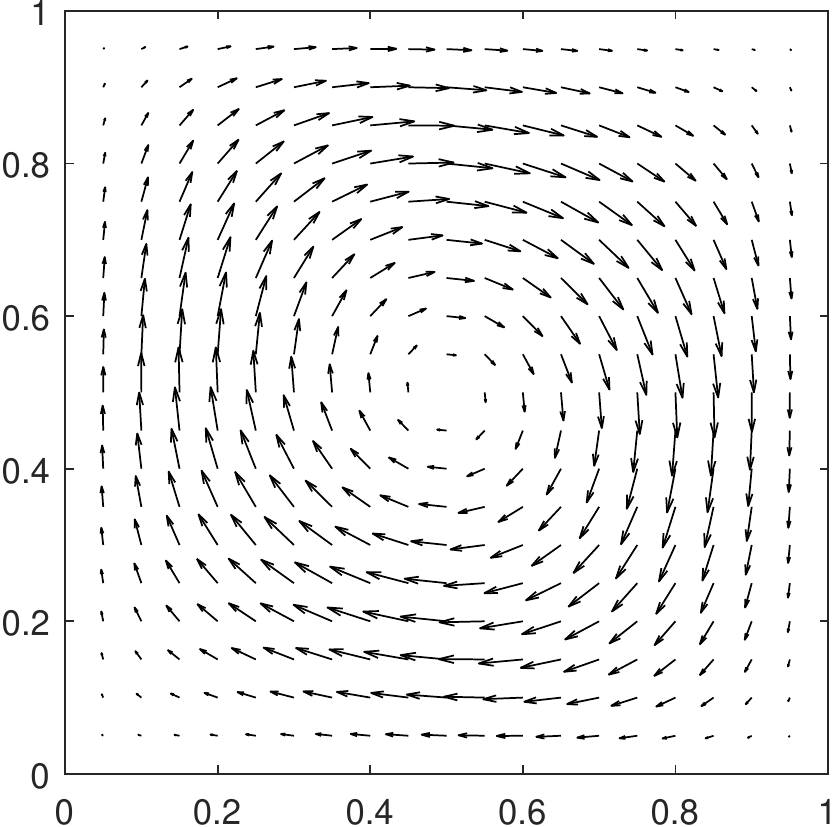}}
\quad
\subfigure[velocity, $\bar{e}=10^6$]{\includegraphics[scale=0.5]{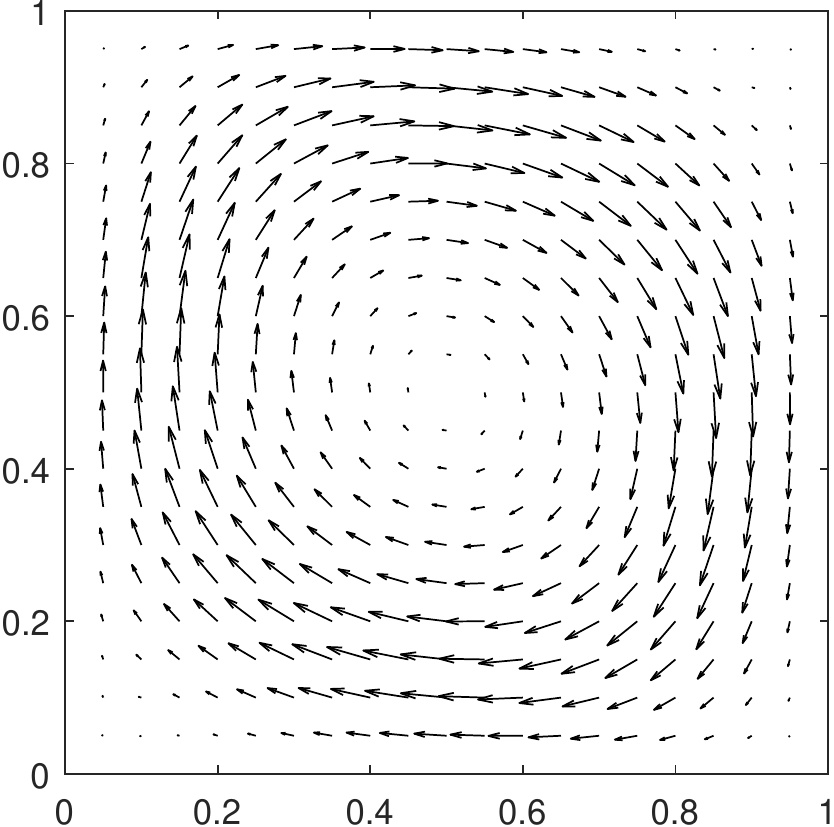}}
\caption{First-order RPC schemes, $P_2-P_1-P_2-P_2$, $\bar{e}=10^4, 10^5, 10^6$, $h$=1/30, $\tau=h^2$, $\T$=1.}
\label{Figure 3.4}
\end{figure}
\vspace{-8mm}
\begin{figure}[H]
\centering
\subfigure[vertical velocity at $x=0.5$]{\includegraphics[width=4.5cm]{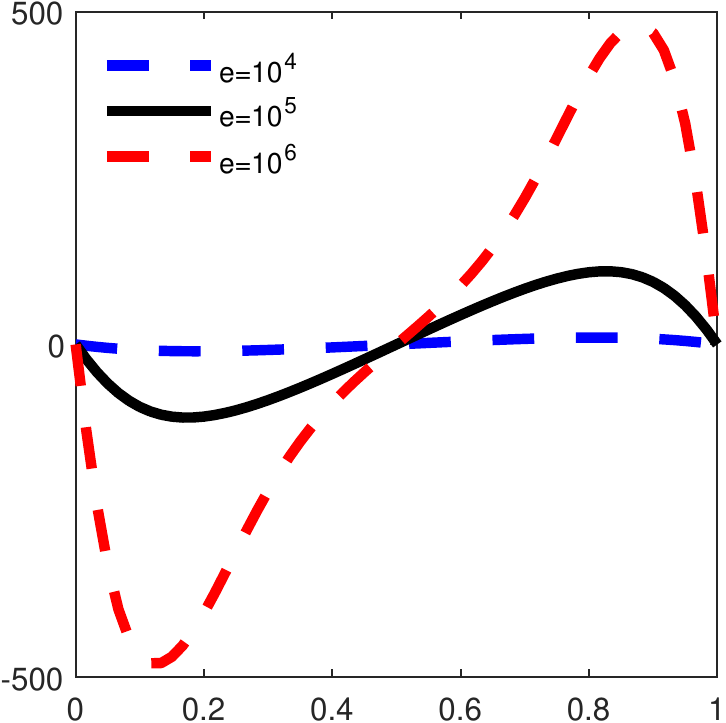}}
\quad
\subfigure[horizontal velocity at $y=0.5$]{\includegraphics[width=4.5cm]{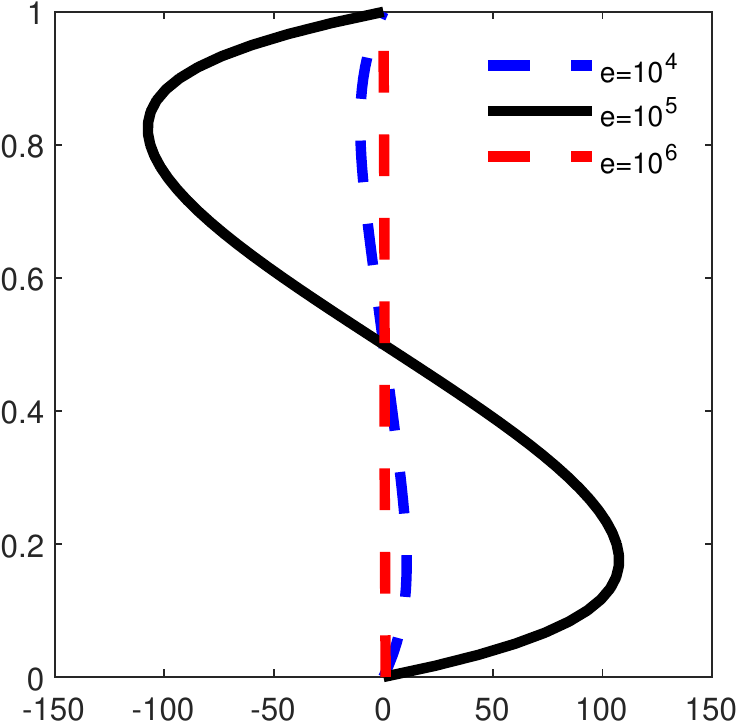}}
\quad
\subfigure[temperature at $y=0.5$]{\includegraphics[width=4.5cm]{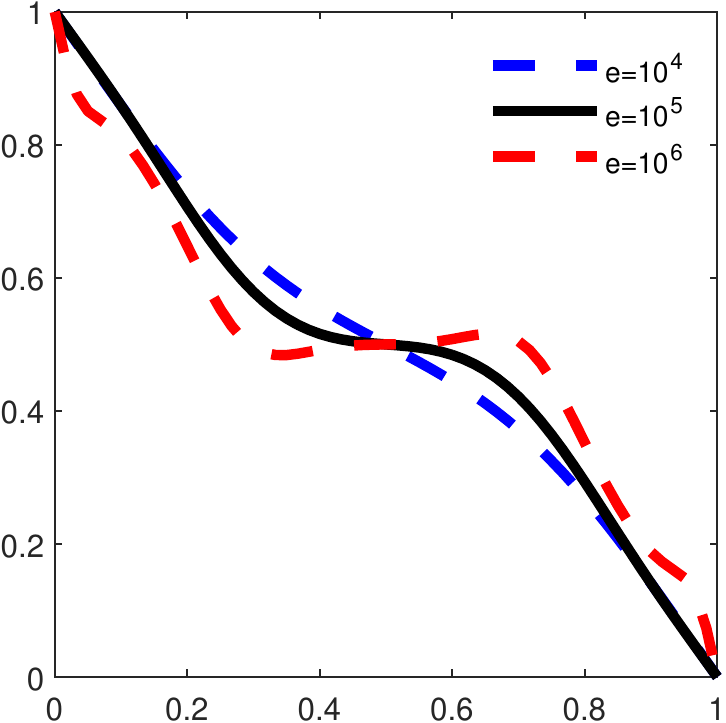}}
\caption{(a) vertical velocity at $x=0.5$; (b) horizontal velocity at $y=0.5$; (c) temperature at y =0.5; for varying thermal expansion numbers for first-order RPC schemes.}
\label{Figure 3.5}
\end{figure}
The Fig \ref{Figure 3.5} shows the vertical velocity distribution at $x=0.5$, the horizontal velocity distribution at $y=0.5$, and the temperature change at $x=0.5$. This is a very popular graphic illustration in buoyancy-driven cavity test research. We found that as the Rayleigh number increases, the differences in the profile diagrams become larger and larger, which is consistent with the previous study~\cite{2017Pressure}.

It can be easily seen from the isotherm of the Figs.~\ref{Figure 3.1}-\ref{Figure 3.4}. The hot fluid is transported to the cold wall. It can be seen from the isobaric and streamline diagrams that the pressure difference in the center of the cavity gradually increases, and the angular velocity and flow field gradually deform from a circle to an ellipse.
\subsection{Thermal driven cavity flow problem of 3D}
We consider the domain $\Omega=[0, 1]\times[0, 1]\times[0, 1]$ and given $f_1, f_2, f_3=0$, $\hat{e}=10^4, 10^5, 10^6$, $\nu=\nu_r=\a=\b=\k=D=1$, $h$=1/15, $\T=1$, $\tau=0.01$. The boundary conditions are given below, and the initial values are the same as the boundary.
\begin{eqnarray}
\begin{array}{l}
T(0,y,z,t)=1, \quad T(1,y,z,t)=0, \quad other \quad \frac{\partial T}{\partial n}=0,\\
u|_{\partial \Omega}=0, \quad \omega |_{\partial \Omega}=0, \quad f_1=0, \quad f_2=0, \quad f_3=0. \no
\end{array}
\end{eqnarray}
\begin{figure}[H]
   \centering
   \begin{minipage}[l]{0.3\textwidth}
   \centerline{\includegraphics[scale=0.6]{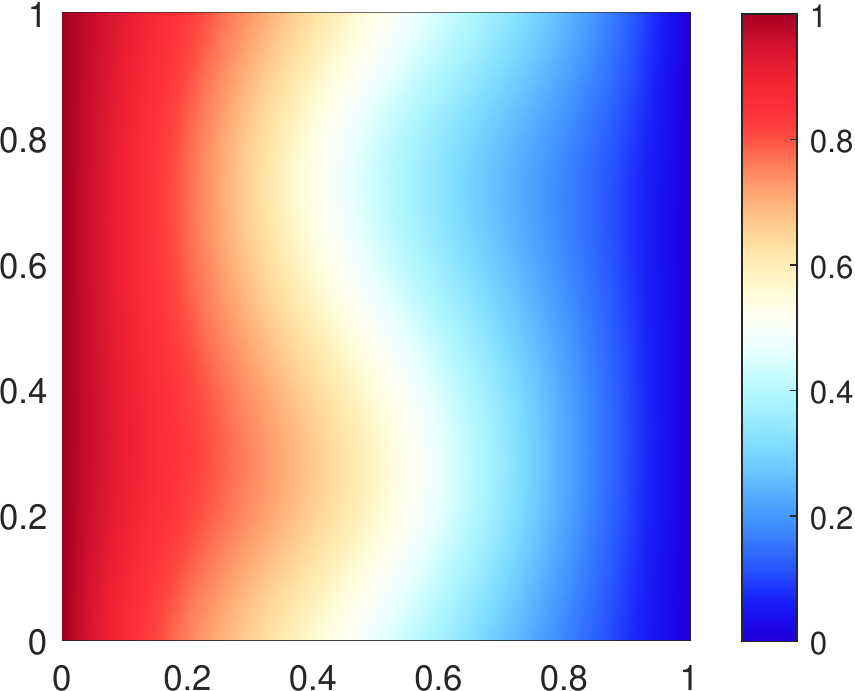}}
   \end{minipage}
   \hfill
   \begin{minipage}[c]{0.3\textwidth}
   \centering
   \centerline{\includegraphics[scale=0.6]{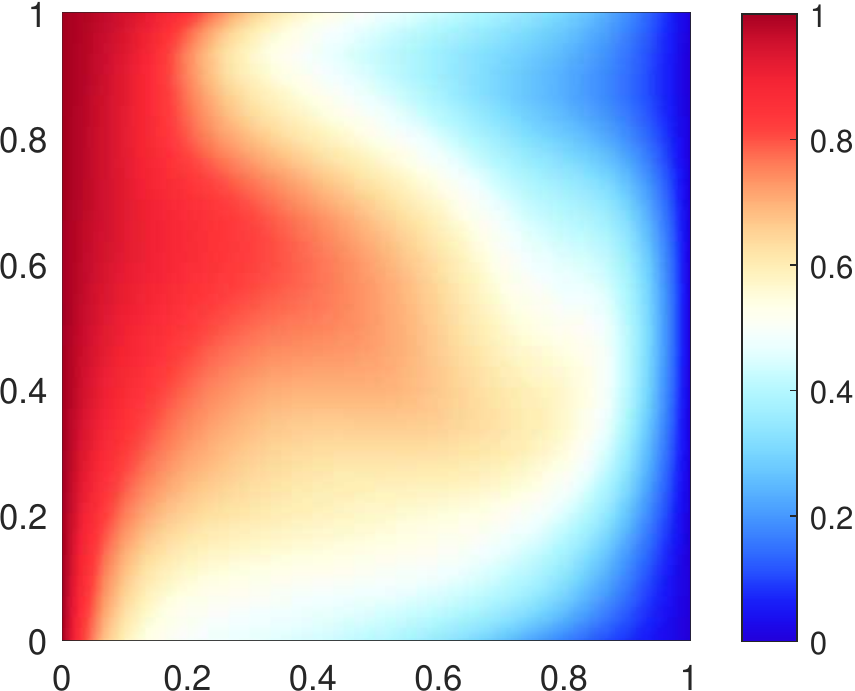}}
   \end{minipage}
    \hfill
   \begin{minipage}[r]{0.3\textwidth}
   \centering
   \centerline{\includegraphics[scale=0.6]{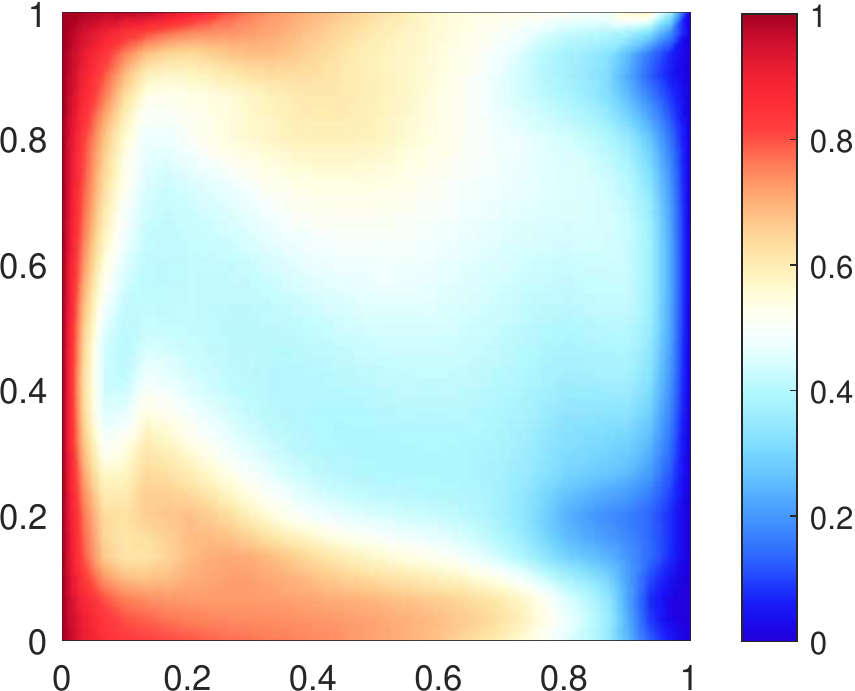}}
   \end{minipage}
   \caption{temperature, $P_2-P_1-P_2-P_2$, $\bar{e}=10^4, 10^5, 10^6$, $h=1/15$, $\tau=h^2$, $\T$=1.}
   \label{Figure 4.1}
\end{figure}
\begin{figure}[H]
   \centering
   \begin{minipage}[l]{0.3\textwidth}
   \centerline{\includegraphics[scale=0.6]{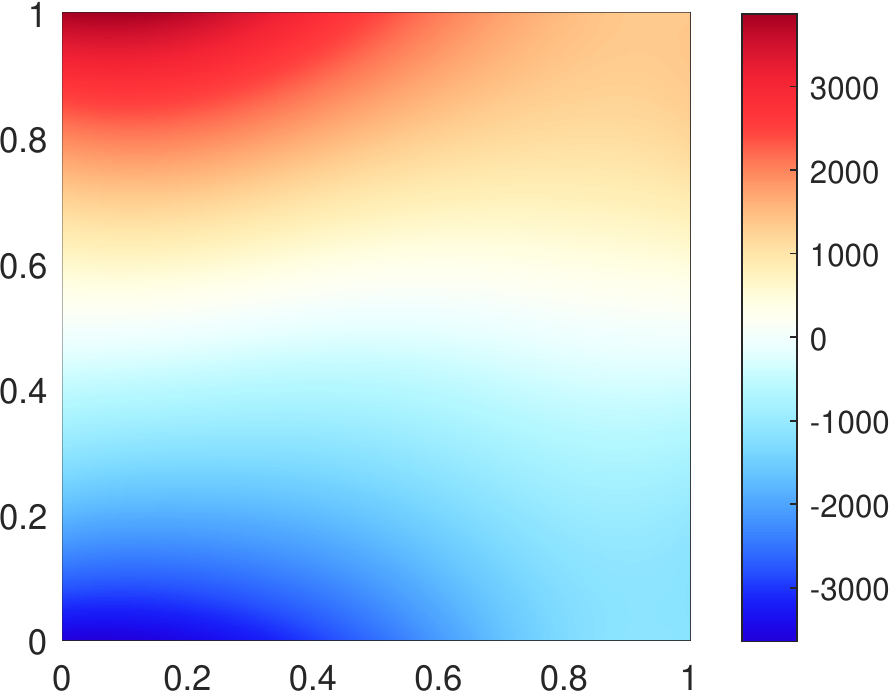}}
   \end{minipage}
   \hfill
   \begin{minipage}[c]{0.3\textwidth}
   \centering
   \centerline{\includegraphics[scale=0.6]{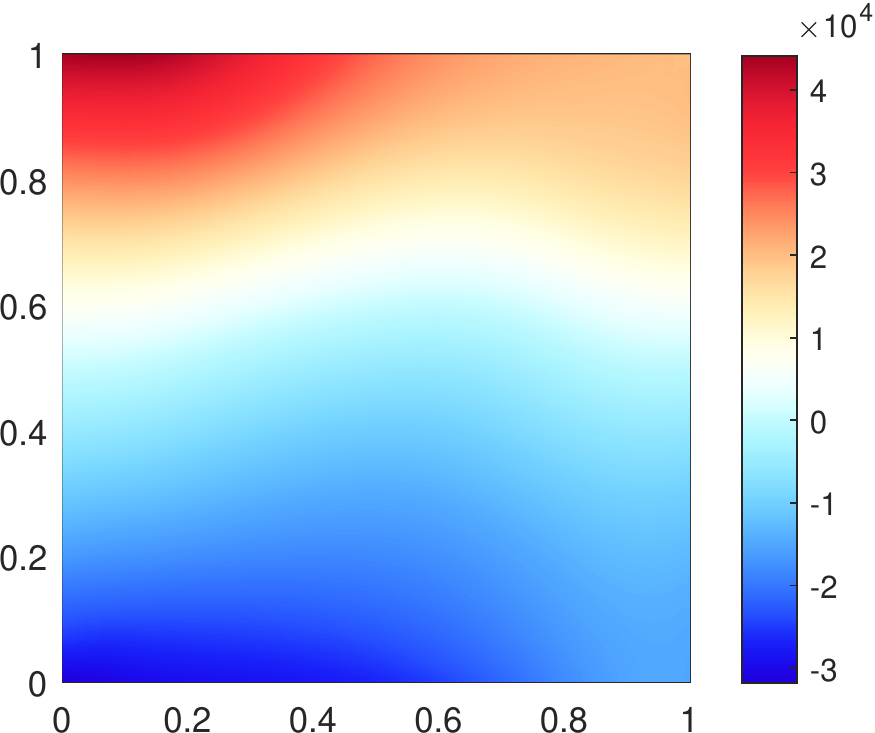}}
   \end{minipage}
    \hfill
   \begin{minipage}[r]{0.3\textwidth}
   \centering
   \centerline{\includegraphics[scale=0.6]{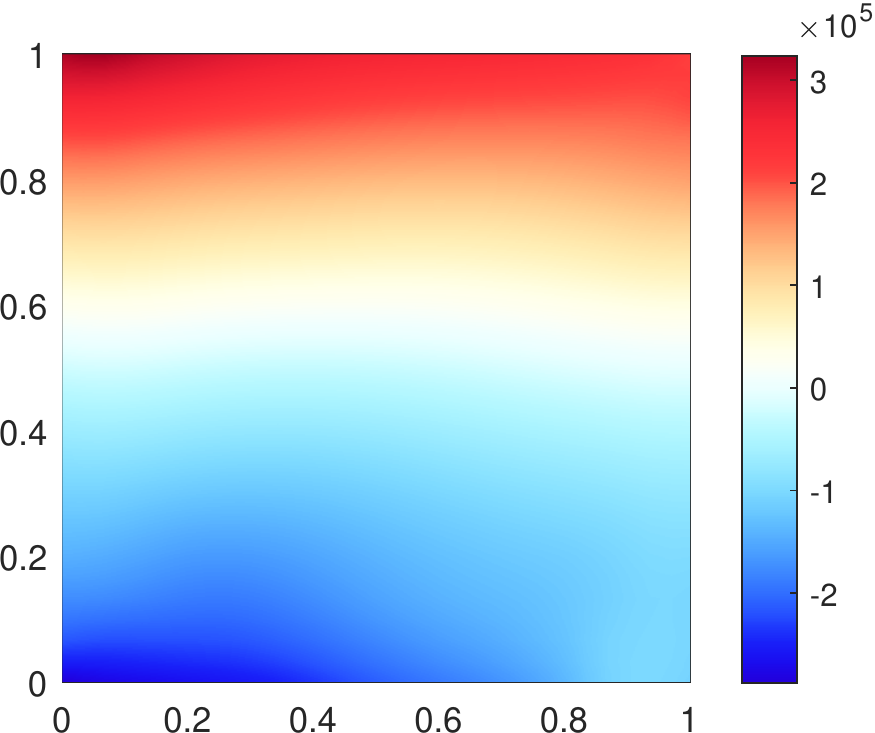}}
   \end{minipage}
   \caption{pressure, $P_2-P_1-P_2-P_2$, $\bar{e}=10^4, 10^5, 10^6$, $h=1/15$, $\tau=h^2$, $\T$=1.}
   \label{Figure 4.2}
\end{figure}
\begin{figure}[H]
   \centering
   \begin{minipage}[l]{0.3\textwidth}
   \centerline{\includegraphics[scale=0.6]{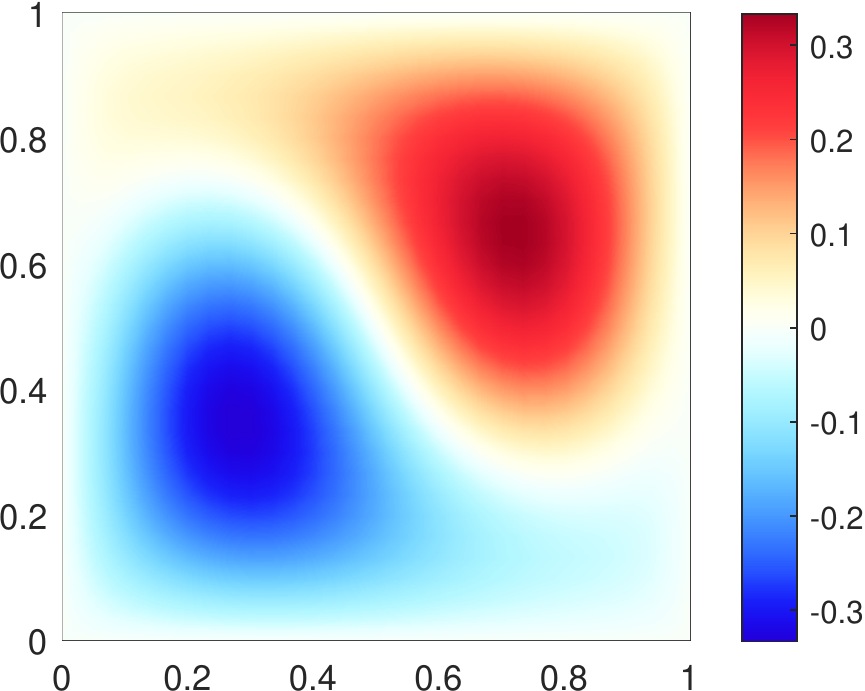}}
   \end{minipage}
   \hfill
   \begin{minipage}[c]{0.3\textwidth}
   \centering
   \centerline{\includegraphics[scale=0.6]{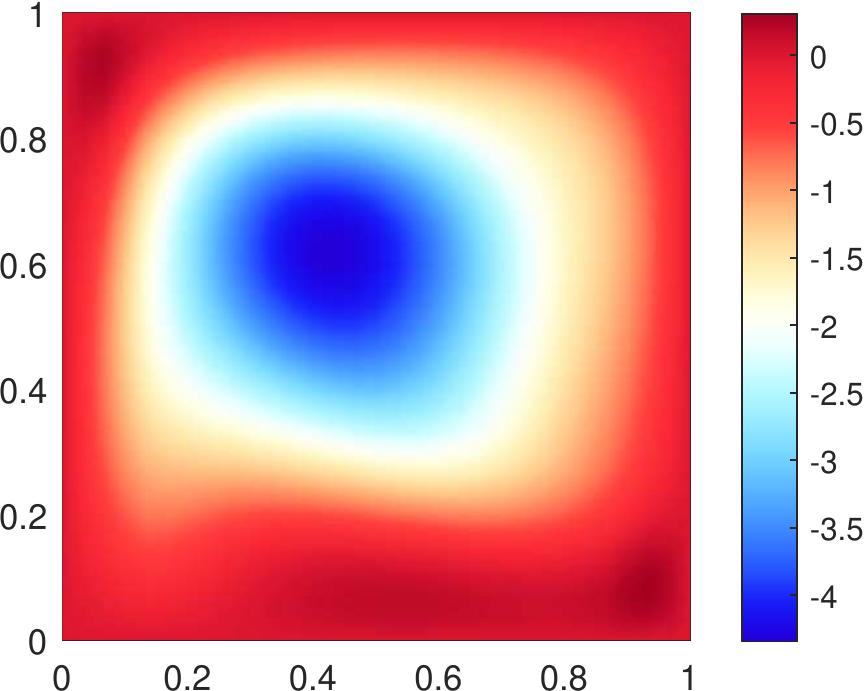}}
   \end{minipage}
    \hfill
   \begin{minipage}[r]{0.3\textwidth}
   \centering
   \centerline{\includegraphics[scale=0.6]{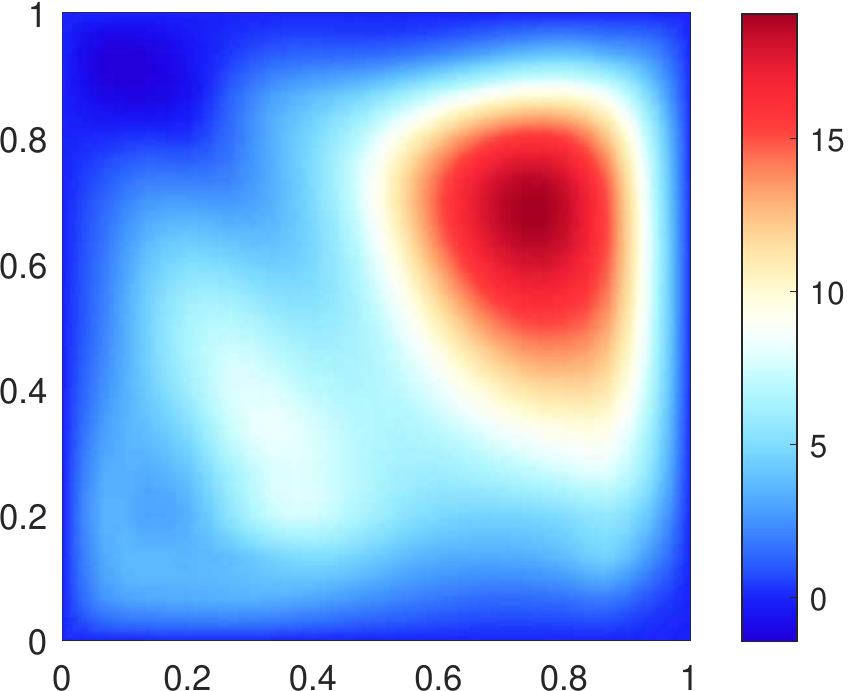}}
   \end{minipage}
   \caption{angular velocity, $P_2-P_1-P_2-P_2$, $\bar{e}=10^4, 10^5, 10^6$, $h=1/15$, $\tau=h^2$, $\T$=1.}
   \label{Figure 4.3}
\end{figure}
\begin{figure}[H]
   \centering
   \begin{minipage}[l]{0.3\textwidth}
   \centerline{\includegraphics[scale=0.35]{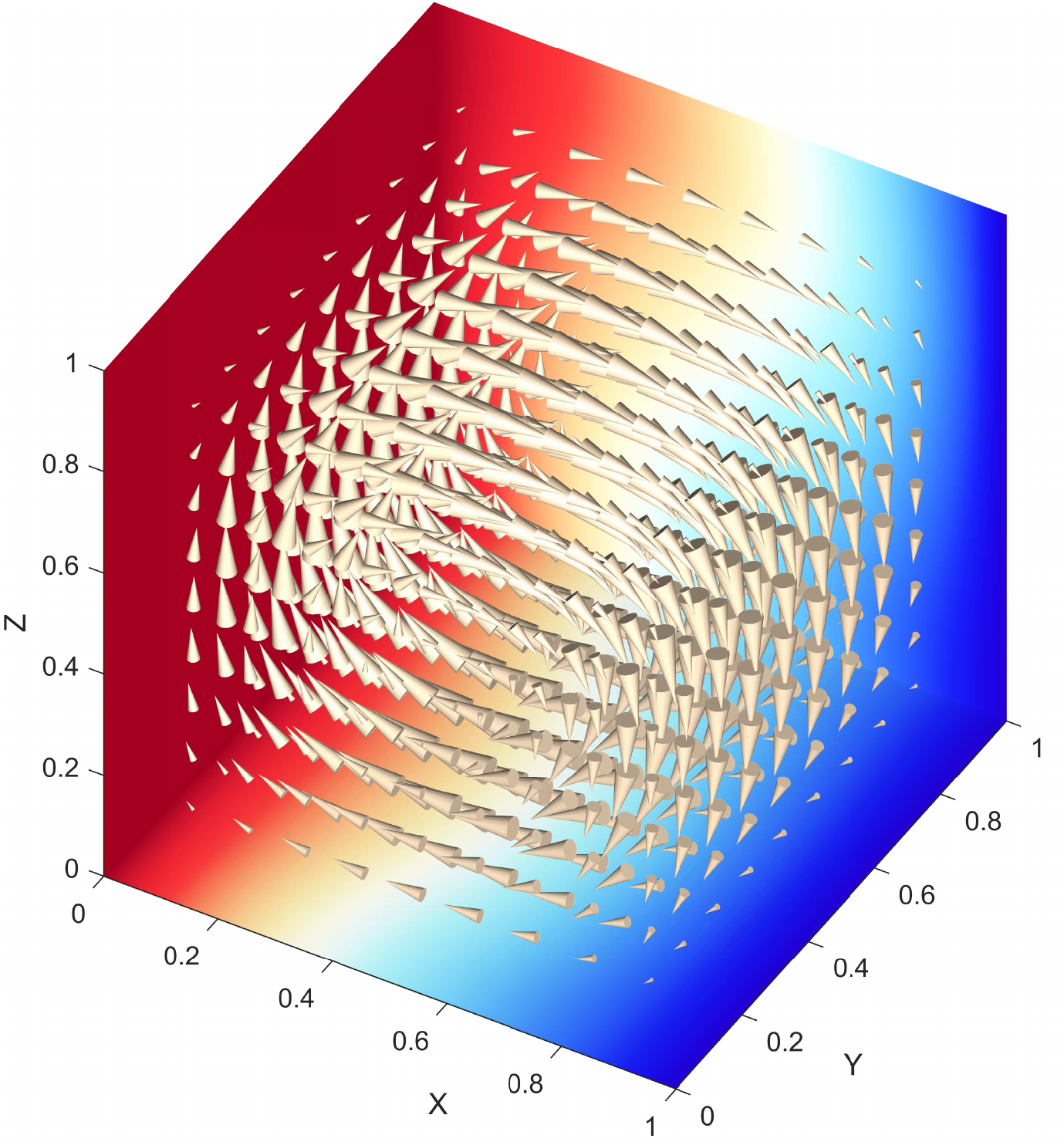}}
   \end{minipage}
   \hfill
   \begin{minipage}[c]{0.3\textwidth}
   \centering
   \centerline{\includegraphics[scale=0.35]{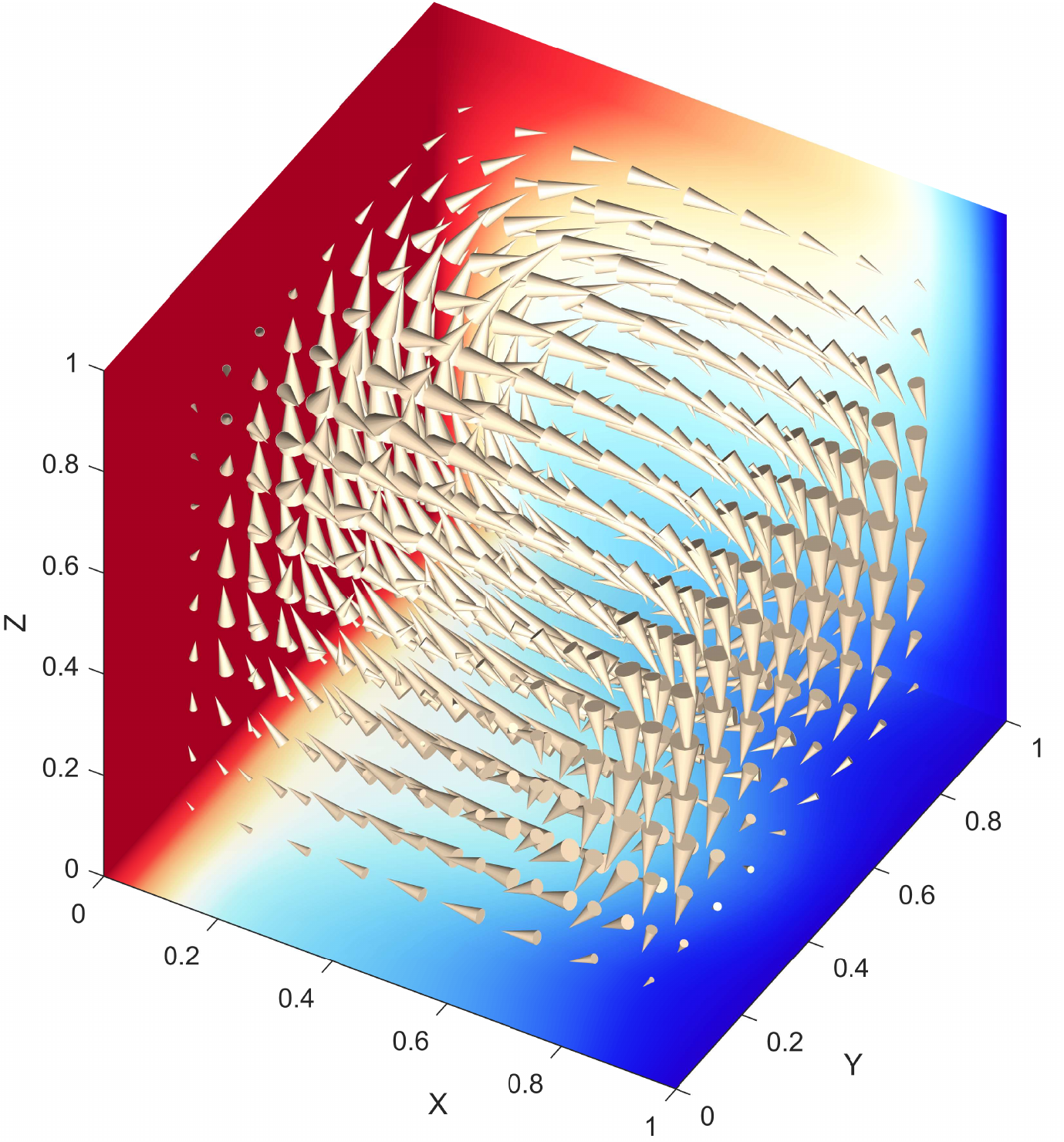}}
   \end{minipage}
    \hfill
   \begin{minipage}[r]{0.3\textwidth}
   \centering
   \centerline{\includegraphics[scale=0.35]{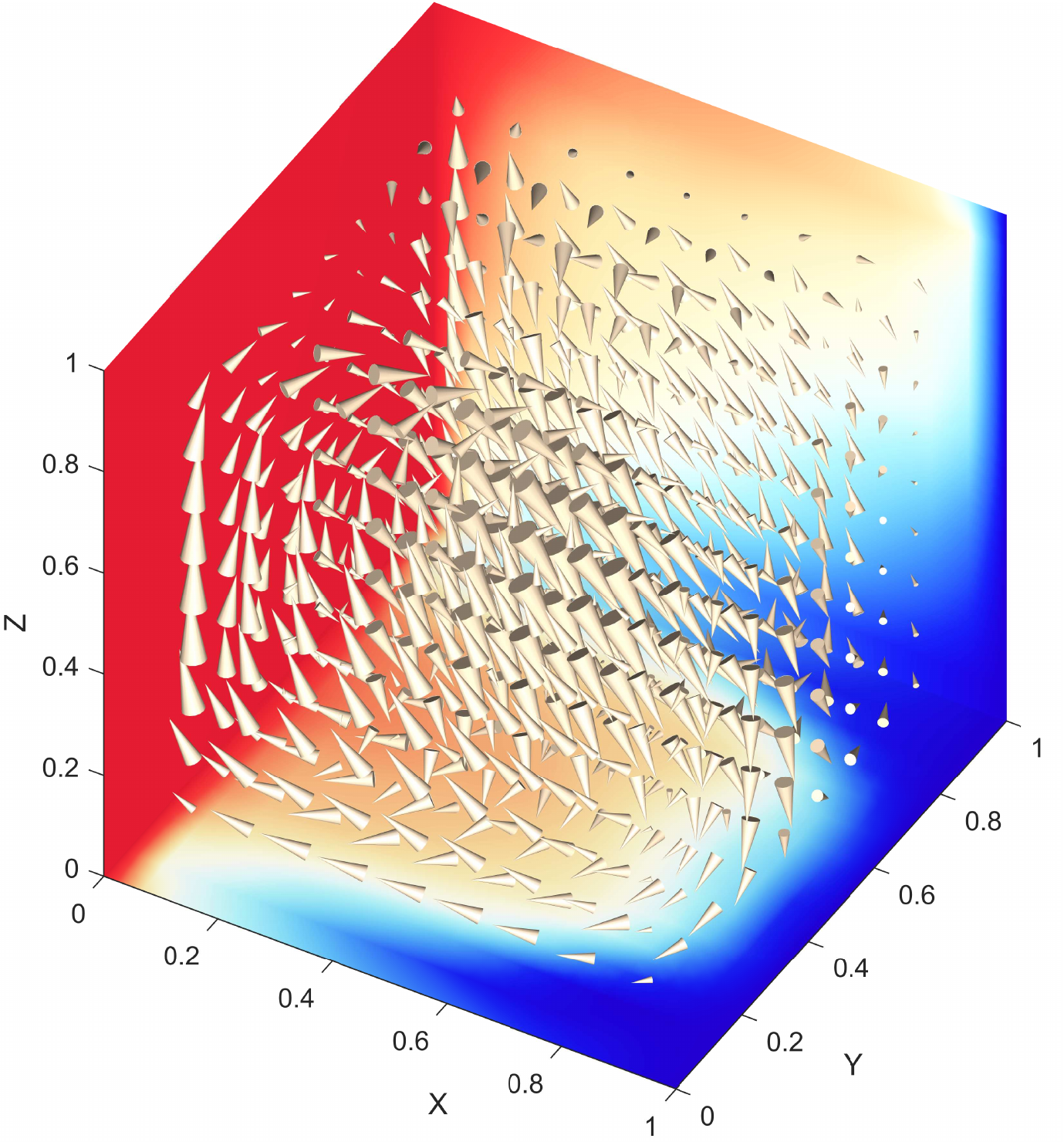}}
   \end{minipage}
   \caption{velocity field, $P_2-P_1-P_2-P_2$, $\bar{e}=10^4, 10^5, 10^6$, $h=1/15$, $\tau=h^2$, $\T$=1.}
   \label{Figure 4.4}
\end{figure}

We found that as the Rayleigh number increases, the differences in the profile diagrams become larger and larger, which can be easily seen from the isotherms in Figs.~\ref{Figure 4.1}-\ref{Figure 4.4}. The hot fluid is transported to the cold wall. It can be seen from the isobaric diagram and streamline diagram that the pressure difference in the center of the cavity gradually increases, which is consistent with the conclusion of the heat-driven cavity flow problem in 2D.
\subsection{Thermal driven cavity flow problem of 3D}
In this experiment, the problem of thermally driven cavity flow is considered. Choose $h=1/15$, $f_1, f_2, f_3=0$, $\hat{e}=10$, $\nu=\nu_r=\a=\b=\k=D=1$, $\tau=h^2$. Different from the research in other articles, the problem is transformed into a three-dimensional cube of $[0, 1]\times[0, 1]\times[0, 1]$. And draw the contour map of velocity, pressure, angular velocity, and temperature at $y=0.6$. The boundary conditions are given below, and the initial values are the same as the boundary.
\begin{eqnarray}
\begin{array}{l}
T(x,y,0,t)=4x(1-x), \quad \frac{\partial T(1,y,z,t)}{\partial n}=0,\quad \frac{\partial T(x,0,z,t)}{\partial n}=0,\quad \frac{\partial T(x,1,z,t)}{\partial n}=0,\quad other \quad T|_{\partial \Omega}=0, \\
u|_{\partial \Omega}=0 ,\quad \omega |_{\partial \Omega}=0 ,\quad f_1=0, \quad f_2=0, \quad f_3=0. \no
\end{array}
\end{eqnarray}
\begin{figure}[H]
\centering
\subfigure[$\T$=1, velocity field]{\includegraphics[height=3cm]{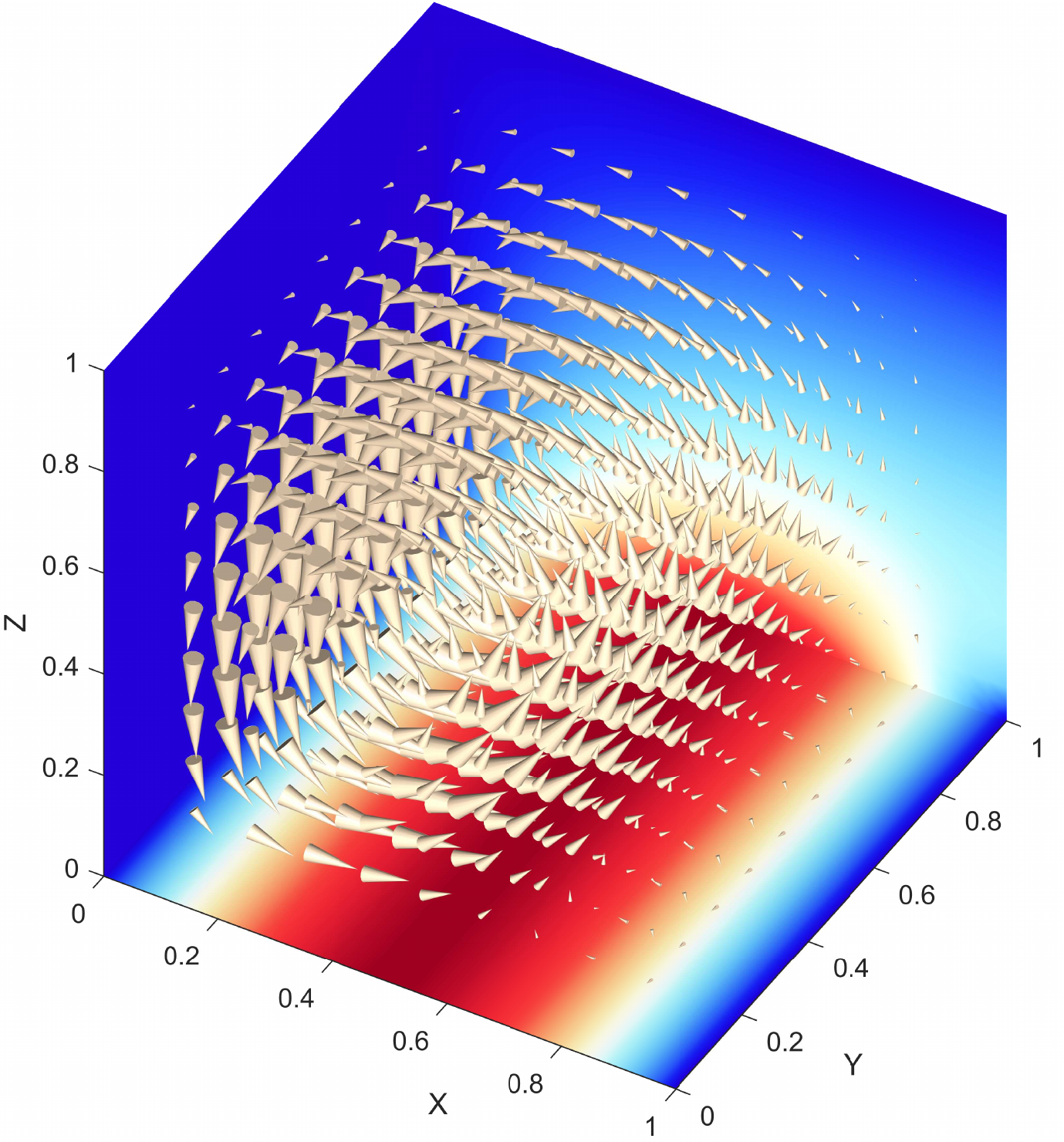}}
\quad
\subfigure[$\T$=1, temperature]{\includegraphics[height=3cm]{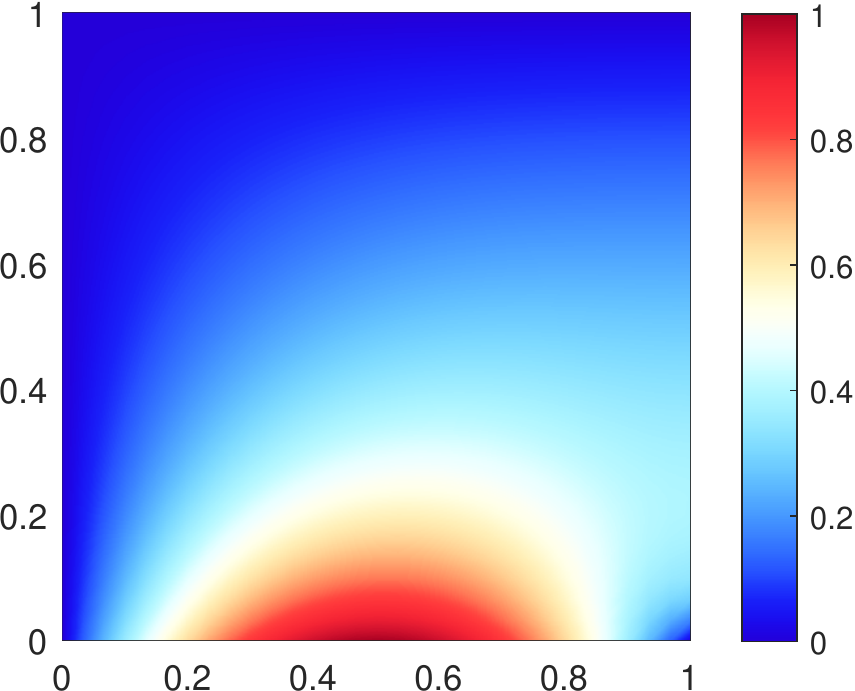}}
\quad
\subfigure[$\T$=1, pressure]{\includegraphics[height=3cm]{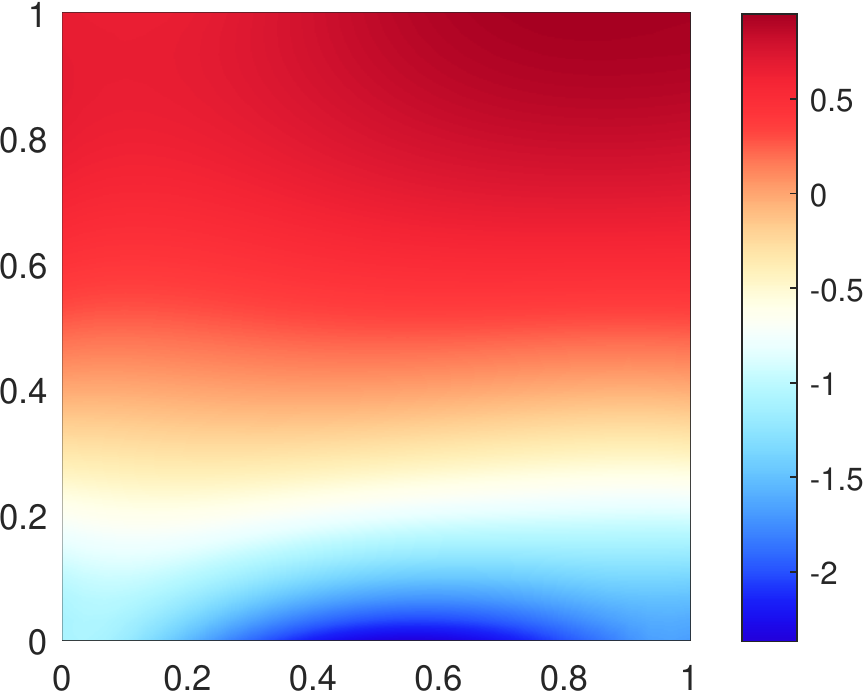}}
\quad
\subfigure[$\T$=1, angular velocity]{\includegraphics[height=3cm]{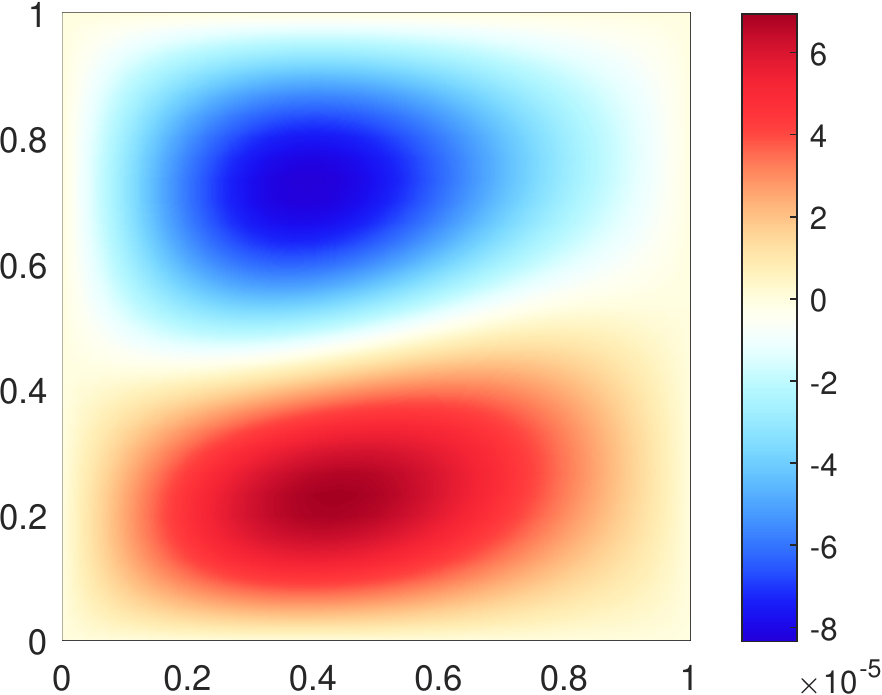}}
\caption{First-order RPC schemes, $P_2-P_1-P_2-P_2$, $h=1/15$, $\tau=h^2$, $\T$=1.}
\label{Figure 5.1}
\end{figure}
Due to the large amount of calculation for the 3D problem, it is noted that the use of the second-order RPC scheme can increase the time step, significantly reduce the time layer iteration, and greatly shorten the CPU time. In order to verify the validity of the second-order format at the same time, we compare the results obtained from the first-order format and the second-order format.
\begin{figure}[H]
\centering
\subfigure[$\T$=1, velocity field]{\includegraphics[height=3cm]{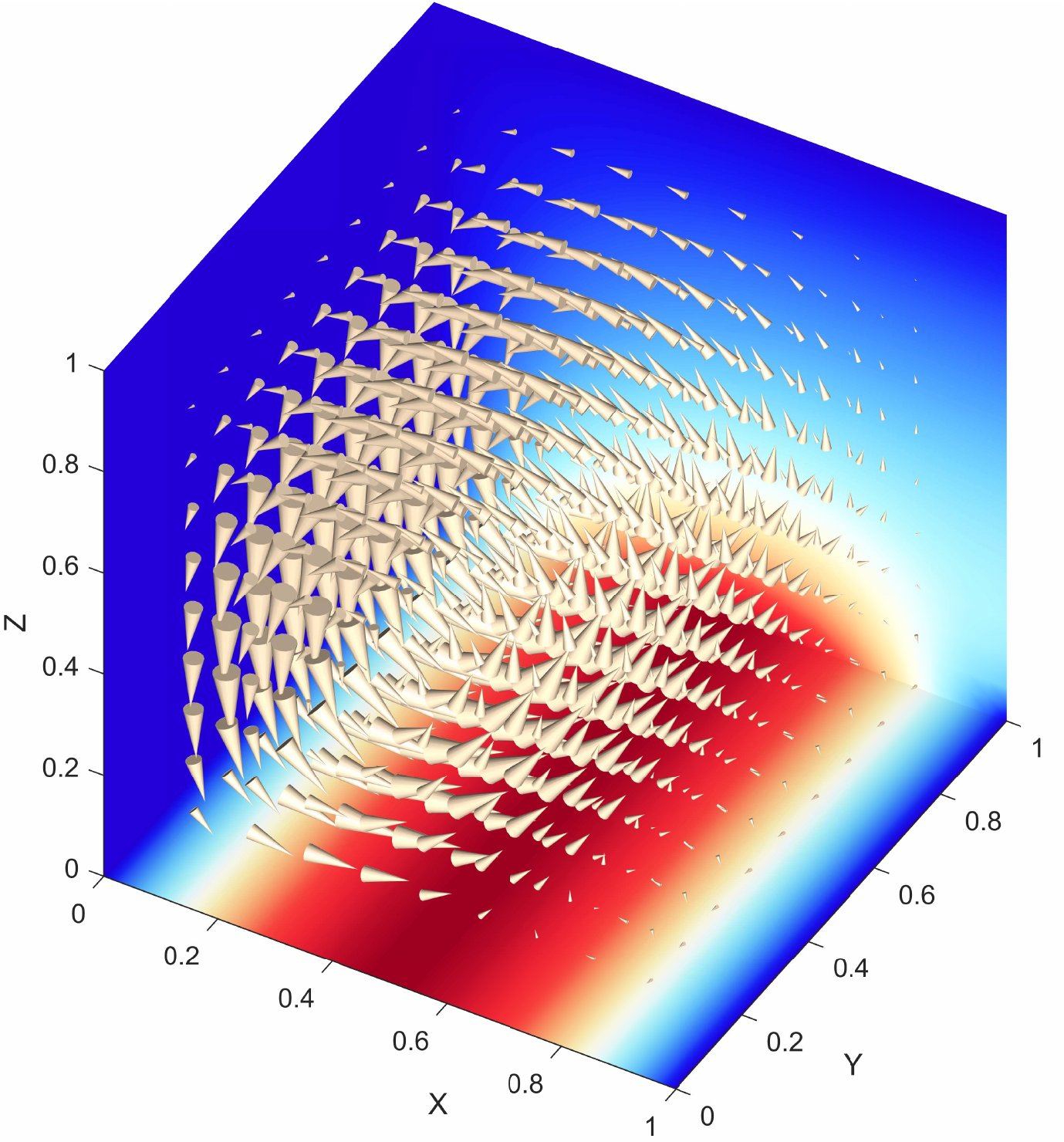}}
\quad
\subfigure[$\T$=1, temperature]{\includegraphics[height=3cm]{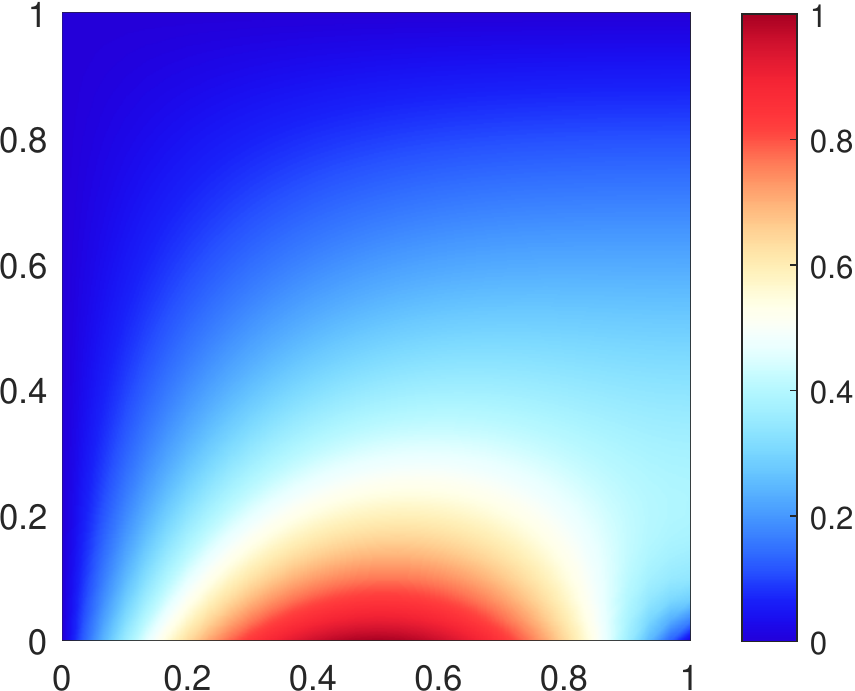}}
\quad
\subfigure[$\T$=1, pressure]{\includegraphics[height=3cm]{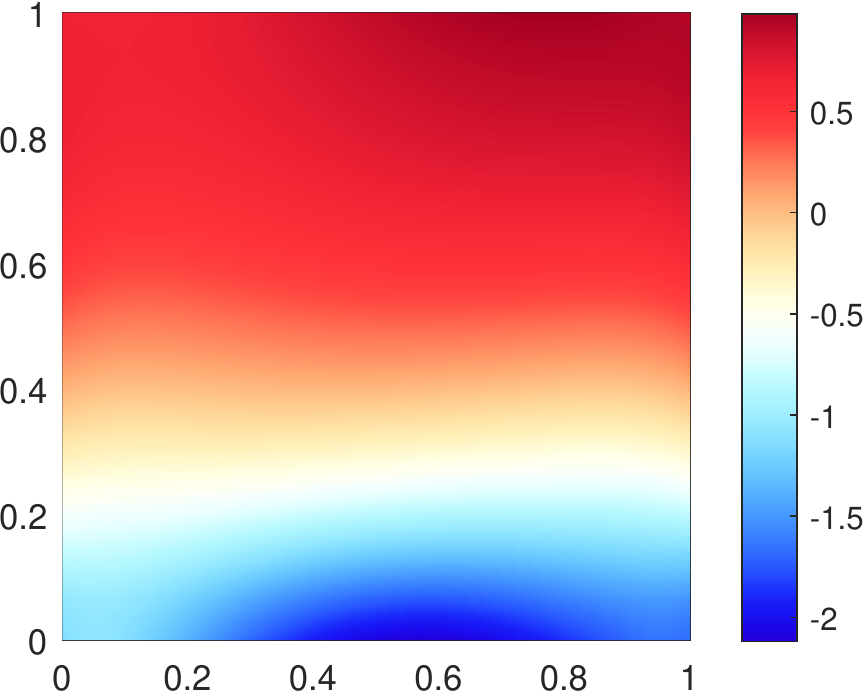}}
\quad
\subfigure[$\T$=1, angular velocity]{\includegraphics[height=3cm]{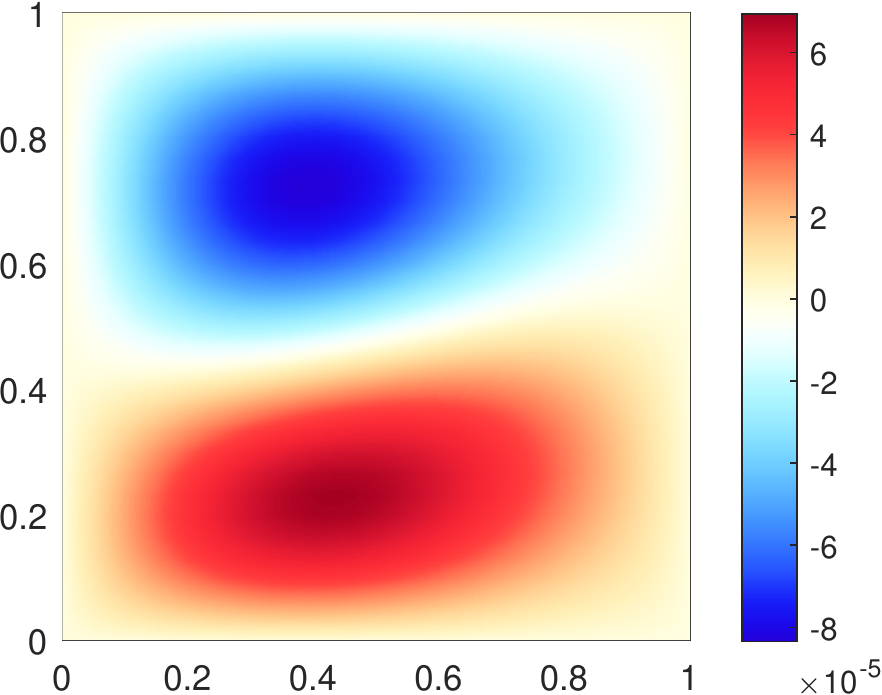}}
\caption{Second-order RPC schemes, $P_2-P_1-P_2-P_2$, $h=1/15$, $\tau=h$, $\T$=1.}
\label{Figure 5.2}
\end{figure}
It can be observed that when $\T$=1, the numerical results of the first-order format and the second-order format are almost the same, see Figs.\ref{Figure 5.1}-\ref{Figure 5.2}, but the first-order format requires 8152.5s, while the second-order format only needs 625.6s. The second-order format is appropriately selected, which will save a lot of CPU time.
\section{Conclusions}
In this paper, the first-order and second-order RPC methods for the problem of time-dependent thermal micropolar fluids are proposed. It also proves the stability and convergence of the first-order time semi-discrete scheme. In addition, a large number of numerical experiments have been carried out to assist in verifying the effectiveness of the RPC method in solving the actual thermomicropolar fluid model. At the same time, The projection method is still limited by the LBB condition, which will result in additional computational cost. This problem may be further optimized and discussed in future papers.
\\
{\bf Data Availability Statement}\\
Data sharing not applicable to this article as no datasets were generated or analyzed during the current study.

\end{document}